\documentclass[11pt]{article}
\usepackage{amsmath}
\usepackage{graphicx}
\usepackage{multirow}
\usepackage{enumerate}
\usepackage[authoryear]{natbib}
\usepackage{url}
\usepackage[margin=1in]{geometry}

\usepackage{amsmath,amssymb,dsfont,color,bm,mathtools,enumitem}
\mathtoolsset{showonlyrefs}
\usepackage{amsthm}
\usepackage{subcaption}

\newtheorem{thm}{Theorem}
\newtheorem{lem}{Lemma}
\newtheorem{prop}{Proposition}

\newtheorem{conj}{Conjecture}

\theoremstyle{definition}
\newtheorem{defn}{Definition}

\newtheorem{example}{Example}
\newtheorem{rmk}{Remark}

\usepackage[parfill]{parskip}
\usepackage{stmaryrd} 

\def\mbeta{\boldsymbol{\beta}}
\def\mbeta{\boldsymbol{\beta}}
\def\malpha{\boldsymbol{\alpha}}
\def\ma{\bm{a}}

\def\mc{\bm{c}}
\def\md{\bm{d}}
\def\me{\bm{e}}

\def\mj{\bm{j}}
\def\mk{\bm{k}}

\def\mp{\bm{p}}
\def\mq{\bm{q}}

\def\bmu{\bm{u}}
\def\mv{\bm{v}}

\def\mx{\bm{x}}
\def\my{\bm{y}}

\def\mA{\bm{A}}

\def\mP{\bm{P}}

\def\mU{\bm{U}}

\def\mX{\bm{X}}

\def\mx{\bm x}

\def\mA{\bm A}

\def\mX{\bm X}

\def\tZ{\mathcal{Z}}

\def\tB{\mathcal{B}}

\def\tE{\mathcal{E}}
\def\tF{\mathcal{F}}
\def\tG{\mathcal{G}}
\def\tH{\mathcal{H}}

\def\tM{\mathcal{M}}
\def\tN{\mathcal{N}}
\def\tO{\mathcal{O}}
\def\tP{\mathcal{P}}

\def\tS{\mathcal{S}}

\def\tY{\mathcal{Y}}
\def\tZ{\mathcal{Z}}

\newcommand{\normSize}[2]{#1\lVert#2#1\rVert_{F,I}}

\def\bbR{\mathbb{R}}
\def\bbE{\mathbb{E}}

\def\bbN{\mathbb{N}}

\newcommand{\mnorm}[1]{\left\lVert#1\right\rVert_\infty}

\newcommand{\vnormSize}[2]{#1\lVert#2#1\rVert_2}

\newcommand{\FnormSize}[2]{#1\lVert#2#1\rVert_F}

\DeclareMathOperator*{\argmin}{arg\,min}

\usepackage{mathrsfs}  
\def\caliB{\mathscr{B}}
\usepackage{wrapfig}
\newcommand*{\KeepStyleUnderBrace}[1]{
  \mathop{%
    \mathchoice
    {\underbrace{\displaystyle#1}}%
    {\underbrace{\textstyle#1}}%
    {\underbrace{\scriptstyle#1}}%
    {\underbrace{\scriptscriptstyle#1}}%
  }\limits
}

\allowdisplaybreaks
\usepackage[colorlinks,citecolor=blue]{hyperref}

\title{Statistical and Computational Efficiency for Smooth Tensor Estimation with Unknown Permutations}
\date{}
\author{%
Chanwoo Lee \\
University of Wisconsin -- Madison\\
\texttt{chanwoo.lee@wisc.edu} \\
\and
Miaoyan Wang\thanks{
    To whom correspondence should be addressed: miaoyan.wang@wisc.edu. The authors gratefully acknowledge NSF CAREER DMS-2141865, DMS-1915978, DMS-2023239, EF-2133740, and funding from the Wisconsin Alumni Research foundation.} \\
University of Wisconsin -- Madison\\
\texttt{miaoyan.wang@wisc.edu} \\
}

\begin{document}
\maketitle

\begin{abstract}%
We consider the problem of structured tensor denoising in the presence of unknown permutations. Such data problems arise commonly in recommendation systems, neuroimaging, community detection, and multiway comparison applications. Here, we develop a general family of smooth tensor models up to arbitrary index permutations; the model incorporates the popular tensor block models and Lipschitz hypergraphon models as special cases. We show that a constrained least-squares estimator in the block-wise polynomial family achieves the minimax error bound. A phase transition phenomenon is revealed with respect to the smoothness threshold needed for optimal recovery. In particular, we find that a polynomial of degree up to {\footnotesize $(m-2)(m+1)/2$} is sufficient for accurate recovery of order-$m$ tensors, whereas higher degrees exhibit no further benefits. This phenomenon reveals the intrinsic distinction for smooth tensor estimation problems with and without unknown permutations. Furthermore, we provide an efficient polynomial-time Borda count algorithm that provably achieves the optimal rate under monotonicity assumptions. The efficacy of our procedure is demonstrated through both simulations and Chicago crime data analysis. 
  \end{abstract}

\noindent%
{\it Keywords:} Tensor estimation, latent permutation, phase transition, statistical-computational efficiency, smooth tensor.
\vfill

\section{Introduction}
\label{sec:intro}
Higher-order tensor datasets are becoming ubiquitously in modern data science applications, for instance, recommendation systems \citep{bi2018multilayer}, social networks \citep{bickel2009nonparametric},
and genomics \citep{wang2019three}. Tensors provide effective representations of data structure that classical vector- and matrix-based methods fail to capture. For example, the music recommendation system~\citep{baltrunas2011incarmusic} records ratings of songs from users in various contexts. This three-way tensor of user $\times$ song $\times$ context allows us to investigate interactions between users and songs in a context-specific manner.

Tensor estimation problems cannot be solved without imposing structures. An appropriate reordering of tensor entries often provides an effective representation of the hidden structure. In the music recommendation example, suppose that we have certain criteria available (such as, similarities of music genres, ages of users, and importance of contexts) to reorder the songs, users, and contexts. Then, the sorted tensor will exhibit smooth structure, because entries from similar groups tend to have similar values.

\begin{figure}[ht]
    \centering
    \includegraphics[width = .95\textwidth]{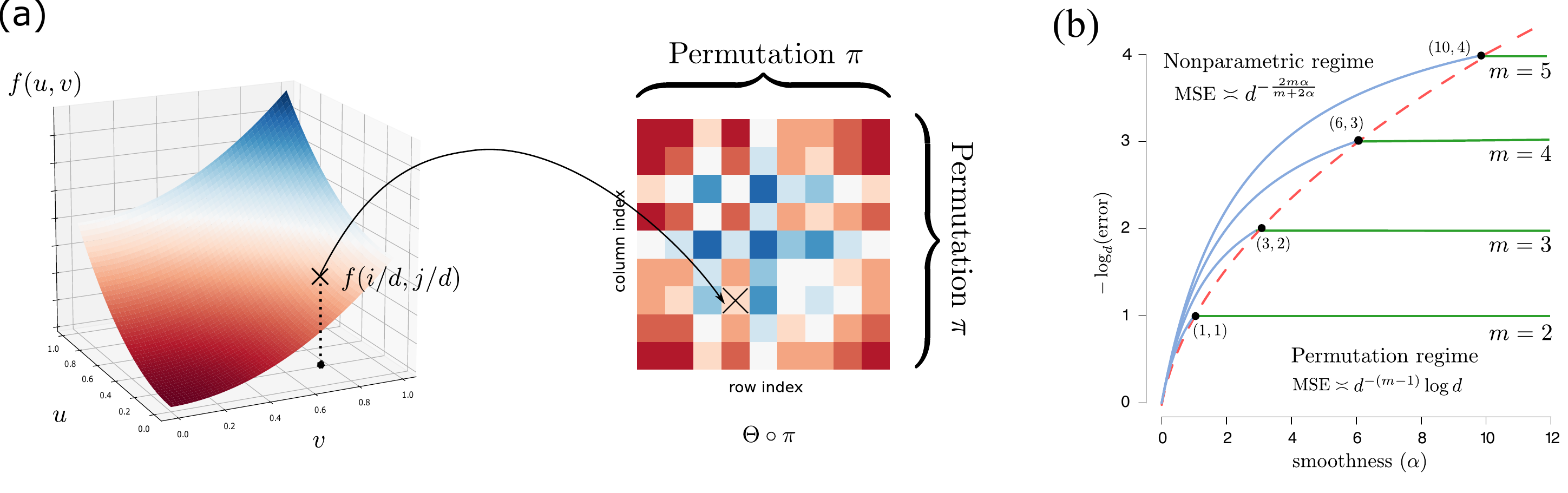}
    \caption{(a): Illustration of order-$m$ $d$-dimensional permuted smooth tensor models with $m=2$. (b): Phase transition of the mean squared error (MSE) (on the $-\log_d$ scale) as a function of smoothness $\alpha$ and tensor order $m$. Bold dots correspond to the critical smoothness level above which higher smoothness exhibits no further benefits for tensor estimation. }
    \label{fig:rate}
\end{figure}

In this article, we develop a \emph{permuted} smooth tensor model based on the aforementioned  motivation. We study a class of structured tensors, called \emph{permuted smooth tensor model}, of the following form:
\begin{align}\label{eq:gmd}
    \tY = \Theta\circ \pi + \text{noise}, \quad \text{where}\quad \Theta(i_1,\ldots,i_m)=f\left({i_1\over d},\ldots,{i_m\over d}\right),
\end{align}
where $\pi\colon[d]\rightarrow[d]$ is an \emph{unknown} latent permutation, $\Theta$ is an \emph{unknown} order-$m$ $d$-dimensional signal tensor, and $f$ is an \emph{unknown} multivariate function with a certain notion of smoothness, and $\Theta \circ \pi $ denotes the permuted tensor after reordering the indices along each of the $m$ modes. Figure~\ref{fig:rate}(a) shows an example of this generative model for the matrix case $m=2$. Our primary goal is to estimate the permuted smooth signal tensor $\Theta \circ \pi$ from the noisy tensor observation $\tY$ of arbitrary order $m$.

\subsection{Our contributions}
We develop a suite of statistical theory, efficient algorithms, and related applications for permuted smooth tensor models~\eqref{eq:gmd}.  
Our contributions are summarized below. 

 \begin{table*}[h]
    \centering
    \resizebox{\textwidth}{!}{%
    \begin{tabular}{c|@{\hskip4pt}c@{\hskip5pt}c@{\hskip4pt}c@{\hskip5pt}c@{\hskip4pt}c@{\hskip5pt}c@{\hskip4pt}c@{\hskip5pt}c@{\hskip5pt}c}
    & \citet{pananjady2022isotonic}&  \citet{balasubramanian2021nonparametric}&  \citet{li2019nearest}&\textbf{Ours} (MLE) & \textbf{Ours} (Borda count)\\
    \hline
       Model structure& monotonic & Lipschitz & Lipschitz &  $\alpha$-smooth  & $\alpha$-smooth \& monotonic\\
             Fixed grid design & $\surd$ &$\times$ & $\times$ & $\times$& $\surd$\\
     Error rate for order-$m$ tensor & \multirow{2}{*}{$d^{-1}$} & $d^{-{2m\over m+2}}$ & $d^{-\lfloor m/3\rfloor }$ & $d^{-{2m\alpha\over m+2\alpha}}\vee d^{-(m-1)}$ & $d^{-{2m\alpha\over m+2\alpha}}\vee d^{-(m-1)}$ \\
       (e.g., when $(m,\alpha)=(3,1)$) & & ($d^{-6/5}$) & $(d^{-1})$ & $(d^{-6/5})$ & $(d^{-6/5})$\\
          Minimax optimality& $\surd$  & $\times$ & $\times$ & $\surd$ & $-*$ \\
     Polynomial algorithm& $\surd$ &$\times$ & $\surd$ & $\times$& $\surd$ 
    \end{tabular}
    }
    \caption{Comparison of our results with previous work. For simplicity, we omit the log term in the rate. $^{*}$The optimality is achieved under extra Lipchitz monotonicity conditions. }
 \label{tab:comp}
\end{table*}

First, we develop a general permuted $\alpha$-smooth tensor model of arbitrary smoothness level $\alpha>0$. We establish the statistically optimal error rate and its dependence on model complexity. Specifically, we express the optimal rate as a function of tensor order $m$, tensor dimension $d$, and smoothness level $\alpha$, given by
\begin{equation}\label{eq:rate}
\text{Rate}(d):=d^{-{2m\alpha \over m+2\alpha}} \vee d^{-(m-1)}\log d.
\end{equation}

Table~\ref{tab:comp} summarizes the comparison of our work with previous results.  Our framework substantially generalizes earlier works that focus on only matrices with $m=2$~\citep{gao2015rate,klopp2017oracle} or Lipschitz function with $\alpha=1$~\citep{balasubramanian2021nonparametric,li2019nearest}. The generalization enables us to obtain results previously impossible: i) As tensor order $m$ increases, we demonstrate the failure of previous clustering-based algorithms~\citep{balasubramanian2021nonparametric,gao2015rate}, and we develop a new block-wise polynomial algorithm for tensors of order $m\geq 3$; ii) As smoothness $\alpha$ increases, we demonstrate that the error rate converges to a fast rate of $\tO(d^{-(m-1)})$, thereby disproving the conjectured lower bound $\tO(d^{-2m/(m+2)})$ posed by earlier work~\citep{balasubramanian2021nonparametric}. The results showcase the accuracy gain of our new approach, as well as the intrinsic distinction between matrices and higher-order tensors. 

Second, we discover a phase transition phenomenon with respect to the smoothness needed for optimal recovery in model~\eqref{eq:gmd}. Figure~\ref{fig:rate}(b) plots the dependence of estimation error in terms of smoothness level $\alpha$ for tensors of order $m$. We characterize two distinct error behaviors determined by a critical smoothness threshold; see Theorems~\ref{thm:LSE}-\ref{thm:minimax} in Section~\ref{sec:lse}. Specifically, the accuracy improves with $\alpha$ in the regime $\alpha\leq m(m-1)/2$, whereas the accuracy becomes a constant of $\alpha$ in the regime $\alpha>m(m-1)/2$. The results imply a polynomial of degree $(m-2)(m+1)/2$ (derived from $m(m-1)/2-1$) is sufficient for accurate recovery of order-$m$ tensors of arbitrary smoothness in model~\eqref{eq:gmd}, whereas higher degrees bring no further benefits. The phenomenon is distinct from matrix problems~\citep{klopp2017oracle,gao2015rate} and classical \emph{non-permuted} smooth function estimation~\citep{tsybakov2009introduction}, thereby highlighting the fundamental challenges in our new setting. These statistical contributions, to our best knowledge, are new to the literature of permuted smooth tensor problems. 

Third, we propose two estimation algorithms with accuracy guarantees: the least-squares estimation and Borda count estimation. The least-squares estimation, although being computationally hard, reveals the fundamental model complexity in the problem. The result serves as the benchmark and a useful guide to the algorithm design. Furthermore, we develop an efficient polynomial-time Borda count algorithm that provably achieves a minimax optimal rate under an extra Lipschitz monotonicity assumption. The algorithm handles a broad range of data types, including continuous and binary observations. 

Lastly, we illustrate the efficacy of our method through both simulations and data applications. A range of practical settings are investigated in simulations, and we show the outperformance of our method compared to alternative approaches. Application to Chicago crime data is presented to showcase the its usefulness. We identify the key global pattern and pinpoint local smooth structure in the denoised tensor. Our method will help practitioners efficiently analyze tensor datasets in various areas. Toward this end, the package and all data used are released at CRAN.

\subsection{Related work}\label{sec:priorwork}
 Our work is closely related to but also clearly distinctive from several lines of existing research. We review related literature for comparison.
\paragraph{Structure learning with latent permutations.} The estimation problem of~\eqref{eq:gmd} falls into the general category of structured learning with \emph{latent permutation}. Models involving latent permutations have recently received a surge of interest, including graphons~\citep{chan2014consistent,klopp2017oracle}, stochastic transitivity models~\citep{chatterjee2015matrix,shah2019low}, statistical seriation~\citep{flammarion2019optimal,hutter2020estimation}, and graph matching~\citep{ding2021efficient,livi2013graph}. These methods, however, are developed for matrices; the tensor counterparts are understood much less well. Table~\ref{tab:comp} summarizes the most related works to ours. Pananjady and Samworth~\citep{pananjady2022isotonic} studied the permuted tensor estimation under isotonic constraints.  We find that our smooth model achieves a much faster rate $\tO(d^{-(m-1)})$ than the rate $\tO(d^{-1})$ for isotonic models. The works~\citep{balasubramanian2021nonparametric,li2019nearest} studied similar smooth models as ours, but we gain substantial improvement in both statistics and computations. \citet{balasubramanian2021nonparametric} developed a (non-polynomial-time) clustering-based algorithm with a rate $\tO(d^{-2m/(m+2)})$. \citet{li2019nearest} developed a (polynomial-time) nearest neighbor estimation with a rate $\tO(d^{-\lfloor m/3\rfloor})$. Neither approach investigates the minimax optimality. By contrast, we develop a polynomial-time algorithm with a fast rate $\tO(d^{-(m-1)})$ under mild conditions. The optimality of our estimator is safeguarded by matching a minimax lower bound.

\paragraph{Low-rank tensor models.} There is a huge literature on structured tensor estimation under low-rank models, including CP models~\citep{kolda2009tensor}, Tucker models~\citep{zhang2018tensor}, and block models~\citep{wang2019multiway}. These models belong to parametric approaches, because they aim to explain the data with a finite number of parameters (i.e., decomposed factors). Our permuted smooth tensor model utilizes a different measure of model complexity than the usual low-rankness. We use \emph{infinite} number of parameters (i.e., smooth functions) to allow increasing model complexity. In this sense, our method belongs to nonparametric approaches. We emphasize that our permuted smooth tensor model does not necessarily include low-rank models. Compared to low-rank models, we utilize a different measure of \emph{model complexity}. When the underlying signal is precisely low-rank, then low-rank methods are preferred. However, if the underlying signal is high rank but has a specific structural pattern, then our nonparametric approach better captures the intrinsic model complexity. 
 
\paragraph{Nonparametric regression.} Our model is also related to nonparametric regression~\citep{tsybakov2009introduction}. One may view the problem~\eqref{eq:gmd} as a nonparametric regression, where the goal is to learn the function $f$ based on scalar response $\tY(i_1,\ldots,i_m)$ and design points $(\pi({i_1}),\ldots,\pi({i_m}))$ in $\bbR^m$; see Figure~\ref{fig:rate}(a). However, the \emph{unknown} permutation $\pi$ significantly influences the statistical and computational hardness of the problem. This latent $\pi$ leads to a phase transition in the estimation error; see Figure~\ref{fig:rate}(b) and Section~\ref{sec:lse}. We reveal two components of error for the problem, one for nonparametric error and the other for permutation error. The impact of unknown permutation hinges on tensor order and smoothness in an intriguing way, as shown in~\eqref{eq:rate}. This is clearly contrary to classical nonparametric regression. 

\paragraph{Graphon and hypergraphon.} Our work is also connected to graphons and hypergraphons. Graphon is a measurable function representing the limit of a sequence of exchangeable random graphs (matrices)~\citep{klopp2017oracle,gao2015rate,chan2014consistent}. Similarly, hypergraphon~\citep{zhao2015hypergraph,lovasz2012large} is introduced as a limiting function of $m$-uniform hypergraphs, i.e., a generalization of graphs in which edges can join $m$ vertices with $m\geq 3$. While both our model~\eqref{eq:gmd} and hypergraphon focus on function representations, there are two remarkable differences. First, unlike the matrix case where graphon is represented by bivariate functions~\citep{lovasz2012large}, hypergraphons for $m$-uniform hypergraphs should be represented as $(2^m-2)$-multivariate functions; see~\citet[Section 1.2]{zhao2015hypergraph}. Our framework~\eqref{eq:gmd} represents the function using $m$ coordinates only, and in this sense, the model shares the common ground as \emph{simple hypergraphons}~\citep{balasubramanian2021nonparametric}. We compare our method to earlier work in theory (Table~\ref{tab:comp} and Sections~\ref{sec:lse}-\ref{sec:borda}) and in numerical studies (Section~\ref{sec:sim}). Second, unlike typical simple hypergraphons where the design points are random, our generative model uses equal-sized fixed design points. The comparison of two approaches will be discussed in Sections~\ref{sec:model} and~\ref{sec:lse}.

\subsection{Notation and organization}
Let $\bbN, \bbN_{+}$ denote the set of non-negative integers and positive integers, respectively. We use $[d]=\{1,\ldots,d\}$ to denote the $d$-set for $d\in\bbN_{+}$. For a set $S$, $|S|$ denotes its cardinality and $\mathds{1}_S$ denotes the indicator function. For two positive sequences $\{a_d\},\{b_d\}$,  we denote $a_d\lesssim b_d$ if $\lim_{d\to\infty} a_d/b_d\leq c$ for some constant $c>0$, and $a_d\asymp b_d$ if $c_1\leq \lim_{d\to \infty} a_d/b_d\leq c_2$ for some constants $c_1,c_2>0$. Given a number $a\in\bbR$, the  function $\lfloor a\rfloor$ is the largest integer strictly smaller than $a$ and the ceiling function $\lceil a\rceil$ is the smallest integer no less than $a$. We use $\tO(\cdot)$ to denote big-O notation hiding logarithmic factors, and $\circ$ the function composition. 

Let $\Theta\in\bbR^{d\times \cdots \times d}$ be an order-$m$ $d$-dimensional tensor, $\pi\colon[d]\rightarrow[d]$ be an index permutation, and $\Theta(i_1,\ldots,i_m)$ the tensor entry indexed by $(i_1,\ldots,i_m)$. We sometimes also use shorthand notation $\Theta(\omega)$ for tensor entries with indices $\omega=(i_1,\ldots,i_m)$. We call a tensor \emph{a binary-valued tensor} if its entries take value on $\{0,1\}$-labels, and \emph{a continues-valued tensor} if its entries take values on a continuous scale. We define the Frobenius norm $\FnormSize{}{\Theta}^2=\sum_{\omega\in[d]^m}|\Theta(\omega)|^2$ for a tensor $\Theta$ and the $\infty$-norm $\mnorm{\mx}=\max_{i\in[d]}|x_i|$ for a vector $\mx=(x_1,\ldots,x_d)^T$.
We use $\Pi(d,d)=\{\pi\colon [d]\to[d]\}$ to denote all permutations on $[d]$ and use $\Pi(d,k)=\{\pi\colon [d]\to[k]\}$ to denote the collection of all onto mappings from $[d]$ to $[k]$. Given $\pi\in\Pi(d,k)$ and $\Theta\in\mathbb{R}^{k\times \cdots \times k}$, we use $\Theta\circ\pi$ to denote the $d$-dimensional tensor such that $(\Theta\circ\pi)(i_1,\ldots,i_m) = \Theta(\pi(i_1),\ldots,\pi(i_m))$ for all $(i_1,\ldots,i_m)\in[d]^m$. An event $A$ is said to occur \emph{with high probability} if $\mathbb{P}(A)$ tends to 1 as the tensor dimension $d\to\infty$.

The rest of the paper is organized as follows. Section~\ref{sec:model} presents the permuted smooth tensor model and its connection to smooth function representation. In Section~\ref{sec:tba}, we establish the approximation error based on block-wise polynomial approximation. Then, we develop two estimation algorithms with accuracy guarantees: the least-squares estimation and Borda count estimation. Section~\ref{sec:lse} presents a statistically optimal but computationally challenging least-squares estimator. Section~\ref{sec:borda} presents a polynomial-time Borda count algorithm with the same minimax optimal rate under an extra Lipschitz monotonicity assumption. Simulations and data analyses are presented in Section~\ref{sec:sim}. We conclude the paper with a discussion in Section~\ref{sec:discussion}. All proofs and extensions are deferred to the Appendix.

\section{Smooth tensor model with unknown permutation}\label{sec:model}
Suppose we observe an order-$m$ $d$-dimensional data tensor from the following model,
\begin{equation}\label{eq:obs}
\tY=\Theta\circ \pi+\tE,
\end{equation}
where $\pi\colon[d]\rightarrow[d]$ is an unknown latent permutation,  $\Theta\in \bbR^{d\times \cdots \times d}$ is an unknown signal tensor under certain smoothness (to be specified in the next paragraph), and $\tE$ is a noise tensor consisting of zero-mean, independent sub-Gaussian entries with variance bounded by $\sigma^2$. The general model allows both continuous- and binary-valued tensors. For instance, in binary tensor problems, the entries in $\tY$ are $\{0, 1\}$-labels from a Bernoulli distribution, and the entrywise variance of $\tE$ depends on the mean. For simplicity, we assume $\sigma=1$ throughout the paper. We call~\eqref{eq:obs} the Gaussian model if $\tE$ consists of i.i.d.\ $\tN(0,1)$ entries, and call~\eqref{eq:obs} the sub-Gaussian model if $\tE$ consists of independently (but not necessarily identically) distributed sub-Gaussian entries.

We now describe the smooth model on the signal $\Theta$. Assume there exists a multivariate function $f\colon [0,1]^m\rightarrow \bbR$ underlying the signal tensor, such that
\begin{align}\label{eq:rep}
\Theta(i_1,\ldots,i_m) = f\left({i_1\over d},\ldots,{i_m\over d}\right), \quad \text{for all }(i_1,\ldots,i_m)\in[d]^m.
\end{align}
For a multi-index $\kappa=(\kappa_1,\ldots,\kappa_m)\in\mathbb{N}^m$ and a vector $\mx=(x_1,\ldots,x_m)^T$, we denote $|\kappa|=\sum_{i\in[m]}\kappa_i$, $\kappa!=\prod_{i\in[m]}\kappa_i!$, $\mx^\kappa=\prod_{i\in[m]}x_i^{\kappa_i}$, and the derivative operator $\nabla_{\kappa}={\partial^{|\kappa|} \over \partial x^{\kappa_1}_1\cdots \partial x_m^{\kappa_m}}$. Assume the generative function $f$ in~\eqref{eq:rep} is in the $\alpha$-H\"older smooth family~\citep{wasserman2006all,tsybakov2009introduction}. 

\begin{defn}[$\alpha$-H\"older smooth]\label{eq:holder}
Let $\alpha> 0$ and $L>0$ be two positive constants. A function $f\colon [0,1]^m\rightarrow \bbR$ is called $\alpha$-H\"older smooth, denoted as $f\in\tF(\alpha,L)$, if 
\begin{align}\label{eq:defn}
\sum_{\kappa:|\kappa|=\lceil \alpha-1 \rceil}{1\over \kappa!}|\nabla_{\kappa} f (\mx)-\nabla_{\kappa} f(\mx_0)|\leq
L\| \mx-\mx_0\|_{\infty}^{\alpha-\lceil \alpha-1 \rceil}
\end{align}
holds for every $\mx,\mx_0\in[0,1]^m$.
\end{defn}

The H\"older smooth function class is one of the most popular function classes considered in the nonparametric regression literature~\citep{klopp2017oracle,gao2015rate}. 
In addition to the function class $\tF(\alpha, L)$, we also define the smooth tensor class based on discretization~\eqref{eq:rep}, 
\begin{equation}\label{eq:parameterP}
{\small \tP(\alpha,L)= \left\{\Theta\in\mathbb{R}^{d\times \cdots \times d} \colon\Theta \text{ is generated from~\eqref{eq:rep} and } f\in\tF(\alpha,L)\right\}.}
\end{equation}
Combining~\eqref{eq:obs} and~\eqref{eq:rep} yields our proposed \emph{permuted smooth tensor model}. 
The unknown parameters are the smooth tensor $\Theta \in \tP(\alpha, L)$ and latent permutation $\pi \in \Pi(d,d)$. The generative model is illustrated in Figure~\ref{fig:rate}(a) for the case $m=2$ (matrices). 
We provide an example to show the applicability of our permuted smooth tensor model.

\begin{example}[Co-authorship networks] Consider a co-authorship network consisting of $d$ nodes (authors) in total. We say there exists a hyperedge of size $m$ between nodes $(i_1,\ldots,i_m)$ if the authors $i_1,\ldots,i_m$ have co-authored at least one paper. The resulting $m$-uniform hypergraph is represented as an order-$m$ (symmetric) adjacency tensor. Our model is then expressed as
\begin{align}
    \bbE\tY(i_1,\ldots,i_m)&=\mathbb{P}(\text{$(i_1,\ldots,i_m)$ have co-authored together})
=f\left({\pi(i_1)\over d},\cdots, {\pi(i_m)\over d}\right),
\end{align}
for all $i_1< \cdots < i_m$. We interpret the permutation $\pi$ as the affinity measures of authors, and the function $f$ represents the $m$-way interaction between authors. The parametric model~\citep{wang2018learning} imposes logistic function $f(x_1,\ldots,x_m)=(1+\exp(-\beta x_1x_2\cdots x_m))^{-1}$. By contrast, our nonparametric model allows unknown $f$ and learns the function directly from data. 
\end{example}

Our model~\eqref{eq:rep} assumes an equal-sized grid design $\{1/d,2/d,\ldots,d/d\}$ from the generative function $f$. One can also extend the model to unbalanced designs; i.e., the signal tensor is generated from $f$ based on 
\begin{align}\label{eq:randommodel}
\Theta(i_1,\ldots,i_m) = f(x_{i_1},\ldots,x_{i_m}), \quad \text{for all} \quad (i_1,\ldots,i_m)\in[d]^m. 
\end{align}
where the design points $\{x_i\}_{i\in[d]}$ may be modeled as either fixed latent variables or i.i.d.\ random variables from a probability distribution supported on $[0,1]$. The random design model has been developed in the literature of graphons and hypergraphons~\citep{chan2014consistent,gao2015rate,klopp2017oracle,balasubramanian2021nonparametric}. Our paper will focus on the grid design~\eqref{eq:rep}, and we will also discuss the extension to~\eqref{eq:randommodel} whenever possible.

Our model has equal dimension and the same permutation along $m$ modes. The results for non-symmetric tensors with $m$ distinct permutations are similar but require extra notations; for brevity exclude this case from the theoretical analysis and instead evaluate its empirical performance in the Appendix A.2.

\section{Block-wise tensor estimation}\label{sec:tba}
Our general strategy for estimating the permuted smooth tensor is based on the block-wise tensor approximation. In this section, we first introduce the tensor block model~\citep{wang2019multiway,han2022exact}. Then, we extend this model to block-wise polynomial approximation.

\subsection{Tensor block model}\label{subsec:bm}
Tensor block models describe a checkerboard pattern in the signal tensor. The block model provides a meta structure to many popular models including the low-rankness~\citep{gao2016optimal} and isotonic tensors~\citep{pananjady2022isotonic}. Here, we use tensor block models as a foundation for estimating permuted smooth models. 

Specifically, suppose that there are $k$ clusters among $d$ entities, and the cluster assignment is represented by a clustering function $z \colon[d]\rightarrow[k]$. Then, the tensor block model assumes that the entries of signal tensor $\Theta\in\bbR^{d\times \cdots \times d}$ take values from a core tensor $\tS\in\bbR^{k\times\cdots\times k}$ according to the clustering function $z$; that is,
\begin{align}\label{eq:block}
    \Theta(i_1,\ldots,i_m) = \tS(z(i_1),\ldots,z(i_m)), \quad \text{ for all } (i_1,\ldots,i_m)\in[d]^m.
\end{align}
Here, the core tensor $\tS$ collects the entry values of $m$-way blocks; the core tensor $\tS$ and clustering function $z\in \Pi(d,k)$ are parameters of interest. A tensor $\Theta$ satisfying~\eqref{eq:block} is called a block-$k$ tensor. Tensor block models allow various data types, as shown below. 

\begin{example}[Gaussian tensor block model] Let $\tY$ be a continuous-valued tensor. The Gaussian tensor block model draws independent normal entries according to $\tY(i_1,\ldots,i_m)\stackrel{\text{ind}}{\sim} N(\tS(z(i_1),\ldots,$\\ $z(i_m)),\sigma^2)$. The mean model belongs to~\eqref{eq:block}, and noise entries are i.i.d.\ $N(0,\sigma^2)$. The Gaussian tensor block model has served as the statistical foundation for many tensor clustering algorithms \citep{wang2019multiway,han2022exact}.
\end{example}

Tensor block models have shown great success in discovering hidden group structures for many applications~\citep{wang2019multiway,han2022exact}. Despite the popularity, the constant block assumption is insufficient to capture delicate structure when the signal tensor is complicated. 
This parametric model aims to explain data with a finite number of blocks; such an approach is useful when the sample outsizes the parameters. Our nonparametric model~\eqref{eq:rep}, by contrast, uses an infinite number of parameters (i.e., smooth functions) to allow growing model complexity. We change the goal of the tensor block model from clustering to approximating the generative function $f$ in~\eqref{eq:rep}. In our setting, the number of blocks $k$ should be interpreted as a resolution parameter (i.e., a bandwidth) of the approximation, similar to the notion of the number of bins in histogram and polynomial regressions~\citep{wasserman2006all}. 

\subsection{Block-wise polynomial approximation}\label{sec:poly}
The block tensor~\eqref{eq:block} can be viewed as a discrete version of piece-wise \emph{constant} function. This connection motivates the use of block-wise \emph{polynomial} tensors to approximate $\alpha$-H\"older functions. Now we extend~\eqref{eq:block} to the block-wise polynomial models. 

We introduce some additional notation. For a given block number $k$, we use $z\in\Pi(d,k)$ to denote the balanced clustering function that partitions $[d]$ into $k$ equally-sized clusters, such that 
$z(i) = \lceil ki/d\rceil$, for all $i\in[d]$.
The collection of inverse images $\{z^{-1}(j)\colon j\in[k]\}$ is a partition of $[d]$ into $k$ disjoint and equal-sized subsets. We use $\tE_k$ to denote the $m$-way balanced partition, i.e., a collection of $k^m$ disjoint and equal-sized subsets in $[d]^m$, such that 
\begin{align}\label{eq:blockind}
    \tE_k = \{z^{-1}(j_1)\times\cdots\times z^{-1}(j_m)\colon (j_1,\ldots,j_m)\in [k]^m\}.
\end{align}
Let $\Delta\in \tE_k$ denote an element in $\tE_k$. We propose to approximate the signal tensor $\Theta$ in~\eqref{eq:rep} using degree-$\ell$ polynomial tensors within each block $\Delta\in \tE_k$. Specifically, let $\caliB(k,\ell)$ denote the class of block-$k$ degree-$\ell$ polynomial tensors, i.e., 
\begin{align}\label{eq:polynomial}
    \caliB(k,\ell) = \bigg\{&\tB\in \bbR^{d\times \cdots \times d}\colon \tB(\omega) = \sum_{\Delta\in\tE_k}\text{Poly}_{\ell,\Delta}(\omega)\mathds{1}\{\omega\in\Delta\}\text{ for all } \omega\in[d]^m\bigg\},
\end{align}
where $\text{Poly}_{\ell,\Delta}(\cdot)$ denotes a degree-$\ell$ polynomial function in $\bbR^m$, with coefficients depending on block $\Delta$; that is, a constant function $\text{Poly}_{0,\Delta}(\omega)= \beta^0_{\Delta}$ for $\ell=0$, a linear function $\text{Poly}_{1,\Delta}(\omega)=\langle \boldsymbol{\beta}_{\Delta},\omega\rangle+\beta_{\Delta}^0$ for $\ell=1$, and so forth. Here $\beta^0_{\Delta}$ and $\boldsymbol{\beta}_{\Delta}$ denote unknown coefficients in polynomial function. Note that the degree-0 polynomial block tensor reduces to the equal-sized constant block model~\eqref{eq:block}. We generalize the constant block model to degree-$\ell$ polynomial block tensors~\eqref{eq:polynomial},  analogous to the generalization from $k$-bin histogram to $k$-piece-wise polynomial regression in nonparametric statistics~\citep{wasserman2006all}.

The smoothness of the function $f$ in \eqref{eq:rep} plays an important role in the block-wise polynomial approximation. The following proposition explains the role of smoothness in the approximation. 

\begin{prop}[Block-wise polynomial tensor approximation]\label{lem:approx}
Suppose $\Theta\in\tP(\alpha,L)$. Then, for every block number $k\leq d$ and degree $\ell\in \bbN$, we have the approximation error
\begin{align}
   \inf_{\tB\in\caliB(k,\ell)} \frac{1}{d^m}\FnormSize{}{\Theta-\tB}^2\lesssim \frac{L^2}{k^{2\min(\alpha,\ell+1)}}.
\end{align}
\end{prop}
Proposition~\ref{lem:approx} implies that we can always find a block-wise polynomial tensor close to the signal tensor generated from $\alpha$-H\"older smooth function $f$. The approximation error decays with the block number $k$ and degree $\min(\alpha, \ell +1)$. 

\section{Statistical limits via least-squares estimation}\label{sec:lse}
We develop two estimation methods based on the block-wise polynomial approximation. We first introduce a statistically optimal but computationally inefficient least-squares estimator. The least-squares estimation serves as a statistical benchmark because of its minimax optimality. 
In Section~\ref{sec:borda}, we present a polynomial-time algorithm with a provably same optimal rate under monotonicity assumptions.

We propose the least-squares estimation for model~\eqref{eq:obs} by minimizing the Frobenius loss over the block-$k$ degree-$\ell$ polynomial tensor family $\caliB(k,\ell)$ up to permutations, 
\begin{align}\label{eq:lseopt}
    (\hat\Theta^{\text{LSE}},\hat \pi^{\text{LSE}}) &= \argmin_{\Theta\in\caliB(k,\ell), \pi\in \Pi(d,d)}\FnormSize{}{\tY-\Theta\circ\pi}.
\end{align}
The least-squares estimator $(\hat\Theta^{\text{LSE}},\hat\pi^{\text{LSE}})$ depends on two tuning parameters: the number of blocks $k$ and the polynomial degree $\ell$. The optimal choice $(k^*,\ell^*)$ will be provided below. 

Theorem~\ref{thm:LSE} establishes the error bound for the least-squares estimator~\eqref{eq:lseopt}. Note that $\Theta$ and $\pi$ are in general not separably identifiable; for example, when the true signal is a constant tensor, then every permutation $\pi\in \Pi(d,d)$ gives equally good fit in statistics. We assess the estimation error on the composition $\Theta \circ \pi$ to avoid this identifiability issue. For two order-$m$ $d$-dimensional tensors $\Theta_1, \Theta_2$, define the mean squared error (MSE) as $\textup{MSE}(\Theta_1,\Theta_2) = d^{-m}\FnormSize{}{\Theta_1-\Theta_2}^2$. 
 
\begin{thm}[Least-squares estimation error]\label{thm:LSE} Let $m\geq 2$. Consider the sub-Gaussian model~\eqref{eq:obs} with $\Theta\in\tP(\alpha,L)$. Let $(\hat\Theta^{\textup{LSE}},\hat\pi^{\textup{LSE}})$ denote the least-squares estimator in \eqref{eq:lseopt} with a given $(k,\ell)$. Then, for every $k\leq d$ and degree $\ell\in \bbN$, we have
\begin{align}\label{eq:rateMSE}
\textup{MSE}(\hat\Theta^{\textup{LSE}}\circ\hat\pi^{\textup{LSE}},\ \Theta\circ\pi)&\lesssim \KeepStyleUnderBrace{\frac{L^2}{k^{2\min(\alpha,\ell+1)}}}_{\textup{approximation error}}+ \KeepStyleUnderBrace{\frac{k^m(\ell+m)^\ell}{d^m}}_{\textup{nonparametric error}}+\KeepStyleUnderBrace{\frac{\log d}{d^{m-1}}}_{\textup{permutation error}}
\end{align}
with high probability. In particular, setting {\small $\ell^* = \min(\lceil \alpha-1\rceil,(m-2)(m+1)/2)$} and $k^*=c_1 d^{m/(m+2\min(\alpha,\ell^*+1))}$
yields the optimized error rate
\begin{align}\label{eq:rates}
     \eqref{eq:rateMSE} \lesssim   \textup{Rate}(d) =
     \begin{cases} 
    c_2 d^{-\frac{2m\alpha}{m+2\alpha}}, & \text{ when } \alpha < m(m-1)/2,\\
     c_3 d^{-(m-1)}\log d, &\text{ when } \alpha \geq m(m-1)/2.
    \end{cases}
\end{align}
Here, the function $\textup{Rate}(d)$ is given in~\eqref{eq:rate}, and the constants $c_1,c_2,c_3>0$ depend on the model configuration $(m, L,\alpha)$ but not on the tensor dimension $d$. See Appendix~\ref{sec:proofLSE} for the constant values. 

\end{thm}

We discuss the asymptotic error rate as $d\rightarrow \infty$ while treating other model configurations fixed. The final least-squares estimation rate~\eqref{eq:rates} has two sources of error: the nonparametric error $d^{-{2m\alpha\over m+2\alpha}}$ and the permutation error $d^{-(m-1)}\log d$. Intuitively, in the tensor data analysis problem, we can view each tensor entry as a data point, so sample size is the total number of entries, $d^m$. The unknown permutation results in $\log(d!)\approx d\log d$ complexity, whereas the unknown generative function results in $d^{-2m\alpha/(m+2\alpha)}$ nonparametric complexity. When the function $f$ is smooth enough, estimating the function $f$ becomes relatively easier compared to estimating the permutation $\pi$. This intuition coincides with the fact that the permutation error dominates the nonparametric error when  $\alpha\geq m(m-1)/2$.

\begin{rmk}[Comparison to non-parametric regression]
We now compare our results with existing work in the literature. In the vector case with $m = 1$,  our model reduces to the one-dimensional regression problem such that
$y_i  = \theta_{\pi(i)}+\epsilon_i, \text{ for all } i\in [d],$
where $\theta_i=f(i/d)$ and unknown $\pi \in \Pi(d,d)$.
A similar analysis of Theorem~\ref{thm:LSE} shows the error rate
\begin{align}\label{eq:m1}
    \frac{1}{d}\sum_{i\in[d]}(\hat\theta_i^{\text{LSE}}-\theta_i)^2 \lesssim \left(d^{-{2\alpha\over 2\alpha+1}}+\log d\right),
\end{align}
under the choice of $\ell^* = 0$ and $k^* \asymp d^{1\over 1+2\min(\alpha,1)}$.
Notice that  $d^{-2\alpha/(2\alpha+1)}$ is the  classical nonparametric minimax rate for $\alpha$-Hölder smooth functions~\citep{tsybakov2009introduction} with \emph{known} permuted design points $\{\pi(i)\}_{i=1}^d$. By contrast, our model involves \emph{unknown} $\pi$, which results in the non-vanishing permutation rate $\log d$ in~\eqref{eq:m1}.
\end{rmk}

\begin{rmk}[Breaking previous limits on matrices/tensors]
In the matrix case with $m=2$, Theorem~\ref{thm:LSE} implies that the best rate is obtained under $\ell^*=0$, i.e., the block-wise \emph{constant} approximation. This result is consistent with existing literature on smooth graphons~\citep{bickel2009nonparametric,gao2015rate,klopp2017oracle}, where constant block model (see Section~\ref{subsec:bm}) is developed for estimation. In the tensor case with $m\geq 3$, earlier work~\citep{balasubramanian2021nonparametric} suggests that constant block approximation ($\ell^*=0$) may remain minimax optima. Our Theorem~\ref{thm:LSE} disproves this conjecture, and we reveal a much faster rate $d^{-(m-1)}$ compared to the conjectured lower bound $d^{-2m/( m+2)}$~\citep{balasubramanian2021nonparametric} for sufficiently smooth tensors. We demonstrate that a polynomial up to degree $(m-2)(m+1)/2$ is sufficient (and necessary; see Theorem~\ref{thm:minimax} below) for accurate estimation of order-$m$ permuted smooth tensors. 
For example, permuted $\alpha$-smooth tensors of order-3 require quadratic approximation $(\ell^*=2)$ with $k^*\asymp d^{1/3}$ blocks, for all $\alpha\geq 2$. The results show the clear difference between matrices and higher-order tensors. 
\end{rmk}

We now show that the rate in~\eqref{eq:rates} cannot be improved. The lower bound is information-theoretical and thus applies to all estimators including, but not limited to, the least-squares estimator~\eqref{eq:lseopt} and Borda count estimator introduced in next section. 
\begin{thm}[Minimax lower bound]\label{thm:minimax}The estimation problem based on the Gaussian model~\eqref{eq:gmd} obeys the minimax lower bound 
\begin{equation}\label{eq:minimax}
\inf_{(\hat \Theta,\hat \pi)}\sup_{
\Theta\in \tP(\alpha,L), \pi\in \Pi(d,d)} \mathbb{P}\left(\textup{MSE}(\hat\Theta\circ\hat\pi,\ \Theta\circ\pi) \gtrsim \textup{Rate}(d) \right) \geq 0.8,
\end{equation}
where the function $\textup{Rate}(d)$ is given in~\eqref{eq:rate}.
\end{thm}
The lower bound in~\eqref{eq:minimax} matches the upper bound in~\eqref{eq:rates}, demonstrating the statistical optimality of  the convergence speed $\text{Rate}(d)=d^{-{2m\alpha  \over 2\alpha+m}}\vee d^{-(m-1)}\log d$. The two-component error reveals the intrinsic model complexity. In particular, the permutation error $d^{-(m-1)}$ dominates nonparametric error $d^{-{2m\alpha/(m+2\alpha)}}$ for sufficiently smooth tensors. This phenomenon is contrary to classical nonparametric regression. 

\begin{rmk}[Phase transition]\label{rmk:phase} We conclude this section by summarizing an interesting phase transition phenomenon.  Figure~\ref{fig:rate}(b) plots the optimal convergence speed $\text{Rate}(d)$. The impact of unknown permutation hinges on the tensor order and smoothness. The accuracy improves with respect to smoothness in the regime $\alpha\leq m(m-1)/2$; however, in the regime $\alpha>m(m-1)/2$, the accuracy becomes a constant with respect to smoothness.  
The result implies a polynomial of degree $\approx (m-2)(m+1)/2$ is sufficient for accurate recovery of order-$m$ tensors, whereas higher degree brings no further benefits. This full picture of error dependence, to our best knowledge, is new to the literature. 
\end{rmk}

\begin{rmk}[Extension to random designs] 
Our theory assumes fixed balanced designs~\eqref{eq:rep}. We can extend the results to random designs~\eqref{eq:randommodel} by imposing extra conditions on the distribution $\{x_i\}_{i\in[d]}$ and generalizing the definition of $\tP(\alpha,L)$, $\caliB(k,\ell)$ and $(\hat \Theta^{\text{LSE}}, \hat \pi^{\text{LSE}})$. Theorems~\ref{thm:LSE} and~\ref{thm:minimax} remain valid under this modification; see the Appendix F for details. 
\end{rmk}

\section{Computational limits and polynomial-time algorithms}\label{sec:borda}

We point out that the least-squares estimator $(\hat\Theta^{\text{LSE}},\hat\pi^{\text{LSE}})$ usually requires exponential-time algorithms, even in the simple matrix case~\citep{gao2015rate}. Specifically, this section will show the general non-existence of \emph{polynomial-time algorithms} with rate~\eqref{eq:rate} under certain conditions. This fundamental computational limit highlights the need of extra assumptions to achieve the optimal rate~\eqref{eq:rate} from a computational standpoint.

\subsection{Computational limits under HPC detection conjecture}\label{sec:complimit}
The hypergraphic planted clique detection conjecture plays an important role in constructing the computational limits of our problem. We briefly introduce the HPC hardness conjecture. 

Consider an $m$-uniform hypergraph $G = (V,E)$, where $V$ is a set of vertices and $E$ is a set of hyperedges. An Erd\"os-R\'enyi random hypergraph, denoted by $\tG_m(d,1/2)$, is a random $m$-uniform hypergraph with $d$ vertices and probability 1/2 for each of the hyperedge connections. The hypergraphic planted clique (HPC) with clique size $\kappa>0$, denoted by $\tG_m(d,1/2,\kappa)$, is generated from Erd\"os-R\'enyi random hypergraph in the following way. First we generate an Erd\"os-R\'enyi random hypergraph from $\tG_m(d,1/2)$. Then we independently pick $\kappa$ vertices with uniform probability from $d$ vertices. Finally, we obtain the $\tG_m(d,1/2,\kappa)$ by including only the hyperedges whose vertices all belong to the picked $\tau$ vertices. The HPC detection refers to the hypothesis testing problem, 
\begin{align}\label{eq:HPC}
    H_0\colon G\sim \tG_m(d,1/2)\quad\text{v.s.}\quad H_1\colon G\sim \tG_m(d,1/2,\kappa).
\end{align}
The earlier work \citep{luo2022tensor} presents the following hardness conjecture for testing \eqref{eq:HPC}.
\begin{conj}[HPC detection conjecture \citep{luo2022tensor}]\label{conj:1} Consider the HPC problem in \eqref{eq:HPC} and let $m\geq 2$ is be fixed integer. Suppose 
$
   \limsup_{d\rightarrow\infty} \log\kappa/\log\sqrt{d}\leq 1-\tau,\quad \text{for any }\tau >0.
   $
Then, for any polynomial-time test sequence $\{\phi\}_d\colon G \mapsto \{0,1\}$, we have 
\[
\liminf_{d\rightarrow\infty} \left\{\mathbb{P}_{H_0}(\phi(G) = 1)+ \mathbb{P}_{H_1}(\phi(G) = 0)\right\}\geq \frac{1}{2},
\]
\end{conj}

Now we construct the computational lower bound based on Conjecture~\ref{conj:1}. 

\begin{thm}[Computational lower bound under general designs]\label{thm:nopoly} Assume Conjecture~\ref{conj:1} holds. Define the smooth tensor class under general designs,
$
\tP_{\textup{gen}}=\{\Theta: \Theta \text{ is generated from~\eqref{eq:randommodel}}$
with fixed latent variables $\{x_i\}_{i\in[d]}$ and $f\in \tF(\alpha,L)\}$. 
Consider the Gaussian model~\eqref{eq:obs} with $\alpha> m/2$. There exists no polynomial algorithm that achieves the statistical optimal convergence $\text{Rate}(d)$; i.e., 
\begin{align}\label{eq:compthm}
\frac{1}{\textup{Rate}(d)}\inf_{\hat \Theta \textup{ Polynomial-time}}\sup_{\Theta\in\tP_{\textup{gen}}}\textup{MSE}(\hat\Theta,\Theta)\rightarrow \infty, \quad \text{as}\ d\to\infty.
\end{align}
\end{thm}
Theorem~\ref{thm:nopoly} shows the \emph{impossibility} of a polynomial-time estimator to achieve the optimal statistical rate in the general model. The intuition in the proof is to show the best bound for \emph{polynomial-time} tensor estimation as $d^{-m/2}$, in the absence of extra model structures. The condition $\alpha>m/2$ is a technical assumption to facilitate the proof. 
Theorem~\ref{thm:nopoly} is not a weakness of our proposed estimator; rather, \eqref{eq:compthm} reveals the non-avoidable statistical-computational gap as a nature of the smooth tensor estimation problem. 

\subsection{Borda count algorithm}
The earlier section has shown the impossibility of polynomial-time algorithms in the general model. In this section, we restrict ourselves to a sub-model with extra monotonicity structures; this structure makes polynomial-time algorithms possible. We introduce a notion of marginal monotonicity for the generative functions.

\begin{defn}[Marginal monotonicity]\label{eq:bdefn}
Let $\beta\geq 0$ be a non-negative constant. A function $f\colon[0,1]^m \rightarrow \bbR$ is called $\beta$-monotonic, denoted as $f\in\tM(\beta)$, if the following holds for all $d$:
\begin{align}\label{eq:monotonic}
   \left({i-j\over d}\right)^{1/\beta}  \lesssim g(i)-g(j),\quad \text{for all } j\leq i\in[d], 
   \end{align}
where we define the \emph{score function} $g(i) ={d^{-(m-1)}}\sum_{(i_2,\ldots,i_m)\in [d]^{m-1}}f(i,i_2,\ldots,i_m) $ for $i\in[d]$.
\end{defn}

We refer to $\tF(\alpha,L)\cap \tM(\beta)$ as the monotonic-plus-smooth function class. This class was initially proposed in previous literature of graphons. The work~\citep{chan2014consistent} proposes the Lipschitz monotonic function to facilitate the analysis of sorting-merging algorithm for matrix estimation; their setting is a special case of our Definition 2 with $(\alpha,\beta,m)=(1,1,2)$. Inspired by earlier work, we consider the similar monotonic-plus-smooth function class $\tF(\alpha,L)\cap \tM(\beta)$ under general configuration $\{(\alpha,\beta,m): \alpha>0,\ 0<\beta\min(\alpha,1)\leq 1, \ m\geq 2\}$. Note that the constraint $\beta \min(\alpha, 1) \leq 1$ is due to the automatic constraint between joint smoothness and marginal smoothness. A large value of $\beta$ in~\eqref{eq:monotonic} implies the steepness of $g$.

Now we introduce the \emph{Borda count} estimation that consists of two stages. The procedure is illustrated in Figure~\ref{fig:borda}. Define the empirical score function $\tau\colon [d]\rightarrow\bbR$ as 
\[
\tau(i) = \frac{1}{d^{m-1}}\sum_{(i_2,\ldots,i_m)\in[d]^{m-1}
} \tY(i,i_2,\ldots,i_m).
\] 

\begin{enumerate}[wide, labelwidth=0pt, labelindent=0pt]
  \item {\bf Sorting stage:} The sorting stage is to rearrange the observed tensor $\tY$ so that the score function $\tau$ of sorted tensor is  monotonically increasing. Define a permutation $\hat\pi^{\text{BC}}$ such that
    \begin{align}\label{eq:permute}
        \tau \circ (\hat\pi^{\text{BC}})^{-1}(1) \leq  \tau\circ (\hat\pi^{\text{BC}})^{-1}(2)
        \leq \cdots 
        \leq
        \tau\circ (\hat\pi^{\text{BC}})^{-1}(d).
    \end{align}
    Then, we obtain sorted observation $\tilde\tY = \tY \circ (\hat\pi^{\text{BC}})^{-1}$, illustrated in Figure~\ref{fig:borda}(b).
\begin{figure}[t]
    \centering
    \includegraphics[width = \textwidth]{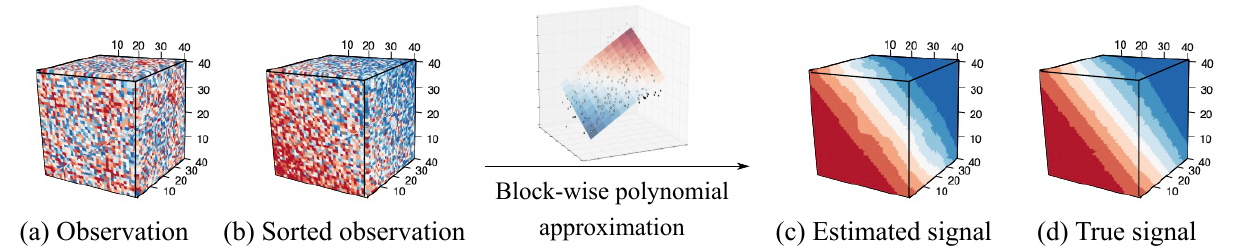}
    \caption{Illustration of Borda count estimation. We first sort tensor entries using the proposed procedure, and then estimate the signal by block-wise polynomial approximation.}
    \label{fig:borda}
\end{figure}

 \item {\bf  Block-wise polynomial approximation stage:} Given sorted observation $\tilde \tY $, we estimate the signal tensor by block-wise polynomial tensor based on the following optimization,
    \begin{align}\label{eq:bclse}
        \hat\Theta^{\text{BC}} = \argmin_{\Theta\in\caliB(k,\ell)}\FnormSize{}{\tilde\tY-\Theta},
    \end{align}
    where $\tB(k,\ell)$ denotes the block-$k$ degree-$\ell$ tensor class  in~\eqref{eq:polynomial}. An example of this procedure is shown in Figure~\ref{fig:borda}(c). The estimator $\hat\Theta^{\text{BC}}$ depends on two tuning parameters: the number of blocks $k$ and polynomial degree $\ell$. The optimal choice of $(k^*,\ell^*)$ is provided in Theorem~\ref{thm:BC}.
  \end{enumerate}

The following theorem ensures the statistical accuracy of the Borda count estimator.
\begin{thm}[Estimation error for Borda count algorithm under marginal monotonicity]\label{thm:BC} Consider the sub-Gaussian model~\eqref{eq:obs} with $f\in \tF(\alpha,L)\cap \tM(\beta)$.
Denote a constant $c(\alpha,\beta,m):= \frac{m(m-1)\beta\min(\alpha,1)}{2(m-(m-1)\beta\min(\alpha,1))}$. Let $(\hat\Theta^{\textup{BC}},\hat\pi^{\textup{BC}})$ be the Borda count estimator with {\footnotesize $
\ell^* = \min\left(\lceil\alpha-1\rceil,\lfloor c(\alpha,\beta,m)\rfloor\right)$} and $k^* \asymp d^{m/ (m+2\min(\alpha,\ell^*+1))}$ in~\eqref{eq:bclse}. Then, we have
\begin{align}\label{eq:rateBC}
\textup{MSE}(\hat\Theta^{\textup{BC}}\circ\hat\pi^{\textup{BC}},\ \Theta\circ\pi)  \lesssim  
   \begin{cases}
   d^{-{2m\alpha\over m+2\alpha}} & \text{ when } \alpha <c(\alpha,\beta,m),\\
   \left({\log d\over d^{m-1}}\right)^{\beta\min(\alpha,1)}&\text{ when } \alpha \geq c(\alpha,\beta,m).
   \end{cases}
\end{align}
\end{thm}
The estimation bound~\eqref{eq:rateBC} comes from the approximation error (Proposition~\ref{lem:approx}), nonparametric error (Theorem~\ref{thm:LSE}), and permutation error (Lemma~\ref{lem:permute}). 

\begin{rmk}[Sufficiently smooth tensors]\label{rm2} When the generative function is infinitely smooth ($\alpha =\infty$) with Lipschitz monotonic score $(\beta=1)$,
our estimation error~\eqref{eq:rateBC} becomes
\begin{equation}\label{eq:infinite}
\textup{MSE}(\hat\Theta^{\textup{BC}}\circ\hat\pi^{\textup{BC}},\ \Theta\circ\pi)\lesssim d^{-(m-1)}\log d,
\end{equation}
under the choice of degree and block number
\begin{equation}\label{eq:choice}
\ell^*= {(m-2)(m+1)/2}\quad \text{ and } \quad k^* \asymp d^{m\over m +2(\ell^*+1)}.
\end{equation}
Now, we compare the rate~\eqref{eq:infinite} with the classical low-rank estimation~\citep{wang2018learning,zhang2018tensor,kolda2009tensor}. The low-rank tensor model with a constant rank is known to have MSE rate $\tO(d^{-(m-1)})$~\citep{wang2018learning}. Our infinitely smooth tensor model achieves the same rate up to the negligible log term. 
\end{rmk}

\paragraph{Hyperparameter tuning.}
Our algorithm has two tuning parameters $(k,\ell)$. The theoretically optimal choices of $(k,\ell)$ are given in Theorems~\ref{thm:LSE} and~\ref{thm:BC}. In practice, since model configuration is unknown, we search $(k,\ell)$ via cross-validation. 
Based on our theorems, a polynomial of degree $\ell^*=(m-2)(m+1)/2$ is sufficient for accurate recovery of order-$m$ tensors, whereas higher degree brings no further benefit. The practical impacts of hyperparameter tuning are investigated in Section~\ref{sec:sim}.

\subsection{Possible relaxation of monotonicity} 

We emphasize that the permutation in our model~\eqref{eq:obs} allows for flexibility in our monotonic-plus-smooth assumptions. Specifically, strict constraint of the target function $f\in \tF(\alpha,L)\cap \tM(\beta)$ is not necessary; our results remain valid if a permutation $\pi$ exists such that $(f\circ \pi)$ belongs to that class. We provide two examples to illustrate the flexibility. 
  
 \begin{example}[Non-monotonic function made monotonic by permutation]\label{ex:free}Consider the quadratic function $f(x,y,z)=(x-0.5)^2+yz$. Although $f$ is non-monotonic, we find that the tensor generated from $f$ can be equivalently represented, up to a small perturbation, as $\bar f\circ \pi$, where $\pi$ is a permutation that reorders the indices to achieve monotonicity, and $\bar f$ is a monotonic-plus-smooth function. The function $\bar f$ is constructed by reordering the indices along the $x$-axis index to ensure marginal monotonicity. The exact expressions of $(\pi, \bar f)$ are provided in the Appendix B. This alternative representation allows us to apply Borta count algorithm and Theorem~\ref{thm:BC} to this non-monotonic function. 
 \end{example}
 
 \begin{example}[Decomposable monotonicity with Lipschitz smooth factors] Let $R\in\mathbb{N}_{+}$ be a constant, and $\{g_{r,i}(\cdot) \colon[0,1]\to \bbR\}$ be a set of Lipschitz smooth ($\alpha=1$) functions for all $(r,i)\in[R]\times [m]$. Based on earlier example~\eqref{ex:free}, every univariate Lipschitz smooth function can be made monotonic-plus-smooth up to permutations. 
Therefore, decomposable smooth functions of the form 
$
 f(x_1,\ldots,x_m)=\sum_{r\in[R]}g_{r,1}(x_1)\cdots g_{r,m}(x_m)
$
and the corresponding low-rank tensors are also monotonic-plus-smooth up to permutations.
\end{example}

\begin{rmk}[Other monotonicity assumptions]\label{rmk:BC2}
One may also consider other assumptions such as isotonic functions~\citep{pananjady2022isotonic}. 
Define an isotonic function class $\tM$ 
\begin{align}\label{eq:iso}
    \tM = \{f\colon[0,1]^m\rightarrow\mathbb{R}\ \big|\ f(x_1,\ldots,x_m)\leq f(x_1',\ldots,x_m')\text{ when } x_i\leq x_i' \text{ for } i\in[m]\}.
\end{align}
The isotonic functions~\eqref{eq:iso} concerns the \emph{joint} monotonicity, where as our monotonicity~\eqref{eq:bdefn} concerns the \emph{marginal} monotonicity. The isotonic functions belong to $\tM(\beta)$ with $\beta=0$ based on our Definition~\ref{eq:bdefn}. In fact, we find that for isotonic functions, Theorem~\ref{thm:BC} can be sharpened to achieve the rate $\text{rate}(g)=d^{-{2m\alpha\over m+2\alpha}} \vee d^{-(m-1)}\log d$. Intuitively, while $\hat{\pi}^{\text{BC}}$ does not consistently estimate $\pi$, the tensor $\Theta \circ \hat{\pi}^{\text{BC}}$ is still a good consistent estimator of $\Theta \circ \pi$. The detailed statement can be found in the Appendix C.6.
\end{rmk}

We show below that, in the special Lipschitz situation $\alpha=\beta=1$, neither the marginal monotonic or the isotonic assumptions alter the minimax rate of Theorem~\ref{thm:minimax}, except for a negligible logarithmic factor. Specifically, let $\tM$ be either the isotonic function class~\eqref{eq:iso} or Lipschitz monotonic function class with $\beta=1$ in~\eqref{eq:monotonic}. We establish the following minimax lower bound for the class $\tF(1,L)\cap \tM$. 

\begin{thm}[Statistical minimax lower bound for monotonic-plus-Lipschitz functions]\label{thm:mminimax}
Consider the Lipschitz smooth tensor class with extra monotonicity 
\[
\tP_{M} = \{\Theta\colon \Theta \text{ is generated from~\eqref{eq:rep}} \text{ and } f\in\tF(1,L)\cap \tM\}.
\] 
Then, the estimation problem based on Gaussian model \eqref{eq:gmd} obeys the minimax lower bound
\begin{align}\label{eq:mconclusion}
    \inf_{(\hat\Theta,\hat\pi)}\sup_{\Theta\in\tP_{M},\pi\in\Pi(d,d)}\mathbb{P}\left(\textup{MSE}(\hat\Theta\circ\hat\pi,\ \Theta\circ\pi) \gtrsim \textup{Rate}(d)/\log d \right) \geq 0.8.
\end{align}
\end{thm}
Theorem~\ref{thm:mminimax} shows that, in the special Lipschitz situation $\alpha=\beta=1$, the extra monotonicity assumption renders no changes to the minimax optimal rate of the problem, except for a negligible logarithmic factor. Therefore, our Borda count algorithm achieves computational and statistical optimality in the region $\tF(1,L)\cap \tM$. The optimality may not be attained in the absence of monotonicity, i.e., in the situation $\tF(\alpha,L)/\tM$. Since the statistical-computational gap is non-avoidable in general, our imposed monotonicity assumptions fill in the gap in several cases. For general $(\alpha,\beta)$, the optimality of Borda count algorithm is unknown; we discuss the proof challenges in the Appendix E.

\section{Numerical analysis}\label{sec:sim}
\subsection{Synthetic data}
We simulate order-3 $d$-dimensional tensors based on the permuted smooth tensor model~\eqref{eq:rep}. Both symmetric and non-symmetric tensors are investigated. The symmetric tensors are generated based on functions $f$ in Table~\ref{tb:md}, and the non-symmetric set-up is described in the Appendix A.2.
\begin{table}[htp]
    \centering
    \resizebox{.9\textwidth}{!}{%
    \begin{tabular}{c|c|c|c|c|c}
        Model ID  &  $f(x,y,z)$ & CP rank & Tucker rank &  $\geq(\alpha,\beta)$ & Isotonic \\\hline
        1 &    $xyz$ & 1 & $(1,1,1)$ &$(\infty, 1)$ &$\surd$ \\
        2 &    $(x+y+z)/3$ & 3 & $(2,2,2)$ &$(\infty,1)$&$\surd$\\
        3 & $(1+\exp(-(3x^2+3y^2+3z^2))^{-1}$ &9& $(4,4,4)$ & $(\infty,1/2)$ &$\surd$\\
        4 & $\log(1+\max(x,y,z))$ &$\geq 100$& $\geq (50,50,50)$ &$(1,1)$&$\surd$\\
        5 &  $\exp\left(-\max(x,y,z)-\sqrt{x}-\sqrt{y}-\sqrt{z}\right)$ &$\geq 100$& $\geq (50,50,50$) & $(1/2,1)$ &$\surd$
    \end{tabular}
    }
    \caption{Smooth functions in simulation. We define the numerical CP/Tucker rank as the minimal rank $r$ for which the relative approximation error is below $10^{-4}$. The reported rank in the table is estimated from a $100\times100\times100$ signal tensor generated by \eqref{eq:rep}.}
    \label{tb:md}
\end{table}
The generative functions involve compositions of operations such as polynomial, logarithm, exponential, square roots, etc. Notice that considered functions cover a reasonable range of model complexities from low rank to high rank. Two types of noise are considered: Gaussian noise and Bernoulli noise.
For the Gaussian model, we simulate continuous-valued tensors with i.i.d.\ noises drawn from $N(0,0.5^2)$. For the Bernoulli model, we generate binary tensors $\tY$ using the success probability tensor $\Theta \circ \pi$. The permutation $\pi$ is randomly chosen. For space consideration, only results for Models 1, 3, and 5 are presented in the main paper. The rest is presented in the Appendix A.1. We first examine impacts of model complexity to estimation accuracy. We then compare Borda count estimation with alternative methods under a range of scenarios. Extensions to non-symmetric tensors and extra simulation results are provided in the Appendix A.

\begin{figure}[t!]
    \centering
    \begin{subfigure}[b]{.85\textwidth}
    \includegraphics[width = \textwidth]{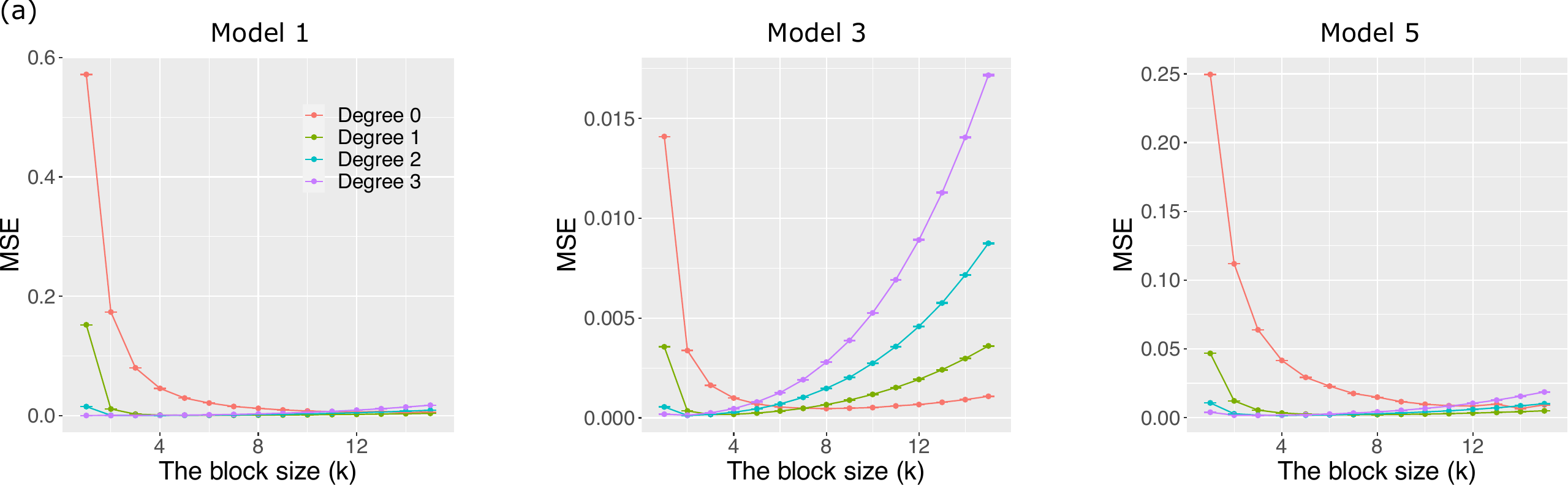}    
    \end{subfigure}
    
    \vspace{.5cm}
    
    \begin{subfigure}[b]{.85\textwidth}
    \includegraphics[width = \textwidth]{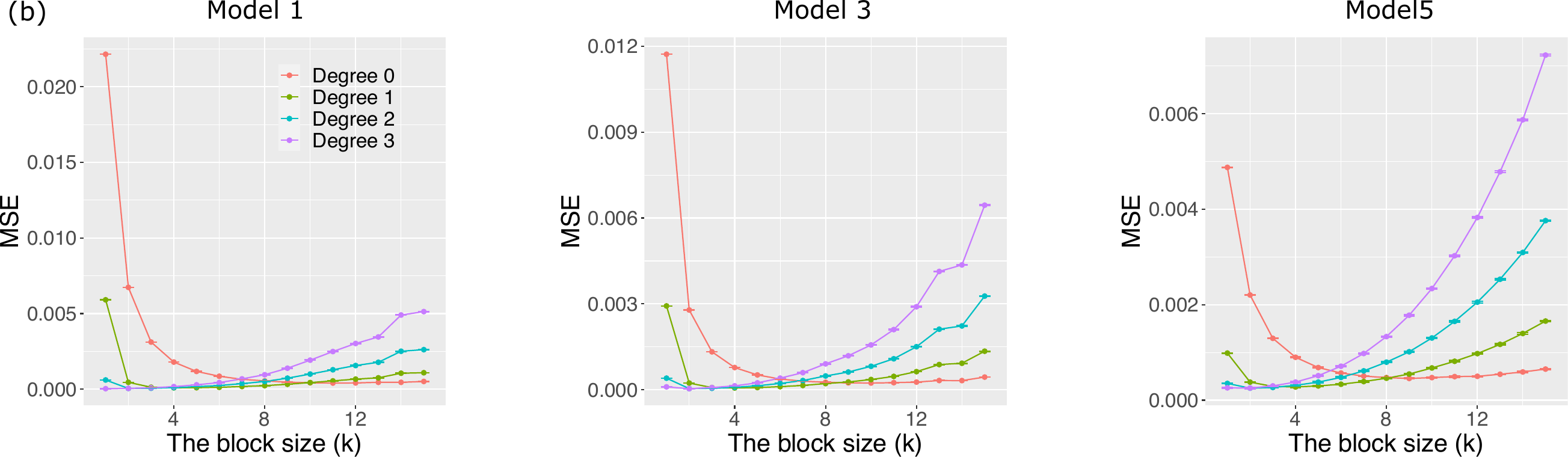}    
    \end{subfigure}
    \caption{MSE versus the number of blocks based on different polynomial approximations. Columns 1-3 consider the Models 1, 3, and 5 respectively. Panel (a) is for continuous tensors, whereas (b) is for the binary tensors.}
    \label{fig:degk}
\end{figure}

\paragraph{Impacts of the number of blocks, tensor dimension, and polynomial degree.}
The first experiment examines the impact of the block number $k$ and degree of polynomial $\ell$ for the approximation. We fix the tensor dimension $d = 100$, and vary the number of blocks $k\in\{1,\ldots,15\}$ and polynomial degree $\ell\in\{0,1,2,3\}.$
Figure~\ref{fig:degk} demonstrates the trade-off in accuracy determined by the number of groups for each polynomial degree. The results confirm our bias-variance analysis in Theorem~\ref{thm:LSE}. While a large block number $k$ provides less biased approximation, this large $k$ renders the signal tensor estimation difficult within each block due to a small sample size. In addition, we find that degree-2 polynomial approximation gives the smallest MSE among considered polynomial approximations for models 1-3. By Remark~\ref{rm2}, plugging $(\alpha,m)=(\infty,3)$ in~\eqref{eq:choice} gives the theoretical choice $(k^*,\ell^*)=(d^{7/3},2)$. The results are consistent with our simulation.

The second experiment investigates the impact of the tensor dimension $d$ for various polynomial degrees. We vary the tensor dimension $d\in\{10,\ldots,100\}$ and polynomial degree $\ell\in\{0,1,2,3\}$ in each model configuration. We set the optimal number of blocks to the configuration that gives the best accuracy. Figure~\ref{fig:degdim} compares the estimation errors among different polynomial approximations. The result verifies that the degree-2 polynomial approximation performs the best under the sufficient tensor dimension, which is consistent with our theoretical results. We emphasize that this phenomenon is different from the matrix case where the degree-0 polynomial approximation gives the best results~\citep{gao2015rate,klopp2017oracle}.

\paragraph{Comparison with alternative methods.} We compare our method ({\bf \small Borda count}) with several popular alternative methods. 
\begin{itemize}[wide, labelwidth=0pt, labelindent=0pt,topsep=0pt,itemsep=-1ex]
     \item SVD ({\bf \small Spectral}) \citep{xu2018rates} performs universal singular-value thresholding~\citep{chatterjee2015matrix} on unfolded tensors.
     \item Least-squares estimation ({\bf \small LSE}) \citep{han2022exact} uses spectral $k$-means to approximately solve the optimization \eqref{eq:lseopt} with constant block approximation ($\ell=0$). 
     \item Least-squares estimation ({\bf \small BAL}) \citep{balasubramanian2021nonparametric} uses count-based statistics to approximately solve the optimization~\eqref{eq:lseopt} with constant block approximation ($\ell=0$). This algorithm is only available for binary observations, so we only use it for comparison under the Bernoulli model.
 \end{itemize}
 
 We choose degree-2 polynomial approximation as our theorems suggested, and vary tensor dimension $d\in\{10,\ldots,100\}$ under each model configuration. Possible hyperparameters are the block number for {\bf \small Borda count}, {\bf \small LSE} and {\bf \small BAL}, and the singular-value threshold for {\bf \small Spectral}. Hyperparametes are set to achieve the best performance in the outputs.

\begin{figure}[t!]
    \centering
    \begin{subfigure}[b]{.85\textwidth}
    \includegraphics[width = \textwidth]{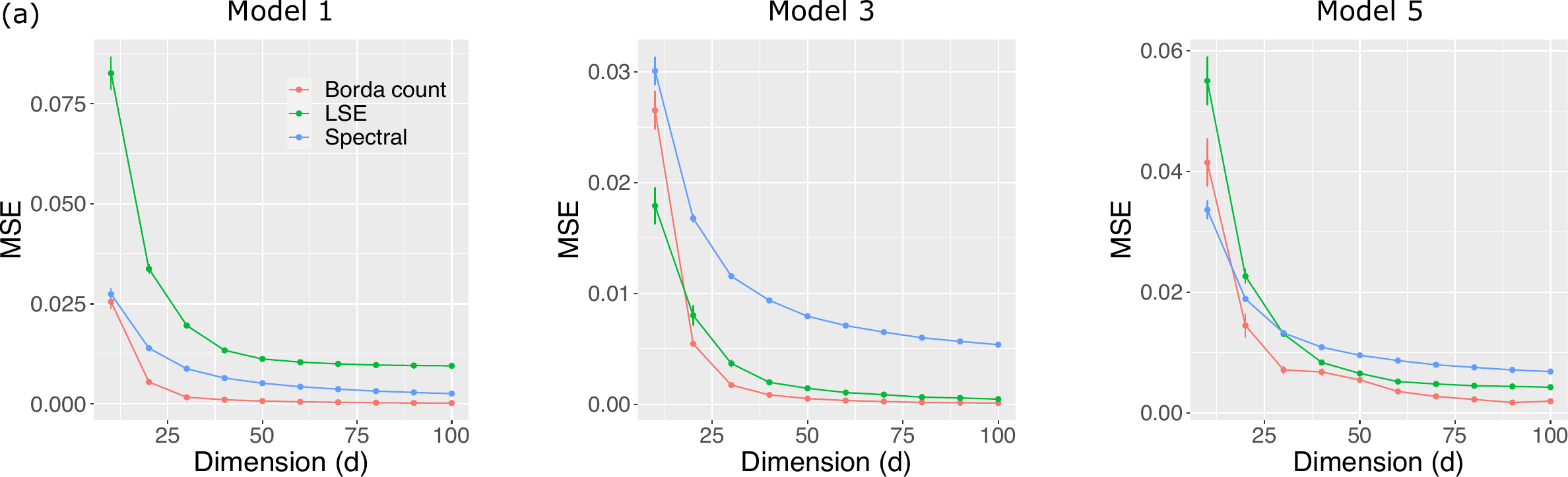}
    \vspace{-.6cm}
    \end{subfigure}
    \begin{subfigure}[b]{.85\textwidth}
    \includegraphics[width = \textwidth]{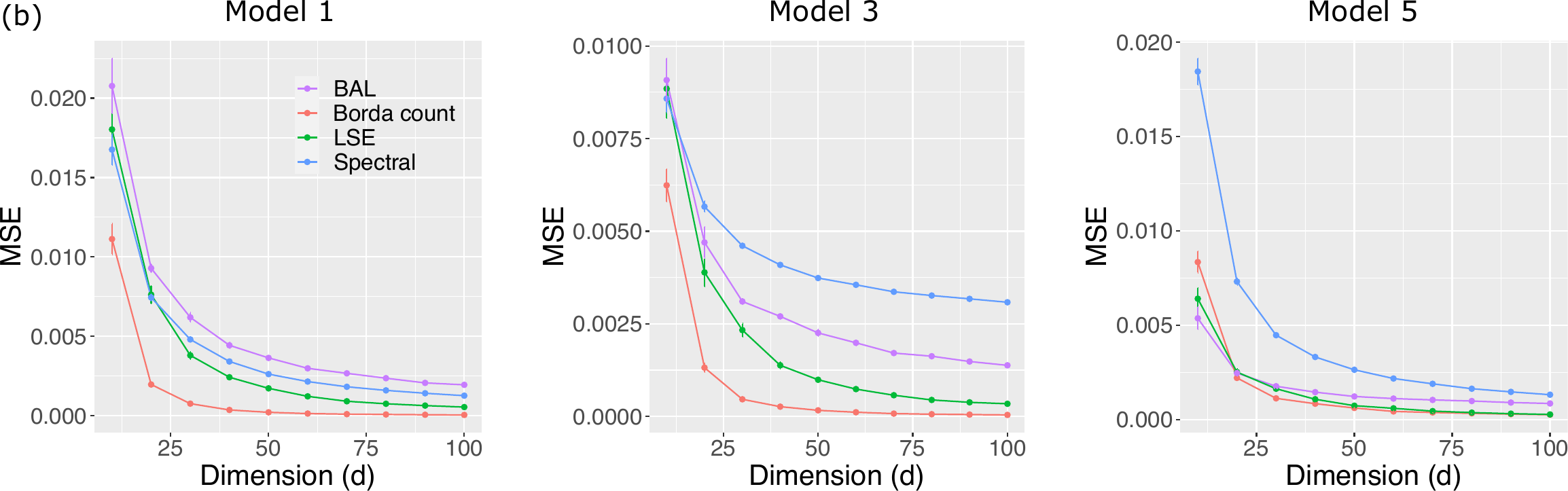}   
    \end{subfigure}
    \caption{MSE versus the tensor dimension based on different estimation methods. Columns 1-3 consider the Models 1, 3, and 5 in Table~\ref{tb:md} respectively. Panel (a) is for continuous tensors, whereas (b) is for the binary tensors.}
    \label{fig:method}
\end{figure}

\begin{figure}[t!]
    \centering
    \includegraphics[width =0.8\textwidth]{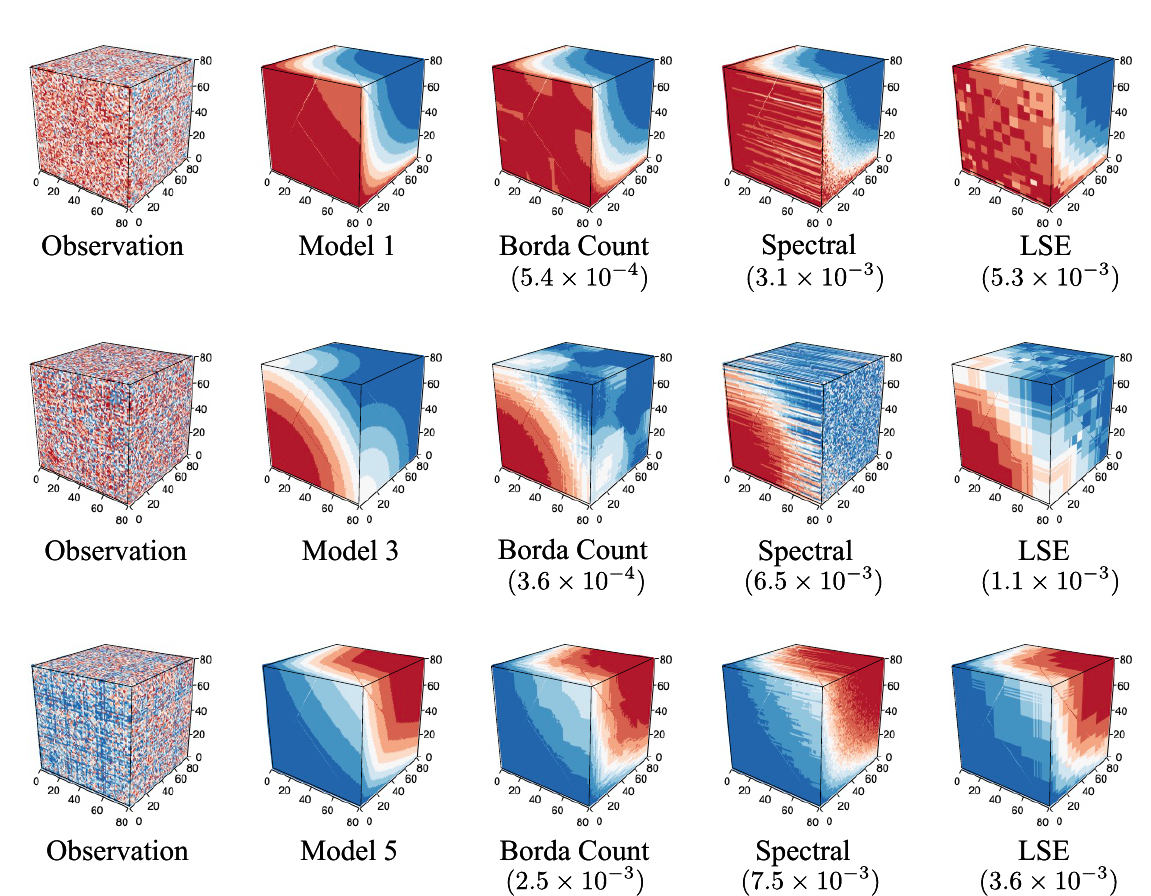}
    \caption{Performance comparison among different methods. The observed data tensors, true signal tensors, and estimated signal tensors are plotted for Models 1, 3 and 5 in Table~\ref{tb:md} with fixed dimension $d = 80$. Numbers in parenthesis indicate the mean squared error.}
    \label{fig:contim}
\end{figure}

Figure~\ref{fig:method} shows that our algorithm {\bf \small Borda count} achieves the best performance in all scenarios as the tensor dimension increases. The poor performance of {\bf \small Spectral} can be explained by the loss of multilinear structure in the tensor unfolding procedure. The sub-optimality of {\bf \small LSE} and {\bf \small BAL} is possibly due to its limits in both statistics and computations. Statistically, constant block approximation results in sub-optimal rates compared to polynomial approximation. Computationally, the least-squares optimization~\eqref{eq:lseopt} is computationally unstable. Figure~\ref{fig:contim} displays true signal tensors of three models and corresponding observed tensors of dimension $d = 80$ with Gaussian noise. We use oracle permutation $\pi$ to obtain the estimated signal tensor from the estimated permuted signal tensor $\hat\Theta\circ\hat\pi$ for the better visualizations. We see that our  {\bf \small Borda count} achieves the best signal recovery, thereby supporting the numerical results in Figure~\ref{fig:method}.

\subsection{Applications to Chicago crime data}
The Chicago crime dataset consists of crime counts reported in the city of Chicago, ranging from January 1st, 2001 to December 11th, 2017. The observed tensor is an order-3 tensor with entries representing the log counts of crimes from 24 hours, 77 community areas, and 32 crime types. We apply our Borda count method to Chicago crime dataset. Because the data tensor is non-symmetric, we allow different number of blocks across the three modes. Cross-validation results suggest the optimal $(k_1,k_2,k_3)=(6,4,10)$, representing the block number for crime hours, community areas, and crime types, respectively.

\begin{figure}[t!]
    \centering
    \includegraphics[width = 0.8\textwidth]{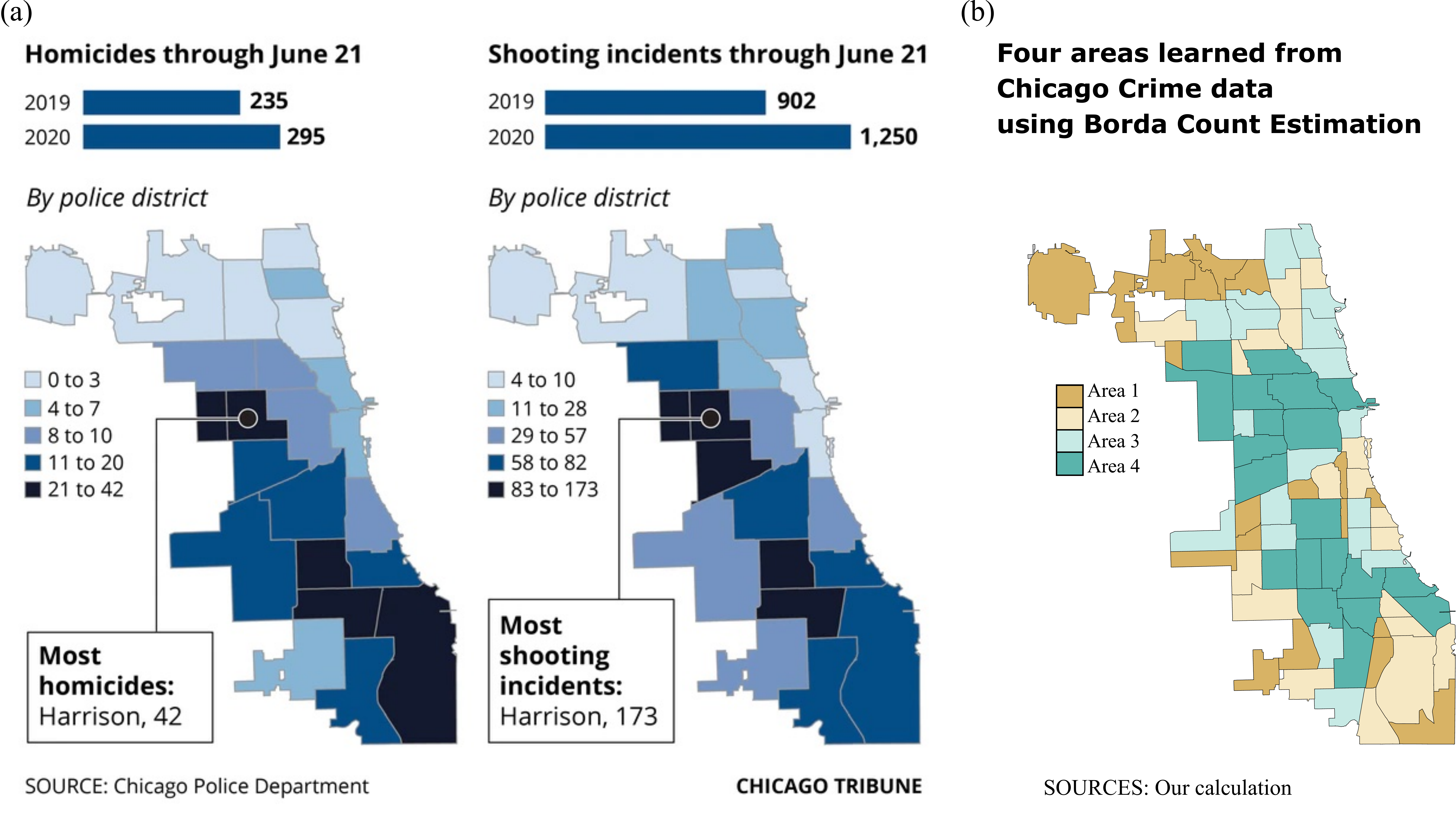}
    \caption{Chicago crime maps. Figure(a) is the benchmark map based on homicides and shooting incidents in community areas in Chicago~\citep{Jeremy.2020}. Figure(b) shows the four clustered areas learned from 32 crime types using our method.}
    \label{fig:area}
\end{figure}
We first investigate the four community areas obtained from Borda count algorithm.  Figure~\ref{fig:area}(b) shows the four areas overlaid on the Chicago map. Interestingly, we find that the clusters are consistent with actual locations, even though our algorithm did not take any geographic information such as longitude or latitude as inputs. In addition, we compare the cluster patterns with benchmark maps based on homicides and shooting incidents in Chicago shown in Figure~\ref{fig:area}(a). Our clusters share similar geographical patterns with Figure~\ref{fig:area}(a). The results demonstrate the power of our approach in detecting patterns from tensor data.

\begin{figure}[t!]
    \centering
    \includegraphics[width = \textwidth]{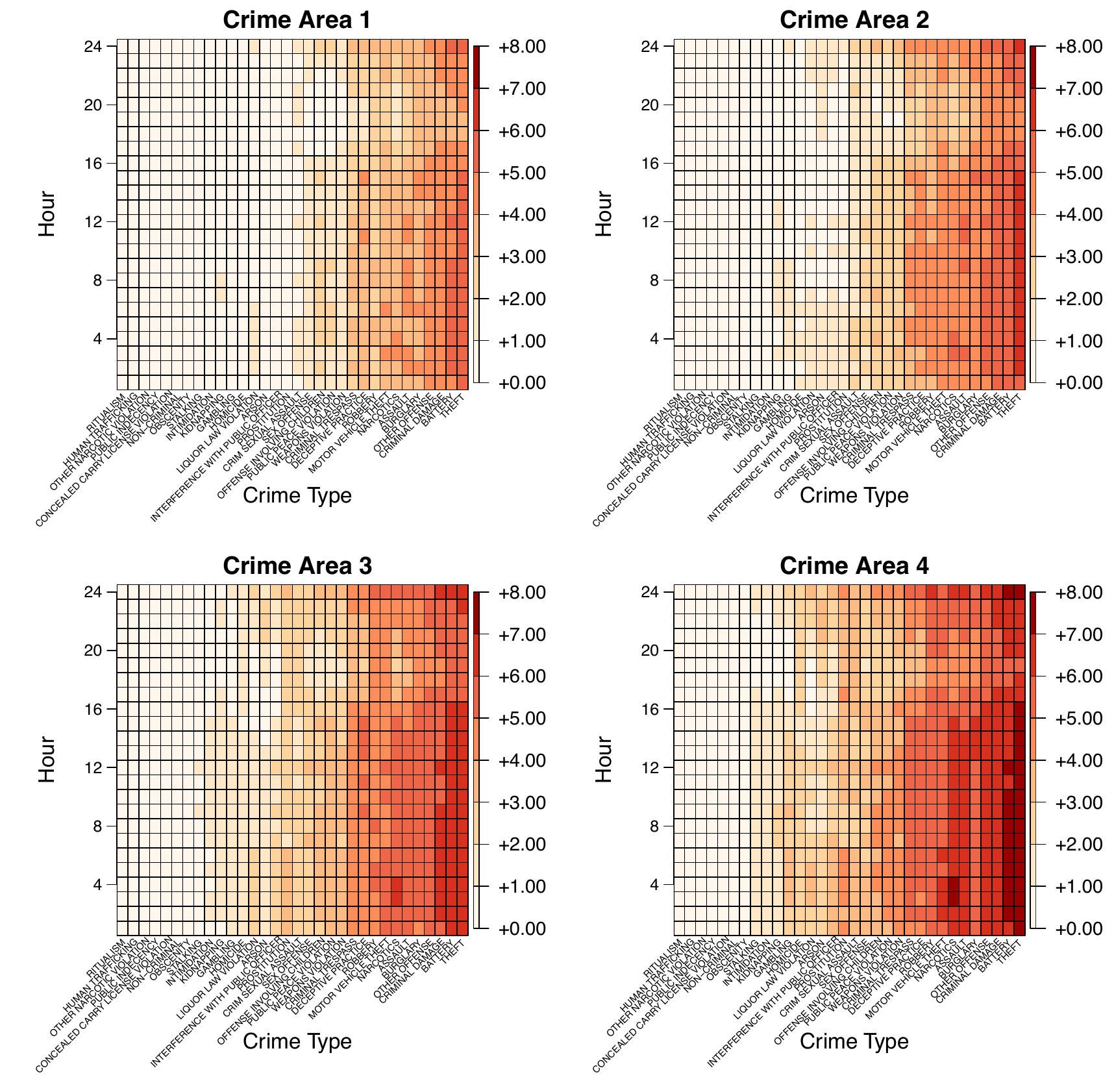}
    \caption{Averaged log counts of crimes according to crime types, hours, and the four areas estimated by our Borda count algorithm. We plot the estimated signal tensor entries averaged within four areas in the heatmap.}
    \label{fig:crimeA}
\end{figure}

Then, we examine the denoised signal tensor obtained from our method and analyze the trends between crime types and crime hours by the four areas in Figure~\ref{fig:area}(b). Figure~\ref{fig:crimeA} shows the averaged log counts of crimes according to crime types and crime hours by four areas. We find that the major difference among four areas is the crime rates. Area 4 has the highest crime rates,  and the crime rates monotonically decrease from Area 4 to Area 1. The variation in crime rates across hour and type, nevertheless, exhibits similarity among the four areas. Figure~\ref{fig:crimeA} shows that the number of crimes increases hourly from 8 p.m., peaks at night hours, and then drops to the lowest at 6 p.m. The identified pattens among the four community areas highlight the applicability of our method in real data.

Finally, we compare the prediction performance based on constant block model and our permuted smooth tensor model. Notice that constant block model uses $\ell=0$ approximation, whereas our permuted smooth tensor model uses $\ell=2$ approximation. We found that the mean squared prediction errors for our model vs.\ constant block model are 0.283 (0.006) vs.\ 0.399 (0.009), respectively. Here, the reported prediction errors are averaged over five runs of cross-validation, with standard errors in parentheses. The block number $(k_1,k_2,k_3)$ with best prediction performance is $(6,4,10)$ for our models, and $(7,11,10)$ for constant block models. We see that the permuted smooth tensor model substantially outperforms the classical constant block models.

\section{Conclusion and Discussions}\label{sec:discussion}
We have presented a suite of statistical theory, estimation methods, and data applications for permuted smooth tensor models. 
 Two estimation algorithms are proposed with accuracy guarantees: the (statistically optimal) least-squares estimation and the (computationally tractable) Borda count estimation.
In particular, we establish an interesting phase transition phenomenon with respect to the critical smoothness level. We demonstrate that a block-wise polynomial of order $(m-2)(m+1)/2$ is sufficient and necessary for accurate recovery of order-$m$ tensors, in contrast to earlier beliefs on constant block approximation. Experiments demonstrate the effectiveness of both theoretical findings and algorithms. 

We discuss several possible extensions from our work.

One limitation of our model is that we consider the 1-dimensional latent space embedding. The extension from the 1-dimensional latent space model to general dimensional latent model is analogous to the extension from the block model to the mixed membership model. Our parallel work~\citep{lee2023statistical} considers the general dimensional latent variable model, by assuming  a set of $s$-dimensional vectors $\ma^{(k)}_{i_k}\in\mathbb{R}^s$ with $s\geq 1$ and a latent function $f\colon[0,1]^{s}\times\cdots\times[0,1]^{s}\rightarrow \mathbb{R}$ such that
    \begin{align}\label{eq:LVM2}
    \Theta(i_1,\ldots,i_m) = f(\ma^{(1)}_{i_1},\ldots,\ma^{(m)}_{i_m}), \text{ for all } (i_1,\ldots,i_m)\in[d_1]\times\cdots\times [d_m].
\end{align}
This generalization extends the latent permutation $\pi\in\Pi(d,d)$ to the set of latent vectors in $[0,1]^s$. However, we find that this extension is not free. We need a stronger analytic function class with $\infty$-smoothness for the theoretical analysis. Compared to this current paper, the analysis of analytic functions uses different techniques and yields new results of its own. We refer readers to~\cite{lee2023statistical} for independent interests. 

Another limitation of our algorithm is the need for hyperparameter tuning. There is a vast literature on nonparametric estimation that focuses on adaptivity. For example, spatially adaptive methods have been developed in the contexts of wavelets~\citep{donoho1994ideal}, splines~\citep{mammen1997locally}, and trend filtering~\citep{tibshirani2014adaptive}; tuning-free algorithms have been proposed for several shape-constrained functions~\citep{chatterjee2019adaptive,feng2022nonparametric,bellec2018sharp}; see~\cite{cai2012minimax} for a review. Our work is orthogonal to these advances, and in principle we can combine these tools in our tensor estimation. In this paper, we choose the standard polynomial algorithm because of its simplicity. The parsimony leads to easier analysis of the critical smoothness level $ (m-2)(m+1)/2$. Exploiting various nonparametric techniques for tensor models warrants future research.

\section*{Acknowledgements}
This research is supported in part by NSF CAREER DMS-2141865, DMS-1915978, DMS-2023239, EF-2133740, and funding from the Wisconsin Alumni Research foundation. 

\section*{Supplementary Materials}
The supplementary materials include all proofs and extra simulation results. The R package for implementing the methods described in this article is released at \url{https://github.com/Chanwoost/Smooth-tensor-estimation-with-unknown-permutations}

\bibliographystyle{Chicago.bst}
\bibliography{tensor_wang}

\begin{thebibliography}{}

\bibitem[\protect\citeauthoryear{Balasubramanian}{Balasubramanian}{2021}]{balasubramanian2021nonparametric}
Balasubramanian, K. (2021).
\newblock Nonparametric modeling of higher-order interactions via
  hypergraphons.
\newblock {\em Journal of Machine Learning Research\/}~{\em 22}, 1--25.

\bibitem[\protect\citeauthoryear{Baltrunas, Kaminskas, Ludwig, Moling, Ricci,
  Aydin, L{\"u}ke, and Schwaiger}{Baltrunas
  et~al.}{2011}]{baltrunas2011incarmusic}
Baltrunas, L., M.~Kaminskas, B.~Ludwig, O.~Moling, F.~Ricci, A.~Aydin, K.-H.
  L{\"u}ke, and R.~Schwaiger (2011).
\newblock {I}n{C}ar{M}usic: Context-aware music recommendations in a car.
\newblock In {\em International Conference on Electronic Commerce and Web
  Technologies}, pp.\  89--100. Springer.

\bibitem[\protect\citeauthoryear{Bellec}{Bellec}{2018}]{bellec2018sharp}
Bellec, P.~C. (2018).
\newblock Sharp oracle inequalities for least squares estimators in shape
  restricted regression.
\newblock {\em The Annals of Statistics\/}~{\em 46\/}(2), 745--780.

\bibitem[\protect\citeauthoryear{Bi, Qu, and Shen}{Bi
  et~al.}{2018}]{bi2018multilayer}
Bi, X., A.~Qu, and X.~Shen (2018).
\newblock Multilayer tensor factorization with applications to recommender
  systems.
\newblock {\em The Annals of Statistics\/}~{\em 46\/}(6B), 3308--3333.

\bibitem[\protect\citeauthoryear{Bickel and Chen}{Bickel and
  Chen}{2009}]{bickel2009nonparametric}
Bickel, P.~J. and A.~Chen (2009).
\newblock A nonparametric view of network models and newman--girvan and other
  modularities.
\newblock {\em Proceedings of the National Academy of Sciences\/}~{\em
  106\/}(50), 21068--21073.

\bibitem[\protect\citeauthoryear{Cai}{Cai}{2012}]{cai2012minimax}
Cai, T.~T. (2012).
\newblock Minimax and adaptive inference in nonparametric function estimation.
\newblock {\em Statistical Science\/}~{\em 27\/}(1), 31--50.

\bibitem[\protect\citeauthoryear{Chan and Airoldi}{Chan and
  Airoldi}{2014}]{chan2014consistent}
Chan, S. and E.~Airoldi (2014).
\newblock A consistent histogram estimator for exchangeable graph models.
\newblock In {\em International Conference on Machine Learning}, pp.\
  208--216.

\bibitem[\protect\citeauthoryear{Chatterjee}{Chatterjee}{2015}]{chatterjee2015matrix}
Chatterjee, S. (2015).
\newblock Matrix estimation by universal singular value thresholding.
\newblock {\em The Annals of Statistics\/}~{\em 43\/}(1), 177--214.

\bibitem[\protect\citeauthoryear{Chatterjee and Lafferty}{Chatterjee and
  Lafferty}{2019}]{chatterjee2019adaptive}
Chatterjee, S. and J.~Lafferty (2019).
\newblock Adaptive risk bounds in unimodal regression.
\newblock {\em Bernoulli\/}~{\em 25\/}(1), 1--25.

\bibitem[\protect\citeauthoryear{Ding, Ma, Wu, and Xu}{Ding
  et~al.}{2021}]{ding2021efficient}
Ding, J., Z.~Ma, Y.~Wu, and J.~Xu (2021).
\newblock Efficient random graph matching via degree profiles.
\newblock {\em Probability Theory and Related Fields\/}~{\em 179\/}(1),
  29--115.

\bibitem[\protect\citeauthoryear{Donoho and Johnstone}{Donoho and
  Johnstone}{1994}]{donoho1994ideal}
Donoho, D.~L. and I.~M. Johnstone (1994).
\newblock Ideal spatial adaptation by wavelet shrinkage.
\newblock {\em biometrika\/}~{\em 81\/}(3), 425--455.

\bibitem[\protect\citeauthoryear{Feng, Chen, Han, Carroll, and Samworth}{Feng
  et~al.}{2022}]{feng2022nonparametric}
Feng, O.~Y., Y.~Chen, Q.~Han, R.~J. Carroll, and R.~J. Samworth (2022).
\newblock Nonparametric, tuning-free estimation of s-shaped functions.
\newblock {\em Journal of the Royal Statistical Society Series B: Statistical
  Methodology\/}~{\em 84\/}(4), 1324--1352.

\bibitem[\protect\citeauthoryear{Flammarion, Mao, and Rigollet}{Flammarion
  et~al.}{2019}]{flammarion2019optimal}
Flammarion, N., C.~Mao, and P.~Rigollet (2019).
\newblock Optimal rates of statistical seriation.
\newblock {\em Bernoulli\/}~{\em 25\/}(1), 623--653.

\bibitem[\protect\citeauthoryear{Gao, Lu, Ma, and Zhou}{Gao
  et~al.}{2016}]{gao2016optimal}
Gao, C., Y.~Lu, Z.~Ma, and H.~H. Zhou (2016).
\newblock Optimal estimation and completion of matrices with biclustering
  structures.
\newblock {\em Journal of Machine Learning Research\/}~{\em 17\/}(1),
  5602--5630.

\bibitem[\protect\citeauthoryear{Gao, Lu, and Zhou}{Gao
  et~al.}{2015}]{gao2015rate}
Gao, C., Y.~Lu, and H.~H. Zhou (2015).
\newblock Rate-optimal graphon estimation.
\newblock {\em The Annals of Statistics\/}~{\em 43\/}(6), 2624--2652.

\bibitem[\protect\citeauthoryear{Gy{\"o}rfi, Kohler, Krzy{\.z}ak, and
  Walk}{Gy{\"o}rfi et~al.}{2002}]{gyorfi2002distribution}
Gy{\"o}rfi, L., M.~Kohler, A.~Krzy{\.z}ak, and H.~Walk (2002).
\newblock {\em A distribution-free theory of nonparametric regression},
  Volume~1.
\newblock Springer.

\bibitem[\protect\citeauthoryear{Han, Luo, Wang, and Zhang}{Han
  et~al.}{2022}]{han2022exact}
Han, R., Y.~Luo, M.~Wang, and A.~R. Zhang (2022).
\newblock Exact clustering in tensor block model: Statistical optimality and
  computational limit.
\newblock {\em Journal of the Royal Statistical Society Series B: Statistical
  Methodology\/}~{\em 84\/}(5), 1666--1698.

\bibitem[\protect\citeauthoryear{H{\"u}tter, Mao, Rigollet, and
  Robeva}{H{\"u}tter et~al.}{2020}]{hutter2020estimation}
H{\"u}tter, J.-C., C.~Mao, P.~Rigollet, and E.~Robeva (2020).
\newblock Estimation of monge matrices.
\newblock {\em Bernoulli\/}~{\em 26\/}(4), 3051--3080.

\bibitem[\protect\citeauthoryear{Jeremy}{Jeremy}{2020}]{Jeremy.2020}
Jeremy, G. (2020).
\newblock A trying first half of 2020 included spike in shootings and homicides
  in chicago.
\newblock {\em Chicago Tribune, 2020-06-26\/}.

\bibitem[\protect\citeauthoryear{Klopp, Tsybakov, and Verzelen}{Klopp
  et~al.}{2017}]{klopp2017oracle}
Klopp, O., A.~B. Tsybakov, and N.~Verzelen (2017).
\newblock Oracle inequalities for network models and sparse graphon estimation.
\newblock {\em The Annals of Statistics\/}~{\em 45\/}(1), 316--354.

\bibitem[\protect\citeauthoryear{Kolda and Bader}{Kolda and
  Bader}{2009}]{kolda2009tensor}
Kolda, T.~G. and B.~W. Bader (2009).
\newblock Tensor decompositions and applications.
\newblock {\em SIAM Review\/}~{\em 51\/}(3), 455--500.

\bibitem[\protect\citeauthoryear{Lee and Wang}{Lee and
  Wang}{2023}]{lee2023statistical}
Lee, C. and M.~Wang (2023).
\newblock Statistical and computational rates in high rank tensor estimation.
\newblock {\em arXiv preprint arXiv:2304.04043\/}.

\bibitem[\protect\citeauthoryear{Li, Shah, Song, and Yu}{Li
  et~al.}{2019}]{li2019nearest}
Li, Y., D.~Shah, D.~Song, and C.~L. Yu (2019).
\newblock Nearest neighbors for matrix estimation interpreted as blind
  regression for latent variable model.
\newblock {\em IEEE Transactions on Information Theory\/}~{\em 66\/}(3),
  1760--1784.

\bibitem[\protect\citeauthoryear{Livi and Rizzi}{Livi and
  Rizzi}{2013}]{livi2013graph}
Livi, L. and A.~Rizzi (2013).
\newblock The graph matching problem.
\newblock {\em Pattern Analysis and Applications\/}~{\em 16\/}(3), 253--283.

\bibitem[\protect\citeauthoryear{Lov{\'a}sz}{Lov{\'a}sz}{2012}]{lovasz2012large}
Lov{\'a}sz, L. (2012).
\newblock {\em Large networks and graph limits}, Volume~60.
\newblock American Mathematical Soc.

\bibitem[\protect\citeauthoryear{Luo and Zhang}{Luo and
  Zhang}{2022}]{luo2022tensor}
Luo, Y. and A.~R. Zhang (2022).
\newblock Tensor clustering with planted structures: Statistical optimality and
  computational limits.
\newblock {\em The Annals of Statistics\/}~{\em 50\/}(1), 584--613.

\bibitem[\protect\citeauthoryear{Mammen and Van De~Geer}{Mammen and Van
  De~Geer}{1997}]{mammen1997locally}
Mammen, E. and S.~Van De~Geer (1997).
\newblock Locally adaptive regression splines.
\newblock {\em The Annals of Statistics\/}~{\em 25\/}(1), 387--413.

\bibitem[\protect\citeauthoryear{Pananjady and Samworth}{Pananjady and
  Samworth}{2022}]{pananjady2022isotonic}
Pananjady, A. and R.~J. Samworth (2022).
\newblock Isotonic regression with unknown permutations: Statistics,
  computation and adaptation.
\newblock {\em The Annals of Statistics\/}~{\em 50\/}(1), 324--350.

\bibitem[\protect\citeauthoryear{Rigollet and H{\"u}tter}{Rigollet and
  H{\"u}tter}{2015}]{rigollet2015high}
Rigollet, P. and J.-C. H{\"u}tter (2015).
\newblock High dimensional statistics.
\newblock {\em Lecture notes for course 18S997\/}~{\em 813\/}(814), 46.

\bibitem[\protect\citeauthoryear{Shah, Balakrishnan, and Wainwright}{Shah
  et~al.}{2019}]{shah2019low}
Shah, N., S.~Balakrishnan, and M.~Wainwright (2019).
\newblock Low permutation-rank matrices: Structural properties and noisy
  completion.
\newblock {\em Journal of Machine Learning Research\/}~{\em 20}, 1--43.

\bibitem[\protect\citeauthoryear{Stone}{Stone}{1982}]{stone1982optimal}
Stone, C.~J. (1982).
\newblock Optimal global rates of convergence for nonparametric regression.
\newblock {\em The Annals of Statistics\/}~{\em 10\/}(4), 1040--1053.

\bibitem[\protect\citeauthoryear{Tibshirani}{Tibshirani}{2014}]{tibshirani2014adaptive}
Tibshirani, R.~J. (2014).
\newblock Adaptive piecewise polynomial estimation via trend filtering.
\newblock {\em The Annals of Statistics\/}~{\em 42\/}(1), 285--323.

\bibitem[\protect\citeauthoryear{Tsybakov}{Tsybakov}{2009}]{tsybakov2009introduction}
Tsybakov, A.~B. (2009).
\newblock {\em Introduction to nonparametric estimation}.
\newblock Springer Science \& Business Media.

\bibitem[\protect\citeauthoryear{Wang, Fischer, and Song}{Wang
  et~al.}{2019}]{wang2019three}
Wang, M., J.~Fischer, and Y.~S. Song (2019).
\newblock Three-way clustering of multi-tissue multi-individual gene expression
  data using semi-nonnegative tensor decomposition.
\newblock {\em The Annals of Applied Statistics\/}~{\em 13\/}(2), 1103--1127.

\bibitem[\protect\citeauthoryear{Wang and Li}{Wang and
  Li}{2020}]{wang2018learning}
Wang, M. and L.~Li (2020).
\newblock Learning from binary multiway data: Probabilistic tensor
  decomposition and its statistical optimality.
\newblock {\em Journal of Machine Learning Research\/}~{\em 21\/}(154), 1--38.

\bibitem[\protect\citeauthoryear{Wang and Zeng}{Wang and
  Zeng}{2019}]{wang2019multiway}
Wang, M. and Y.~Zeng (2019).
\newblock Multiway clustering via tensor block models.
\newblock In {\em Advances in Neural Information Processing Systems}, pp.\
  713--723.

\bibitem[\protect\citeauthoryear{Wasserman}{Wasserman}{2006}]{wasserman2006all}
Wasserman, L. (2006).
\newblock {\em All of nonparametric statistics}.
\newblock Springer Science \& Business Media.

\bibitem[\protect\citeauthoryear{Xu}{Xu}{2018}]{xu2018rates}
Xu, J. (2018).
\newblock Rates of convergence of spectral methods for graphon estimation.
\newblock In {\em International Conference on Machine Learning}, pp.\
  5433--5442.

\bibitem[\protect\citeauthoryear{Zhang and Xia}{Zhang and
  Xia}{2018}]{zhang2018tensor}
Zhang, A. and D.~Xia (2018).
\newblock Tensor {SVD}: Statistical and computational limits.
\newblock {\em IEEE Transactions on Information Theory\/}~{\em 64\/}(11),
  7311--7338.

\bibitem[\protect\citeauthoryear{Zhao}{Zhao}{2015}]{zhao2015hypergraph}
Zhao, Y. (2015).
\newblock Hypergraph limits: a regularity approach.
\newblock {\em Random Structures \& Algorithms\/}~{\em 47\/}(2), 205--226.

\end{thebibliography}

\clearpage
\appendix
\section*{Appendix}
The appendix includes extra simulation results, proofs, and lemmas.

\renewcommand{\thefigure}{S\arabic{figure}}
\setcounter{figure}{0}   
\renewcommand{\thetable}{S\arabic{table}}
\setcounter{table}{0}   

\section{Extra numerical results}\label{app:extra}
\subsection{Results for Models 2 and 4 in Table~\ref{tab:comp}}\label{sec:extra}
\begin{figure}[htp!]
    \centering
    \includegraphics[width =0.9\textwidth]{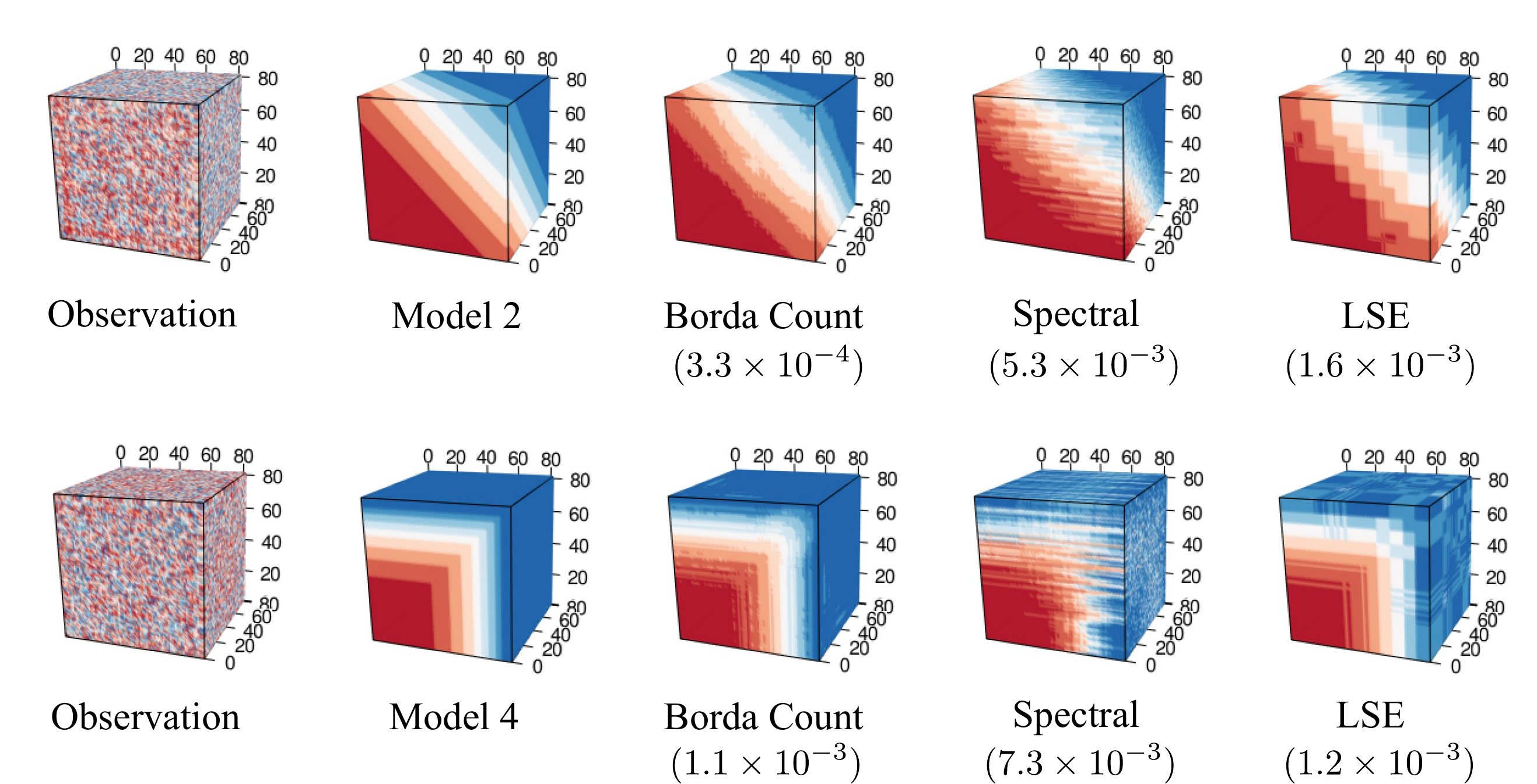}
    \caption{Performance comparison between different methods. The observed data tensors, true signal tensors, and estimated signal tensors are plotted for Models 2 and 4 in Table~\ref{tb:md} with fixed dimension $d = 80$. Numbers in parenthesis indicate the mean squared error.}
    \label{fig:extrav}
\end{figure}

\begin{figure}[htp!]
    \centering
    \begin{subfigure}[b]{\textwidth}
    \includegraphics[width = \textwidth]{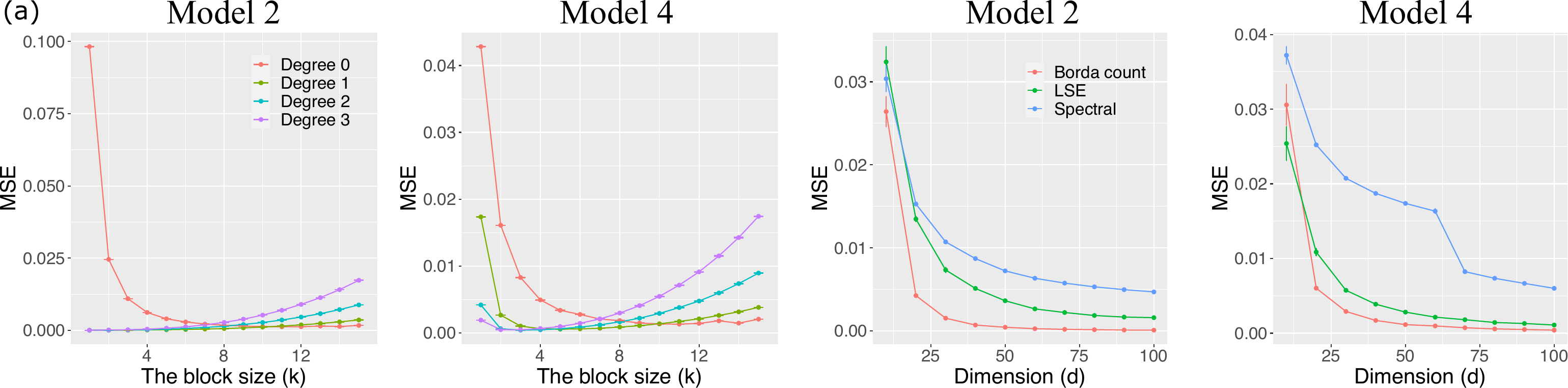}  
        \vspace{0cm}
    \end{subfigure}
    \begin{subfigure}[b]{\textwidth}
    \includegraphics[width = \textwidth]{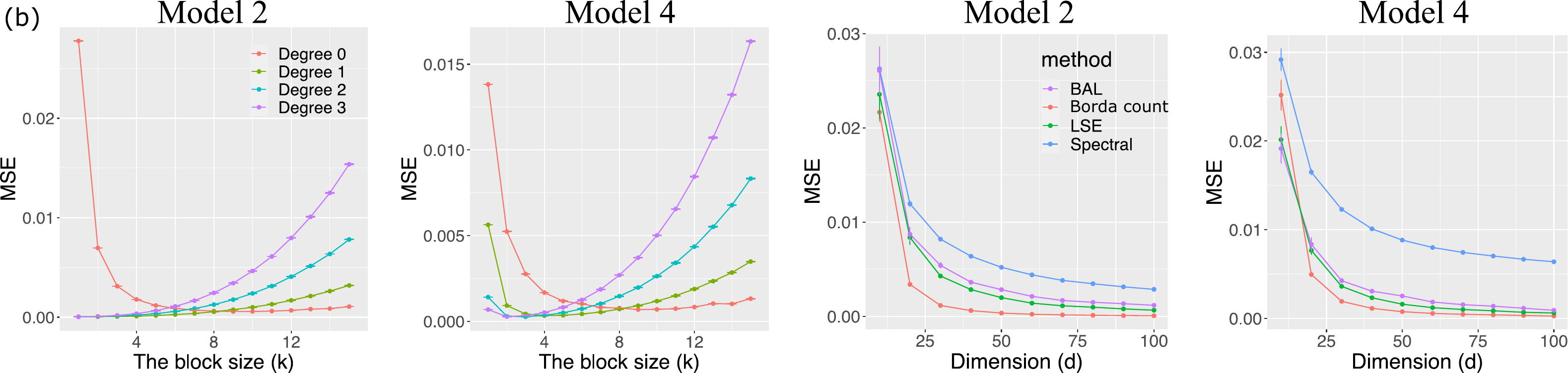}    
    \end{subfigure}
    \caption{Simulation results for Models 2 and 4 in Table~\ref{tb:md}. Columns 1-2 plots MSE versus the number of blocks for different polynomial approximation, while Columns 3-4 shows the MSE versus the tensor dimension according to estimation methods. Panel (a) is for continuous tensors, whereas (b) is for the binary tensors.}
    \label{fig:extrasim1}
\end{figure}

We first present simulation results for Models 2 and 4 omitted in Section~\ref{sec:sim}. Figure~\ref{fig:extrav} compares the estimation performance among the {\bf \small Borda count, LSE}, and {\bf \small Spectral methods}. We find that our Borda count algorithm outperforms others in both models. The first two columns in Figure~\ref{fig:extrasim1} show the impact of the number of blocks $k$ and degree of polynomial $\ell$ for the approximation with fixed dimension $d = 100$. Similar to results for Models 1, 3 and 5 in the main paper, we find the optimal $k$ balances the trade-off between approximation error and signal tensor estimation error within each block. The last two columns compare our {\bf \small Borda count} with other alternative methods. We find that our method still outperforms {\bf \small LSE} and {\bf \small Spectral} in all scenarios under Models 2 and 4.

\begin{figure}[h!]
    \centering
    \begin{subfigure}[b]{\textwidth}
    \vspace{0.5cm}
    
    \includegraphics[width = \textwidth]{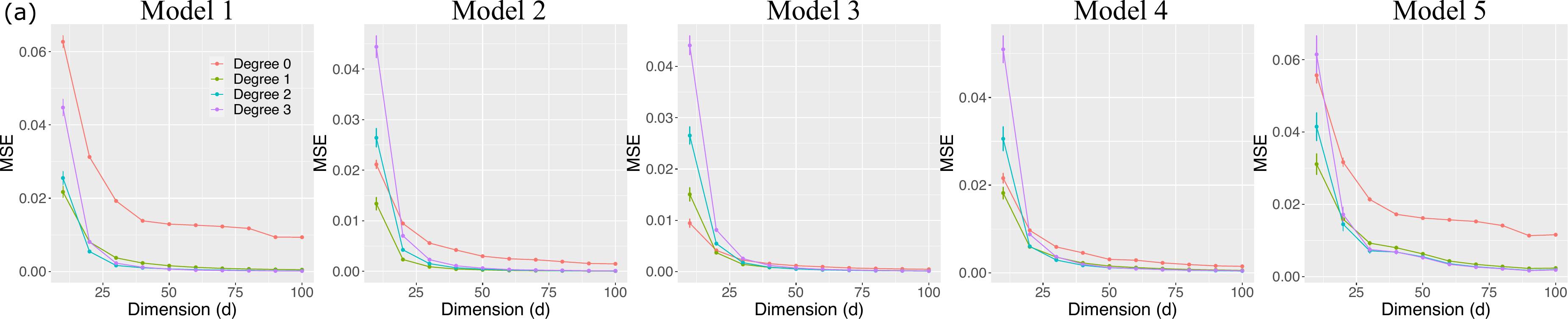} 
    \vspace{0cm}
    \end{subfigure}
    \begin{subfigure}[b]{\textwidth}
    \includegraphics[width =\textwidth]{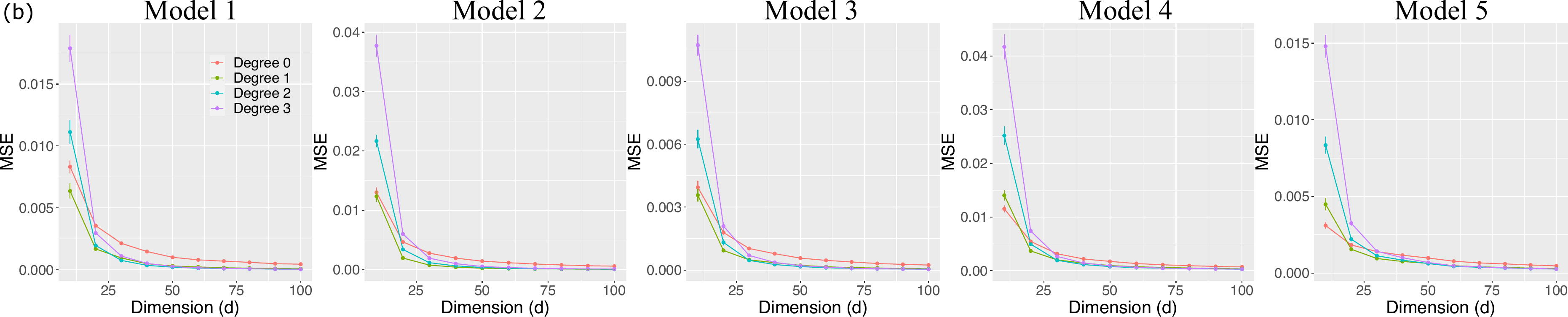}    
    \end{subfigure}
    \caption{MSE versus the tensor dimension based on different polynomial approximations. Columns 1-5 consider the Models 1-5 in Table~\ref{tb:md} respectively. Panel (a) is for continuous tensors, whereas (b) is for the binary tensors.}
    \label{fig:degdim}
\end{figure}

\begin{figure}[h]
    \centering
    \includegraphics[width = \textwidth]{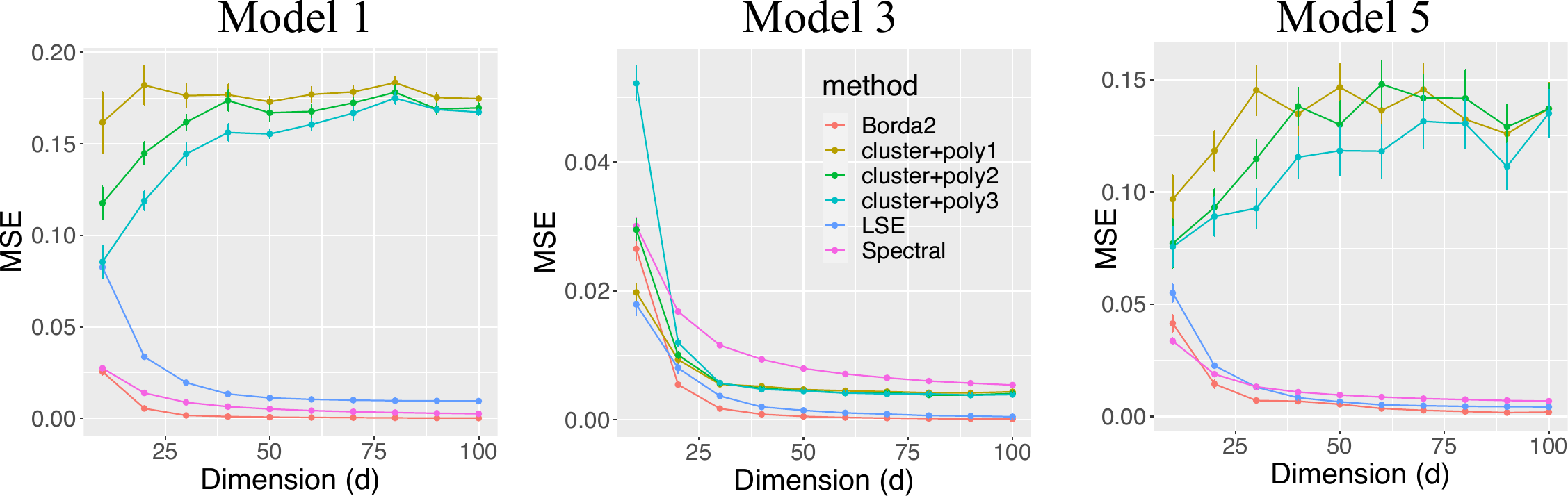}   
    \caption{MSE versus the tensor dimension based on different estimation methods. Columns 1-3 consider the Models 1, 3, and 5 with continuous case in Table~\ref{tb:md} respectively. cluster+poly$\ell$ means polynomial $\ell$-approximation}
    \label{fig:method2}
\end{figure}

Suggested by one reviewer, we consider another competing algorithm (denoted as cluster + poly$\ell$ algorithm) for comparison. The algorithm performs clustering first and then polynomial approximation within clusters. Since we do not know the permutation, we randomly assign the order of nodes within a block. We perform the simulation using the same setting in Section~\ref{sec:sim}. 
Figure~\ref{fig:method2} shows the comparison among all methods, including the new cluster + poly$\ell$ algorithm for $\ell\in\{1,2,3\}$. Notice that cluster + poly0 algorithm is equivalent to LSE algorithm. The result shows that the cluster + poly$\ell$ algorithms are unstable except Model 3. 

\subsection{Investigation of non-symmetric tensors}\label{subsec:asym}
Our models and techniques easily extend to non-symmetric tensors.
We describe the simulation set-up for non-symmetric tensors and results. We simulate order-3 tensors based on the non-symmetric functions in Table~\ref{tb:md2}. 

\begin{table}[ht]
    \centering
    \begin{tabular}{c|c}
        Model ID  &  $f(x,y,z)$  \\\hline
        1 &    $xy+z$ \\
        2&  $x^2+y +yz^2$\\
        3 & $x(1+\exp(-3(x^2+y^2+z^2)))^{-1}$\\
        4 & $\log(1+\max(x,y,z)+x^2+yz)$ \\
        5 &  $\exp\left(-x-\sqrt{y}-z^3\right)$
    \end{tabular}
    \caption{List of non-symmetric smooth functions in simulation.}
    \label{tb:md2}
\end{table}

We fix the tensor dimension ${30\times40\times 50}$ and assume that the noise tensors are from Gaussian distribution. Similar to other simulations, we evaluate the accuracy of the estimation by MSE and report the summary statistics across $n_{\text{sim}} = 20$ replicates. The hyperparameters are chosen via cross-validation that give the best accuracy for each method. Table~\ref{tb:hyper} summarizes the choice of hyperparameters.

\begin{table}[ht]
    \centering
    \begin{tabular}{c|c|c|c|c|c}
        Method &  Model 1 & Model 2 & Model 3& Model 4 & Model 5  \\\hline
        Borda count &   (2,1,2)&(1,2,2)& (1,3,3) &(2,1,2)&(1,4,4)\\
        LSE &(6,2,3)&(8,5,8)&(6,9,6)&(9,5,6)&(7,9,3)\\
        Spectral & (1,24)&(3,48)&(1,48)&(1,28)&(1,22)
    \end{tabular}
    \caption{Hyperparameters for the methods under Models 1-5 in Table~\ref{tb:md2}. For {\bf \small Borda count} and {\bf \small LSE} methods, the values in the table indicate the number of blocks. For {\bf \small Spectral} method, the first value indicates the tensor unfolding mode, while the second one represents the singular value threshold.}
    \label{tb:hyper}
\end{table}

 Table~\ref{tb:asymresult} compares the MSEs from repeated simulations based on different methods under Models 1-5 (see Table~\ref{tb:md2}). We find that Borda count estimation outperforms all alternative methods for non-symmetric tensors. 
The results demonstrate the applicability of our method to general tensors

\begin{table}[ht]
\centering
\resizebox{\textwidth}{!}{%
    \begin{tabular}{c|c|c|c|c|c}
        Method &  Model 1 & Model 2 & Model 3& Model 4 & Model 5  \\\hline
        Borda count &  {\bf 0.57 (0.01)}&{\bf 0.51 (0.02)}& {\bf0.87 (0.02)} &{\bf1.02 (0.02)}& {\bf2.56 (0.21)}\\
        LSE &23.58 (0.03)&7.70 (0.04)&9.45 (0.05)&3.29 (0.05)&9.93 (0.03)\\
        Spectral & 10.76 (0.06)&10.64 (0.05)&6.27 (0.05)&10.90 (0.06)&5.24 (0.04)
    \end{tabular}
    }
    \caption{MSEs from 20 repeated simulations based on different methods. All numbers are displayed on the scales $10^{-3}$. Standard errors are reported in parenthesis.}
    \label{tb:asymresult}
\end{table}

\subsection{Extra results on Chicago crime data analysis}\label{subsec:chicago}
We investigate the ten groups of crime types from our method. Table~\ref{tb:crimetb} shows the similar type of crimes captured by our clustering. We find that group 2 consists of misdemeanors such as public indecency, non-criminal, and concealed carry license violation, while group 6 represents sex-related offenses such as prostitution, sex offense, and crime sexual assault.

\begin{table}[ht!]
    \centering
    \includegraphics[width = \textwidth]{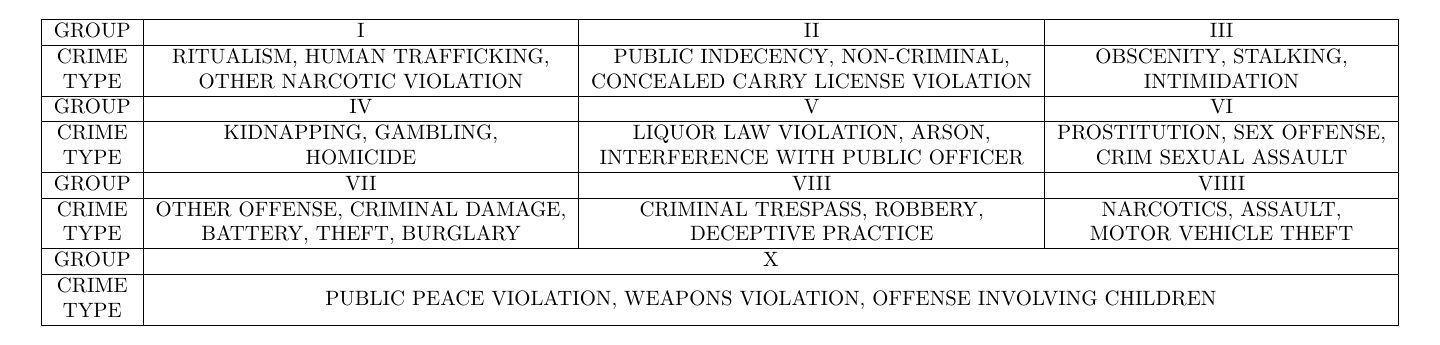}
    \caption{Groups of crime types learned based on the Borda count estimation.}
    \label{tb:crimetb}
\end{table}

\section{Details on Example~\ref{ex:free}}\label{sec:ex}
    We show that Borda count algorithm works in the Example~\ref{ex:free}. For simplicity, assume that $d$ is an even number. The signal tensor is represented by
    \begin{align}\label{eq:model6}
       \Theta(i,j,k) = \left(\frac{i}{d}-0.5\right)^2+\left(\frac{j}{d}\right)\left(\frac{k}{d}\right),\quad \text{for all }(i,j,k)\in[d]^3.
    \end{align}
    We construct permutations $\pi_1,\pi_2,\pi_3$ and a multivariate function $\bar f\colon [0,1]^3\rightarrow \mathbb{R}$ such that
    \begin{align}\label{eq:newrep}
        \Theta(i,j,k) = \bar f\left(\frac{\pi_1(i)}{d},\frac{\pi_2(j)}{d},\frac{\pi_3(k)}{d}\right)\pm{1\over d},\quad \text{for all }(i,j,k)\in[d]^3,
    \end{align}
    where the last term $1\over d$ is the approximation error. Specifically, the permutations in~\eqref{eq:newrep} are defined by
    \begin{align}\label{eq:newperm}
        \pi_1(i) = \begin{cases}
            2i-d, &\text{ if } i> \frac{d}{2},\\ d+1-2i, &\text{ if } i\leq \frac{d}{2},
        \end{cases}\quad\text{and}\quad  \pi_2(i) = \pi_3 (i)= i.
    \end{align}
Define the function $\bar f\colon [0,1]^3\rightarrow \mathbb{R}$ by 
\begin{align}\label{eq:newf}
\bar f&\colon [0,1]^3\to [0,1]\notag \\
&(x,y,z)\mapsto {1\over 4} x^2+yz.
\end{align}
One can verify that $\bar f$ is monotonic-plus-smooth such that $\bar f\in \tF(1,1)\cap \tB(1/2)$. Here, we choose the smoothness $\alpha=1$ to ensure the perturbation term $1/d$ in~\eqref{eq:newrep} is dominated by the approximation error. As a result, the construction \eqref{eq:newperm}-\eqref{eq:newf} satisfies~\eqref{eq:newrep}, and by applying Proposition~\ref{lem:approx} to $\bar f$, we obtain
\[
{1\over d^m}\inf_{\tB\in\caliB(k,\ell)}\FnormSize{}{\Theta-\tB}^2 \leq {L^2 \over k^{2{\min(1,\ell+1)}}}+{1\over d^2}\leq {L+1\over k^{\min(1,\ell+1)}}.
\]
In conclusion, the signal tensor constructed in the Example~\ref{ex:free} can be regarded as a tensor generated from $\tF(1,1)\cap \tB(1/2)$. 
Therefore, our Borda count algorithm is applicable to this case with the claimed accuracy rate~\eqref{eq:rateBC} for $(\alpha,\beta)=(1,1/2)$. This example highlights the flexibility of our permuted monotonic-plus-smooth assumption, which accommodate a broader range of non-monotonic functions through permutations.

\section{Proofs of main theorems}\label{app:theorem}

\subsection{Proof of Proposition~\ref{lem:approx}}
\begin{proof}[Proof of Proposition~\ref{lem:approx}]
We denote $\tE_k$ as the $m$-way partition by
\begin{align}
    \tE_k = \{\bigtimes_{a=1}^m z^{-1}(j_a)\colon (j_1,\ldots,j_m)\in [k]^m\}, 
\end{align}
where $z\colon [d]\rightarrow[k]$ is the balanced clustering function such that $z(i) = \lceil ki/d\rceil,$ for all $i \in[d]$, and we use the shorthand $\bigtimes_{a=1}^m$ to denote the Cartesian product of $m$ sets. For a given $\mj:=(j_1,\ldots,j_m)$ and the related partition  $ \bigtimes_{a=1}^m z^{-1}(j_a)\in\tE_k$, fix any index  $\omega_{\mj} \in \bigtimes_{a=1}^m z^{-1}(j_a)$. Then, we have 
\begin{align}\label{eq:ind}
\|\omega-\omega_{\mj}\|_\infty \leq \frac{d}{k},\quad \text{for all }\omega\in \bigtimes_{a=1}^m z^{-1}(j_a). 
\end{align} 

We define the block-wise degree-$\ell$ polynomial tensor $\tB$ based on the partition $\tE_k$ as
\begin{align}
    \tB(\omega) = \sum_{\mj\in[k]^m}\text{Poly}_{\min(\lfloor \alpha \rfloor,\ell), \omega_{\mj}}\left({\omega\over d}\right)\mathds{1}\{\omega \in \bigtimes_{a=1}^m z^{-1}(j_a)\},
\end{align}
where $\text{Poly}_{\min(\lfloor \alpha \rfloor,\ell),\omega_{\mj}}$ denotes a degree-$\ell$ polynomial function defined by
\begin{align}
\text{Poly}_{\min(\lfloor \alpha \rfloor,\ell),\omega_{\mj}} (\mx)=\sum_{\kappa:|\kappa|\leq \min(\lfloor \alpha\rfloor,\ell+1)}{1\over \kappa!}\left(\mx-{\omega_{\mj}\over d}\right)^\kappa \nabla_\kappa f\left({\omega_{\mj}\over d}\right).
\end{align}
It follows from Taylor's theorem that
\begin{align}\label{eq:polyapp}
\left|f\left({\omega\over d}\right)-\text{Poly}_{\min(\lfloor \alpha \rfloor,\ell),\omega_{\mj}}\left({\omega\over d}\right)\right|\leq L\left\|{\omega-\omega_{\mj}\over d}\right\|_{\infty}^{\min(\alpha,\ell+1)}\leq L \left({1\over k} \right)^{\min(\alpha,\ell+1)},
\end{align}
for all $\omega\in  \bigtimes_{a=1}^m z^{-1}(j_a)$. 

Based on the construction of block-wise degree-$\ell$ polynomial tensor $\tB$, we have 
\begin{align}
   & \frac{1}{d^m}\FnormSize{}{\Theta-\tB}^2\\
=&\ \frac{1}{d^m} \sum_{\omega\in[d]^m}|\Theta(\omega)-\tB(\omega)|^2\\
=&\ \frac{1}{d^m}\sum_{\mj\in[k]^m}\sum_{\omega \in \bigtimes_{a=1}^m z^{-1}(j_a)}\bigg|f\left(\frac{\omega}{d}\right)-\text{Poly}_{\min(\lfloor \alpha \rfloor,\ell),\omega_{\mj}}\left(\frac{\omega}{d}\right)\bigg|^2\\
    \leq
    &\ \frac{L^2}{k^{2\min(\alpha,\ell+1)}},
\end{align}
where the last inequality uses~\eqref{eq:polyapp}.
\end{proof}

\subsection{Proof of Theorem~\ref{thm:LSE}}\label{sec:proofLSE}
\begin{proof}[Proof of Theorem~\ref{thm:LSE}]
By Proposition~\ref{lem:approx}, there exists a block-wise polynomial tensor $\tB\in\caliB(k,\ell)$ such that
\begin{equation}\label{eq:approx}
\FnormSize{}{\tB-\Theta}^2\lesssim {L^2d^m\over k^{2\min(\alpha,\ell+1)}}.
\end{equation}
By the triangle inequality,
\begin{equation}\label{eq:tri}
\FnormSize{}{\hat\Theta^\textup{LSE}\circ\hat\pi^\textup{LSE} -\Theta\circ \pi}^2\leq 2\FnormSize{}{\hat\Theta^\textup{LSE}\circ\hat\pi^\textup{LSE} -\tB\circ \pi}^2+2\underbrace{\FnormSize{}{\tB\circ \pi-\Theta\circ \pi}^2}_{\textup{Proposition~\ref{lem:approx}}}.
\end{equation}
Therefore, it suffices to bound $\FnormSize{}{\hat\Theta^\textup{LSE}\circ\hat\pi^\textup{LSE} -\tB\circ \pi}^2$. By the global optimality of least-square estimator, we have
\begin{align}
\FnormSize{}{\hat\Theta^\textup{LSE}\circ\hat\pi^\textup{LSE} -\tB\circ \pi}&\leq \left\langle {\hat \Theta^\textup{LSE}\circ \hat \pi^\textup{LSE}-\tB\circ \pi \over \FnormSize{}{ \hat \Theta^\textup{LSE}\circ \hat \pi^\textup{LSE}-\tB\circ \pi}},\ \tE+(\Theta\circ \pi-\tB\circ \pi) \right \rangle \\
&\leq \sup_{\pi, \pi'\colon[d]\to[d]}\sup_{\tB, \tB'\in \caliB(k,\ell)} \left\langle {\tB'\circ \pi'-\tB\circ \pi \over \FnormSize{}{\tB'\circ \pi'-\tB\circ \pi}}, \tE \right \rangle+\underbrace{\FnormSize{}{\tB\circ \pi-\Theta\circ \pi}}_{\textup{Proposition~\ref{lem:approx}}}.
\end{align}
Now we bound inner product term. For fixed $\pi,\pi'$, let $\mP$ and $\mP'$ be permutation matrices corresponding to permutations $\pi$ and $\pi'$ respectively. 
We express vectorized block-wise degree-$\ell$ polynomial tensors, $\text{vec}(\tB)$ and $\text{vec}(\tB')$, by discrete polynomial functions. Specifically, denote $\text{vec}(\tB)=\mX\mbeta$ and $\text{vec}(\tB')=\mX\mbeta'$, where $\mX\in\bbR^{d^m\times k^m(\ell+m)^\ell}$ is a design matrix consisting of $m$-multivaraite degree-$\ell$ polynomial basis over grid design $\{1/d,\ldots,d/d\}$ (or $\{x_i\}_{i=1}^d$ for general designs~\eqref{eq:randommodel}), $\mbeta, \mbeta'\in\bbR^{k^m(\ell+m)^\ell}$ are corresponding coefficient vectors. Notice that the number of coefficients for $m$-multivariate polynomial of degree-$\ell$ is ${\ell+m \choose \ell} \approx (\ell+m)^\ell$; we choose to use $(\ell+m)^\ell$ for each block for notational simplicity. Therefore, we rewrite the inner product by
\begin{align}\label{eq:inner}
   \left\langle {\tB'\circ \pi'-\tB\circ \pi \over \FnormSize{}{\tB'\circ \pi'-\tB\circ \pi}}, \tE \right \rangle &=  \left\langle {(\mP')^{\otimes m}\text{vec}(\tB')-(\mP)^{\otimes m}\text{vec}(\tB) \over \FnormSize{}{(\mP')^{\otimes m}\text{vec}(\tB')-(\mP)^{\otimes m}\text{vec}(\tB)}}, \tE \right \rangle\notag \\
    &=\left\langle {(\mP')^{\otimes m}\mX\mbeta'-(\mP)^{\otimes m}\mX\mbeta \over \FnormSize{}{(\mP')^{\otimes m}\mX\mbeta'-(\mP)^{\otimes m}\mX\mbeta}}, \tE \right \rangle\notag
    \\&=
     \left\langle {\mA\mc \over \vnormSize{}{\mA\mc}}, \tE \right \rangle
\end{align}
where we denote $\vnormSize{}{\cdot}$ denotes the vector 2-norm, and we define $\mA := \begin{pmatrix}\mP' &-\mP\end{pmatrix}\begin{pmatrix}\mX & 0\\ 0&\mX\end{pmatrix}\in\bbR^{d^m\times 2k^m(\ell+m)^\ell}$ and $\mc:=\begin{pmatrix} \mbeta'\\\mbeta\end{pmatrix}\in\bbR^{2k^m(\ell+m)^\ell}$. 

By Lemma~\ref{lem:embedding}, we have 
\begin{align}\label{eq:embedding}
\sup_{\mc\in \bbR^{2k^m(\ell+m)^\ell}}\left\langle {\mA\mc \over \vnormSize{}{\mA\mc}}, \tE \right \rangle \leq \sup_{\mc\in \bbR^{2k^m(\ell+m)^\ell } }\left\langle {\mc\over \vnormSize{}{\mc}},\ \me \right\rangle,
\end{align}
where $\me\in\bbR^{2k^m(\ell+m)^\ell}$ is a sub-Gaussian random vector with variance proxy bounded by $\sigma^2$. 
By the union bound of sub-Gaussian maxima over countable set $\{\pi,\pi'\colon [d]\to[d]\}$, we obtain
\begin{align}\label{eq:union}
    \mathbb{P}&\left(\sup_{\pi, \pi'\colon[d]\to[d]}\sup_{\tB, \tB'\in \caliB(k,\ell)} \left\langle {\tB'\circ \pi'-\tB\circ \pi \over \FnormSize{}{\tB'\circ \pi'-\tB\circ \pi}}, \tE \right \rangle\geq t\right)\notag \\
    &\leq \sum_{\pi,\pi'\in[d]^d}\mathbb{P}\left(\sup_{\tB, \tB'\in \caliB(k,\ell)} \left\langle {\tB'\circ \pi'-\tB\circ \pi \over \FnormSize{}{\tB'\circ \pi'-\tB\circ \pi}}, \tE \right \rangle\geq t\right)\notag \\
    &\leq d^d\mathbb{P}\left( \sup_{\mc\in \bbR^{2k^m(\ell+m)^\ell } }\left\langle {\mc\over \vnormSize{}{\mc}},\ \me \right\rangle\geq t\right)\notag\\
    &\leq \exp\left(-\frac{t^2}{8\sigma^2} +k^m(\ell+m)^\ell\log 6 + d\log d\right),
\end{align}
where the second inequality is from \eqref{eq:embedding} and the last inequality is from Lemma~\ref{lem:subga}. Setting $t = C\sigma\sqrt{k^m(\ell+m)^\ell + d\log d}$ in~\eqref{eq:union} for a sufficiently large $C>0$ gives
\begin{align}
    \sup_{\pi, \pi'\colon[d]\to[d]}\sup_{\tB, \tB'\in \caliB(k,\ell)} \left\langle {\tB'\circ \pi'-\tB\circ \pi \over \FnormSize{}{\tB'\circ \pi'-\tB\circ \pi}}, \tE \right \rangle\lesssim \sigma\sqrt{k^m(\ell+m)^\ell +d\log d},
\end{align}
with high probability.

Combining the inequalities~\eqref{eq:approx}, \eqref{eq:tri} and \eqref{eq:union} yields the desired conclusion
\begin{equation}\label{eq:three}
\FnormSize{}{\hat\Theta^\textup{LSE}\circ\hat\pi^\textup{LSE} -\Theta\circ \pi}^2\lesssim \sigma^2\left(k^m(\ell+m)^\ell +d\log d\right)+{L^2d^m \over k^{2\min(\alpha, \ell+1)}} .
\end{equation}

Finally, optimizing~\eqref{eq:three} with respect to $(k,l)$ gives that 
\begin{align}
     \eqref{eq:three}\ \lesssim 
     \begin{cases} 
    L^2\left({\sigma\over L}\right)^{4\alpha\over m+2\alpha} d^{-\frac{2m\alpha}{m+2\alpha}}, & \text{ when } \alpha < m(m-1)/2,\\
     \sigma^2 d^{-(m-1)}\log d, &\text{ when } \alpha \geq m(m-1)/2,
    \end{cases}
\end{align}
under the choice
\[
\ell^* = \min(\lceil \alpha-1\rceil,(m-2)(m+1)/2),\quad k^* = \left\lceil \left( d^m L^2 /\sigma^2\right)^{1\over m+2\min(\alpha,\ell^*+1)} \right\rceil.
\]

\end{proof}

\subsection{Proof of Theorem~\ref{thm:minimax}}
\begin{proof}[Proof of Theorem~\ref{thm:minimax}] By the definition of the tensor space, we seek the minimax rate $\varepsilon^2$ in the following expression
\begin{equation}\label{eq:final}
\inf_{(\hat \Theta,\hat \pi)}\sup_{\Theta\in \tP(\alpha,L)}\sup_{\pi\in\Pi(d,d)} \mathbb{P}\left({1\over d^m}\FnormSize{}{\Theta\circ \pi-\hat \Theta\circ \hat \pi}^2 \geq \varepsilon^2 \right).
\end{equation}
On one hand, if we fix a permutation $\pi\in\Pi(d,d)$, the problem can be viewed as a classical $m$-dimensional $\alpha$-smooth nonparametric regression with $d^m$ sample points. The minimax lower bound is known to be $\varepsilon^2=L^2\left(\sigma\over L\right)^{4\alpha\over m+2\alpha}d^{-{2m\alpha\over m+2\alpha}}$. On the other hand, if we fix $\Theta\in\tP(\alpha,L)$, the problem become a new type of convergence rate due to the unknown permutation. We refer to the resulting error as the permutation rate, and we will prove that $\varepsilon^2=\sigma^2 d^{-(m-1)}\log d$. Since our target is the sum of the two rates, it suffices to prove the two different rates separately. In the following arguments, we will proceed by this strategy. 

\paragraph{Nonparametric rate.} The nonparametric rate for $\alpha$-smooth function is readily available in the literature. We cite the results here for self-completeness. 

\begin{lem}[Minimax rate for $\alpha$-smooth function estimation; see Section 2 in \citet{stone1982optimal}]\label{lem:non} Consider a sample of $N$ data points, $(\mx_1,Y_1), \ldots, (\mx_N,Y_N)$, where $\mx_n=({i_1\over d},\ldots,{i_m\over d})\in[0,1]^m$ is the $m$-dimensional predictor and $Y_n\in\bbR$ is the scalar response for $n\in[N]$. Consider the observation model
\[
Y_n=f(\mx_n)+\varepsilon_n,\quad \text{with}\ \varepsilon_n \sim \text{i.i.d. }  \tN(0,1), \quad\text{ for all }n\in[N].
\]
Assume $f$ is in the $\alpha$-H\"older smooth function class, denoted by $\tF(\alpha,L)$. Then,
\begin{equation}
\inf_{\hat f}\sup_{f\in \tF(\alpha,L)}\mathbb{P}\left(\| f-\hat f\|_{L^2}\geq \sigma^{4\alpha\over m+2\alpha}L^{2m\over m+2\alpha} N^{-{2\alpha\over m+2\alpha}}\right)\geq 0.9.
\end{equation}
\end{lem}
Our desired nonparametric rate readily follows from Lemma~\ref{lem:non} by taking sample size $N=d^m$ and function norm $\| f-\hat f \|_{L^2}:={1\over d^m}\FnormSize{}{\Theta-\hat \Theta}^2$. In summary, for a given permutation $\pi\in\Pi(d,d)$, we have
\begin{align}\label{eq:non}
\inf_{\hat \Theta}\sup_{\Theta\in \tP(\alpha,L)}\mathbb{P}\left({1\over d^m}\FnormSize{}{\hat \Theta\circ \pi-\Theta\circ \pi}^2\geq L^2\left(\sigma\over L\right)^{4\alpha\over m+2\alpha} d^{-{2m\alpha\over m+2\alpha}}\right) \geq 0.9.
\end{align}

\paragraph{Permutation rate.}
Since nonparametric rate dominates permutation rate when $\alpha< 1$, it is sufficient to prove the permutation rate lower bound for $\alpha\geq 1.$
We first show the minimax permutation rate for $k$-block degree-$0$ tensor family $\tB(k,0)$, and then construct a smooth $f\in\tF(\alpha,L)$ to mimic the constant block tensors. 

Let $\Pi(d,k)$ denote the collection of all possible onto mappings from $[d]$ to $[k]$.
Lemma~\ref{lem:permutation} shows the permutation rate over  $k$-block degree-0 tensor family $\tB(k,0)$ is $\sigma^2 d^{-(m-1)}\log k$. 
\begin{lem}[Permutation error for tensor block model]\label{lem:permutation}
For every given integer $k\in[d]$, there exists a core tensor $\tS\in\bbR^{k\times \cdots \times k}$ such that the corresponding $d$-dimensional block-$k$ Gaussian tensor block model estimation problem has the minimax rate:
\begin{equation}\label{eq:givenC}
\inf_{\hat \Theta}\sup_{z \in \Pi(d,k)}\mathbb{P}\left( {1\over d^m}\FnormSize{}{\hat \Theta-\tS\circ z} ^2\gtrsim \frac{\sigma^2\log k} {d^{m-1}}\right)\geq 0.9.
\end{equation}
\end{lem}
The proof of Lemma~\ref{lem:permutation} is constructive and deferred to Section~\ref{sec:tech}. 

To prove the permutation rate, we fix a core tensor $\tS\in\bbR^{k\times \cdots \times k}$ satisfying~\eqref{eq:givenC}, and use it to construct the smooth tensor. We construct a function $f\in\tF(\alpha,L)$ that mimics the core tensor $\tS$ in block tensor family $\tB(k,0).$  For notational simplicity, we do not distinguish the fractional number and its rounding integer. For example, we simply write $d/2$ (instead of $\lceil d/2 \rceil$ or $\lfloor d/2 \rfloor$) to represent its rounding integer. Define $k = d^\delta$ for some $\delta\in(0,1)$, which will be specified later. Consider a smooth function $K(x)$ that is infinitely differentiable,
\begin{align}
    K(x) = C_k\exp\left(-\frac{1}{1-64x^2}\right)\mathds{1}\left\{|x|<\frac{1}{8}\right\},
\end{align}
where $C_k>0$ satisfies $\int K(x)dx = 1.$  Then, we define a smooth convolution function as
\begin{align}
    \psi(x) = \int_{-3/8}^{3/8}K(x-y)dy.
\end{align}
The smooth convolution function has support $[-1/2,1/2]$ and takes value 1 on the interval $[-1/4,1/4]$. 
For a given core tensor $\tS$ in Lemma~\ref{lem:permutation}, we define an $\alpha$-smooth function
\begin{align}\label{eq:constructf}
    f(x_1,\ldots,x_m) = \sum_{(a_1,\ldots,a_m)\in[k]^m}\left[\tS(a_1,\ldots,a_m)\prod_{i=1}^m \psi\left(kx_i-a_i+\frac{1}{2}\right) \right].
\end{align}
One can verify that $f\in\tF(\alpha,L)$ as long as the constant $\delta=\delta(\alpha,L)>0$ is set sufficiently small depending on $\alpha$ and $L$. Notice that for any $(a_1,\ldots,a_m)\in [k]^m$, 
\begin{align}\label{eq:1}
f(x_1,\ldots,x_m) = \tS(a_1,\ldots,a_m), \quad \text{if } (x_1,\ldots,x_m)\in \bigtimes_{i=1}^m\left[\frac{a_i-3/4}{k},\frac{a_i-1/4}{k}\right].
\end{align}
From this observation, we define a sub-domain $I\subset [d]$ such that
\begin{align}\label{eq:2}
    I = \left(\bigcup_{a=1}^k\left[\frac{d(a-3/4)}{k},\frac{d(a-1/4)}{k}\right]\right)\bigcap [d].
\end{align}
We have that $|I|=d/2$ by definition. Let $\Theta(\tS)\in\mathbb{R}^{d\times \cdots \times d}$ denote the tensor induced by $f$ in~\eqref{eq:constructf}. We use subscript $I$ to denote objects when restricted in the indices set $I$. For example, $\Theta_I(\tS)\in\mathbb{R}^{d/2 \times \cdots \times d/2}$ denotes the sub-tensor with indices in $I$, and $\normSize{}{\cdot}$ denotes the sum of squares over indices in $I$. Based on~\eqref{eq:1} and~\eqref{eq:2}, $\Theta_I(\tS)$ has block structure with the core tensor $\tS$. We use $\Pi(d/2,d/2)$ to denote the set of all permutations on $I$ while fixing indices on $[d]\setminus I$; that is, $\Pi(d/2,d/2)=\{\pi\colon I\to I\}\cong \{\pi\in\Pi(d,d)\colon \pi(i) = i \text{ for } i\in[d]\setminus I\})$. Then, we have
\begin{align}\label{eq:lower}
&\inf_{(\hat \Theta,\hat \pi)}\sup_{\pi\in \Pi(d,d)}\mathbb{P}\left({1\over d^m}\FnormSize{}{\hat \Theta\circ \hat \pi-\Theta(\tS)\circ \pi}^2\geq \varepsilon^2\right)\notag \\
\stackrel{(*)}{\geq }&\ \inf_{\hat \Theta}\sup_{\pi\in \Pi(d,d)}\mathbb{P}\left({1\over d^m}\normSize{}{\hat \Theta-\Theta(\tS)\circ \pi}^2\geq \varepsilon^2\right)\notag \\
\stackrel{(**)}{\geq}& \ \inf_{\hat  \Theta}\sup_{\pi \in\Pi(d/2,d/2)}\mathbb{P}\left({1\over (d/2)^m}\FnormSize{}{\hat \Theta_{I}-\Theta_{I}(\tS)\circ \pi}^2\geq 2^m\varepsilon^2\right)\notag \\
\stackrel{(***)}{\geq} & \ \inf_{\hat  \Theta}\sup_{z\in \Pi(d/2,k)}\mathbb{P}\left({1\over (d/2)^m}\FnormSize{}{\hat \Theta_I-\tS\circ z}^2\geq 2^m\varepsilon^2\right),
\end{align}
where $(*)$ absorbs the estimate $\hat \pi$ into the estimate $\hat \Theta$, and $(**)$ uses the permutation collection $\Pi(d/2,d/2)\cong \{\pi: I\to I\}$, and $(***)$ uses the block structure of $\Theta_I(S)$. Therefore, we reduce the problem of estimating $\pi \colon [d]\to[d]$ in the $\alpha$-smooth tensor to estimating clustering $z: I\to[k]$ in the sub-tensor. Applying Lemma~\ref{lem:permutation} to~\eqref{eq:lower} by using $d/2$ in the place of $d$ and $k=d^{\delta} $ for the constant $\delta=\delta(\alpha,L)>0$ yields the desired conclusion.

\paragraph{Combining two rates.} Now, we combine~\eqref{eq:non} and~\eqref{eq:lower} to get the desired lower bound. For any $\Theta$  generated as in  \eqref{eq:rep} with $f\in\tF(\alpha,L)$, by union bound, we have
\begin{align}
&\mathbb{P}\left\{ {1\over d^m}\FnormSize{}{\hat \Theta- \Theta}^2\gtrsim L^2\left(\sigma\over L\right)^{4\alpha\over m+2\alpha}d^{-{2m\alpha\over m+2\alpha}}+\frac{\sigma^2\log d}{d^{m-1}}\right\}\\
\geq &\quad \mathbb{P}\left\{ {1\over d^m}\FnormSize{}{\hat \Theta-\Theta}^2\gtrsim L^2\left(\sigma\over L\right)^{4\alpha\over m+2\alpha}d^{-{2m\alpha\over m+2\alpha}}\right\} +\mathbb{P}\left\{ {1\over d^m}\FnormSize{}{\hat \Theta- \Theta}^2\gtrsim \frac{\sigma^2\log d}{d^{m-1}}\right\}-1.
\end{align}
Taking sup on both sides with the property
\[
\sup_{\substack{\Theta\in \tP(\alpha,L)\\\pi\in \Pi(d,d)}}(f(\pi)+g(\Theta))=\sup_{\pi\in \Pi(d,d)}f(\pi)+\sup_{\Theta\in \tP(\alpha,L)}g(\Theta)
\]
yields the desired rate~\eqref{eq:minimax}. 
\end{proof}

\subsection{Proof of Theorem~\ref{thm:nopoly}}\label{sec:nopoly}
The proof of Theorem~\ref{thm:nopoly} leverages results of hypergraphic planted clique and constant higher-order clustering problems. We first briefly explain the constant higher-order clustering problems. We then prove the main result.

\subsubsection{Constant higher-order clustering and computational lower bound}  
Let $\mk = (k_1,\ldots,k_m)$ and $\md = (d_1,\ldots,d_m)$. We introduce the constant high-order clustering (CHC) problem~\citep{luo2022tensor}. Consider a data tensor $\tY\in\bbR^{d_1\times \cdots \times d_m}$ generated from the signal plus noise model
\begin{equation}\label{eq:gaussian}
\tY=\Theta+\tE, 
\end{equation}
where the entries in $\tE$ are i.i.d. drawn from Gaussian distribution, and the signal tensor $\Theta$ contains a constant planted structure:
\begin{align}\label{eq:chcdef}
    \Theta\in\Theta_{\text{CHC}}(\mk,\md,\lambda) :=\{\lambda'\mathds{1}_{I_1}\otimes \cdots\otimes \mathds{1}_{I_m}\colon|I_i| = k_i,\ \text{for all }i\in[m],\ \lambda' \geq \lambda\}.
\end{align}
Here $I_i\subset [d_i]$ denotes a subset of indices, $|\cdot|$ denotes the cardinality of the set, $\mathds{1}_{I_i}$ is the $d_i$-dimensional indicator vector such that $(\mathds{1}_{I_i})_j = 1$ if $j\in I_i$ and 0 otherwise. The CHC detection problem is to test the following hypothesis based on the observed tensor $\tY$, 
\begin{align}\label{eq:chct}
    H_0\colon \Theta = 0\quad\text{v.s.}\quad H_1\colon \Theta\in\Theta_{\text{CHC}}(\mk,\md,\lambda).
\end{align}

The following proposition provides the asymptotic regime for impossible polynomial-time detection of CHC under Conjecture~\ref{conj:1}. This proposition plays important role to prove our Theorem~\ref{thm:nopoly}.
 \begin{prop}[Theorem 15 in \citep{luo2022tensor}]\label{prop:CHC} 
Consider CHC detection problem in \eqref{eq:chct} in the Gaussian noise model~\eqref{eq:gaussian} under the asymptotic regime $d\rightarrow\infty$ satisfying
\begin{align}\label{eq:equal}
  d =  d_1= \cdots= d_m ,\quad
    k = k_1 = \cdots = k_m = d^\delta,\quad \lambda  = d^{-\gamma},
\end{align}with $0\leq\delta\leq1$ and $\gamma >(m\delta-m/2)\vee 0.$ Then, under Conjecture~\ref{conj:1}, for any polynomial-time test sequence $\{\phi\}_d\colon \tY\mapsto\{0,1\}$, we have
\begin{align}
    \liminf_{d\rightarrow\infty} \left\{\mathbb{P}_{\tH_0}(\phi(\tY) = 1)+ \sup_{\Theta\in\Theta_{\textup{CHC}}(\mk,\md,\lambda)}\mathbb{P}_{\Theta}(\phi(\tY) =0) \right\}\geq\frac{1}{2}.
\end{align}
\end{prop}
 
\subsubsection{Proof of Theorem~\ref{thm:nopoly}}
 \begin{proof}[Proof of Theorem~\ref{thm:nopoly}]
 Assume that  the true signal $\Theta\in\mathbb{R}^{d\times \cdots \times d}$ has constant planted structure with a given $\mk=(k,\ldots,k)$ such that 
\begin{align}
    \Theta \in\Theta_{\text{CHC}}(\mk,\md,\lambda) = \{ \lambda \mathds{1}_{I}\otimes \cdots\otimes \mathds{1}_{I}\colon  |I| = k\}.
\end{align}
We have $\Theta_{\text{CHC}}\subset \tP_{\textup{gen}}$, because we can set the infinitely smooth function $f: [0,1]^m \to[0,1]$ by
\begin{equation}\label{eq:fun}
f(x_1,\ldots,x_m)=\lambda \prod_{i\in[m]}x_i,
\end{equation}
under the choice $x_i=\mathds{1}_{i\in I}$ for all $i\in[m]$. Then, there is one-to-one correspondence between tensors in $\Theta_{\textup{CHC}}(\mk,\md,\lambda)$ and tensors generated by the above $f$ and $\{x_i\}_{i\in[d]}$.

We consider the regime where polynomial-time solvable test is impossible based on Proposition~\ref{prop:CHC}. We set $\delta=1/2$, $k = c_1d^{\delta}$, and $\lambda = c_2d^{-\gamma}$ for any fixed $\gamma \in (0, {2\alpha -m \over m})$, so that any polynomial-time test sequence $\phi$ satisfies 
\[
 \liminf_{d\rightarrow\infty} \left\{\mathbb{P}_{\tH_0}(\phi(\tY) = 1)+\sup_{\Theta \in \Theta_{\textup{CHC}}(\mk,\md,\lambda)}\mathbb{P}_{\Theta}(\phi(\tY) =0) \right\}\geq\frac{1}{2}.
\]
The choice of $\lambda$ is possible given that $\alpha > m/2$. Notice that the choice $\lambda \lesssim O(1)$ ensures the function~\eqref{eq:fun} satisfies the definition~\eqref{eq:holder} for all $\alpha>0$. 

We prove by contradiction. Assume that there exists a hypothetical estimator $\hat\Theta$ from a polynomial-time algorithm  that attains the rate $\text{Rate}(d)$. Specifically, there exists a constant $b>0$ such that 
\begin{align}\label{eq:esterror}
   \limsup_{d\rightarrow\infty}\frac{1}{\textup{Rate}(d)} \sup_{\Theta\in\Theta_{\text{CHC}}(\mk,\md,\lambda)}\frac{1}{d^m}\mathbb{E}\FnormSize{}{\hat\Theta-\Theta}^2\leq b.
\end{align}
By Markov's inequality, the inequality~\eqref{eq:esterror} implies that, when $d$ is sufficiently large, for all $\Theta\in\Theta_{\text{CHC}}(\mk,\md,\lambda)$ and all $u>0$, we have
\begin{align}\label{eq:markov}
    \FnormSize{}{\hat\Theta-\Theta}\leq u \sqrt{\text{Rate}(d)d^{m}},
\end{align}
with probability at least $1-b/u$.
Consider the hypothesis test in \eqref{eq:chct}. We employ the following test 
\begin{align}
    \phi(\tY) = \mathds{1}(\FnormSize{}{\hat\Theta}\geq u\sqrt{\text{Rate}(d)d^m}).
\end{align}
The Type I error of the test $\phi$ is controlled by
\begin{align}
    \mathbb{P}_0(\FnormSize{}{\hat\Theta}\geq u\sqrt{\text{Rate}(d)d^m})  = \mathbb{P}_0(\FnormSize{}{\hat\Theta-\Theta}\geq u\sqrt{\text{Rate}(d)d^m}) \leq b/u.
\end{align}
For Type II error, we obtain,
\begin{align}
\sup_{\Theta\in\Theta_{\text{CHC}}(k,d,\lambda)}\mathbb{P}_\Theta(\phi(\tY) = 0) &= \sup_{\Theta\in\Theta_{\text{CHC}}(k,d,\lambda)}\mathbb{P}_\Theta(\FnormSize{}{\hat\Theta}< u \sqrt{\text{Rate}(d)d^m})\\&\leq\sup_{\Theta\in\Theta_{\text{CHC}}(k,d,\lambda)}\mathbb{P}_\Theta(\FnormSize{}{\hat\Theta-\Theta}^2> \FnormSize{}{\Theta}^2-u^2\text{Rate}(d)d^m)\\&\stackrel{(*)}{\leq}\sup_{\Theta\in\Theta_{\text{CHC}}(\mk,\md,\lambda)}\mathbb{P}_\Theta(\FnormSize{}{\hat\Theta-\Theta}^2> u^2\text{Rate}(d)d^m)\\&\stackrel{(**)}{\leq} b/u.
\end{align}
The inequality $(*)$ holds because 
\begin{align}\label{eq:regime}
    \FnormSize{}{\Theta}^2\geq \lambda^2 k^m = c_1^mc_2^2d^{{m\over 2}-\gamma} \geq 2 u^2\text{Rate}(d)d^m \asymp d^{{m\over 2}-{2\alpha-m \over m}}
\end{align}
where the last inequality is true under the regime $c_1^mc_2>2u^2$. We can always choose constants $c_1$ and $c_2$ given the value $u$. The inequality $(**)$ holds because of the statement~\eqref{eq:markov}. 
Putting Type I and II errors together, we obtain
\begin{align}
 \mathbb{P}_{\tH_0}(\phi(\tY) = 1)+ \sup_{\Theta\in\Theta_{\text{CHC}}(\mk,\md,\lambda)}\mathbb{P}_{\Theta}(\phi(\tY) =0) \leq 2b/u < 1/2,
\end{align}
for $u>4b$. This fact contradicts the Proposition~\ref{prop:CHC}.  Therefore, there is no polynomial-time $\hat\Theta$ satisfying \eqref{eq:esterror}.
\end{proof}

\subsection{Proof of Theorem~\ref{thm:BC}}
We first present a lemma to show the estimation error of $\hat \pi$. The exponent $\beta$ measures the difficulty for permutation estimation. We find that a larger $\beta$ guarantees a faster consistency rate of $\hat\pi^{\textup{BC}}$, which is represented below.

\begin{lem}[Permutation error]\label{lem:permute}
Consider the sub-Gaussian tensor model~\eqref{eq:obs} with $f\in \tM(\beta)$. With high probability, we have
\begin{align}
    \textup{Loss}(\pi,\hat\pi^{\textup{BC}}):= \frac{1}{d}\max_{i\in[d]}|\pi(i)-\hat\pi^{\textup{BC}}(i)|\lesssim \left(d^{-(m-1)/2}\sqrt{\log d}\right)^{\beta}.
\end{align}
\end{lem}
The proof of Lemma~\ref{lem:permute} is provided in Section~\ref{sec:tech}. 

\begin{proof}[Proof of Theorem~\ref{thm:BC}]
By Proposition~\ref{lem:approx}, there exists a block-wise polynomial tensor $\tB\in\caliB(k,\ell)$ satisfying \eqref{eq:approx}.
By the triangle inequality,
we decompose estimation error into three terms,
\begin{align}\label{eq:3decomp}
    &\FnormSize{}{\hat\Theta^{\text{BC}}\circ \hat\pi^{\text{BC}}-\Theta\circ\pi}\notag \\
    \leq&\ \FnormSize{}{\hat\Theta^{\text{BC}}\circ \hat\pi^{\text{BC}}-\tB\circ\hat\pi^{\text{BC}}}+\FnormSize{}{\tB\circ\hat\pi^{\text{BC}}-\Theta\circ\hat\pi^{\text{BC}}}+\FnormSize{}{\Theta\circ\hat\pi^{\text{BC}}-\Theta\circ\pi}\nonumber\\
    =&\underbrace{\FnormSize{}{\hat\Theta^{\text{BC}}-\tB}}_{\textup{Nonparametric error}}+\ \underbrace{\FnormSize{}{\Theta\circ\hat\pi^{\text{BC}}-\Theta\circ\pi}}_{\textup{Permutation error}}+\underbrace{\FnormSize{}{\tB-\Theta}}_{\textup{Proposition~\ref{lem:approx}}}.
\end{align}
Therefore, it suffices to bound two terms $\FnormSize{}{\Theta\circ\hat\pi^{\text{BC}}-\Theta\circ\pi}$ and $\FnormSize{}{\hat\Theta^{\text{BC}}-\tB}$ separately. 

\noindent{\bf Permutation error}.
    For any $(i_1,\ldots,i_m)\in[d]^m$, we have 
    \begin{align}
        &\left|\Theta\left(\hat\pi^{\textup{BC}}(i_1),\ldots,\hat\pi^{\textup{BC}}(i_m)\right)-\Theta\left(\pi(i_1),\ldots,\pi(i_m)\right)\right|\\&\leq \left\|\left(\frac{\hat\pi^{\textup{BC}}(i_1)}{d},\ldots,\frac{\hat\pi^{\textup{BC}}(i_m)}{d}\right)-\left(\frac{\pi(i_1)}{d},\ldots,\frac{\pi(i_m)}{d}\right)\right\|_\infty^{\min(\alpha,1)}\\&\leq \left[\frac{1}{d}\max_{i\in[d]}\left|\hat\pi^{\textup{BC}}(i)-\pi(i)\right|\right]^{\min(\alpha,1)}\\&\lesssim \left(\sigma d^{-(m-1)/2}\sqrt{\log d}\right)^{\beta{\min(\alpha,1)}},
    \end{align}
    where the first inequality is from the $\alpha$-H\"older smoothness of $\Theta$, and the last inequality is from Lemma~\ref{lem:permute}. Therefore, we obtain the upper bound of the permutation error 
    \begin{align}\label{eq:permerror}
        \frac{1}{d^m}\FnormSize{}{\Theta\circ\hat\pi^{\text{BC}}-\Theta\circ\pi}^2\lesssim \left(\sigma^2 \frac{\log d}{d^{m-1}}\right)^{\beta{\min(\alpha,1)}}.
    \end{align}
    
    \noindent{\bf Nonparametric error.} 
    Recall that Borda count estimation is defined by $\hat\Theta^{\textup{BC}} := \argmin_{\Theta\in\caliB(k,\ell)}\FnormSize{}{\tY\circ(\hat\pi^{\textup{BC}})^{-1}-\Theta}^2$. By the optimality of least-square estimator, we have
\begin{align}
\FnormSize{}{\hat\Theta^\textup{BC} -\tB}&\lesssim \left\langle {\hat \Theta^\textup{BC}-\tB \over \FnormSize{}{ \hat \Theta^\textup{BC}-\tB}},\ \tY\circ(\hat\pi^{\textup{BC}})^{-1}-\tB \right \rangle\\ 
&= \left\langle {\hat \Theta^\textup{BC}-\tB \over \FnormSize{}{ \hat \Theta^\textup{BC}-\tB}},\tE\circ \hat \pi^{\textup{BC}}+
\left(\Theta\circ\pi\circ(\hat\pi^{\textup{BC}})^{-1}-\tB\right) \right \rangle \\
&\leq \sup_{\tB, \tB'\in \caliB(k,\ell)}\sup_{\pi \in \Pi(d,d)}\left\langle {\tB'-\tB \over \FnormSize{}{\tB-\tB}} \circ \pi, \ \tE \right \rangle+\FnormSize{}{\Theta\circ \pi-\tB\circ \hat\pi^{\textup{BC}}}\\
&\leq \sup_{\tB, \tB'\in \caliB(k,\ell)} \sup_{\pi\in \Pi(d,d)}\left\langle {\tB'-\tB \over \FnormSize{}{\tB-\tB}}\circ \pi, \tE \right \rangle+\underbrace{\FnormSize{}{\Theta\circ \pi-\Theta\circ \hat\pi^{\textup{BC}}}}_{\textup{Permutation error \eqref{eq:permerror}}}+\underbrace{\FnormSize{}{\Theta-\tB}}_{\textup{Proposition~\ref{lem:approx}}}
\end{align}
The first inner product term can be bounded using the same argument in the proof of Theorem~\ref{thm:LSE}. By counting the degree-of-freedom in the space $\{(\tB'-\tB)\colon \tB, \tB'\in\caliB(k,\ell)\} \subset\bbR^{2(\ell+m)^\ell k^m}$ and $\pi\in\Pi(d,d)$, we obtain
\begin{align}
    \mathbb{P}\left(\sup_{\tB, \tB'\in \caliB(k,\ell)} \left\langle {\tB'-\tB \over \FnormSize{}{\tB'-\tB}}, \tE \right \rangle\geq t\right) \lesssim \exp\left(-\frac{t^2}{\sigma^2} +k^m(\ell+m)^\ell+d\log d \right),
\end{align}
Setting $t \asymp \sigma\sqrt{k^m(\ell+m)^\ell+d\log d}$ gives
\begin{align}\label{eq:nonparae}
    \sup_{\tB, \tB'\in \caliB(k,\ell)} \left\langle {\tB'-\tB \over \FnormSize{}{\tB'-\tB}}, \tE \right \rangle\lesssim \sigma\sqrt{k^m(\ell+m)^\ell+d\log d},
\end{align}
with high probability.

Finally, combining all sources of error from Proposition~\ref{lem:approx} and inequalities \eqref{eq:permerror}, \eqref{eq:nonparae}, \eqref{eq:3decomp} yields
\begin{align}\label{eq:three2}
    \frac{1}{d^m} \FnormSize{}{\hat\Theta^{\text{BC}}\circ \hat\pi^{\text{BC}}-\Theta\circ\pi}\lesssim \left(\sigma^2{\log d\over d^{m-1}}\right)^{\beta{\min(\alpha,1)}}+ \sigma^2 \frac{k^m(\ell+m)^\ell}{d^m}+ \frac{L^2}{k^{2\min(\alpha,\ell+1)}}.
\end{align}
Here the term $d\log d$ in~\eqref{eq:nonparae} is absorbed into first term of~\eqref{eq:three2} due to the property $\beta\min(\alpha,1)\leq 1$.
Finally, optimizing~\eqref{eq:three2} with respect to $(k,l)$ gives that 
\begin{align}
     \eqref{eq:three2}\ \lesssim 
     \begin{cases} 
    L^2\left({\sigma\over L}\right)^{4\alpha\over m+2\alpha} d^{-\frac{2m\alpha}{m+2\alpha}}, & \text{ when } \alpha <c(\alpha,\beta,m),\\
    \left( \sigma^2 \log d\over d^{m-1}\right)^{\beta\min(\alpha,1)}, &\text{ when } \alpha \geq c(\alpha,\beta,m),
    \end{cases}
\end{align}
under the choice
\[
\ell^* = \min\left(\lceil\alpha-1\rceil,\lfloor c(\alpha,\beta,m)\rfloor\right),\quad k^* = c_1d^{m/ (m+2\min(\alpha,\ell^*+1))},
\]
where $c(\alpha,\beta,m):= \frac{m(m-1)\beta\min(\alpha,1)}{2(m-(m-1)\beta\min(\alpha,1))}$.

\end{proof}

\subsection{Proof of Theorem~\ref{thm:mminimax}}\label{sec:genminimax}
\begin{proof}[Proof of Theorem~\ref{thm:mminimax}]
For the case where $\alpha = \beta = 1$, we have $\text{Rate}(d)/\log d \leq d^{-\frac{2m}{m+2}}$. Thus, it suffices to show that the estimation problem under the Gaussian model~\eqref{eq:gmd} satisfies the minimax lower bound
\begin{align}\label{eq:nonpara2}
    \inf_{(\hat\Theta,\hat\pi)}\sup_{\Theta\in\tP_{M},\pi\in\Pi(d,d)}\mathbb{P}\left(\textup{MSE}(\hat\Theta\circ\hat\pi,\ \Theta\circ\pi) \gtrsim d^{-{2m\over m+2}} \right) \geq 0.8.
\end{align}
Since $d^{-\frac{2m}{m+2}}$ is the nonparametric rate established in the proof of Theorem~\ref{thm:minimax}, we only need to show that this rate remains the same under additional Lipschitz monotonicity or isotonicity constraints. 

We introduce the following non-negative function
\begin{align}
    K(x) = (1-2|x|)\mathds{1}\{|x|\leq 1/2\}.
\end{align}
It is easy to see that the function $K$ is Lipschitz smooth over $x\in\mathbb{R}$. For any $\ma = (a_1,\ldots,a_m)\in [k]^m$, define the function
\begin{align}\label{eq:phi}
    \phi_{\ma}(x_1,\ldots,x_m) =  {1\over km} \prod_{i\in[m]}K\left(kx_i-a_i+\frac{1}{2}\right).
\end{align}
The function $\phi_{\ma}(x_1,\ldots,x_m)$ is supported on $\bigtimes_{i\in[m]}[{a_i-1 \over k}, {a_i\over k}]$. By the construction of the function $\phi_{\ma}$, we have the following Lemma~\ref{lem:phi}. The proof of Lemma is deferred to Section~\ref{sec:tech}.
\begin{lem}\label{lem:phi}
For any $\ma\in[k]^m$, the function \eqref{eq:phi} satisfies
\begin{enumerate}
    \item $\phi_{\ma}(x_1,\ldots,x_m)\in\tF(1,2)$.
    \item $\sum_{\omega \in[d]^m}\phi^2_{\ma}(\omega /d)\geq C(m)d^mk^{-2-m}$ for some constant $C(m)>0$. 
\end{enumerate}
\end{lem}

Now we prove~\eqref{eq:nonpara2}. For ease of presentation, we take $L=7$; otherwise, we can globally scale the functions in $\tP_M$. Let us take $k \asymp d^{m\over 2+m}$, and let $\Gamma= \{0,1\}^{k^m}$ be the set of all binary sequences of length $ k^m$, and denote $\ma=(a_1,\ldots,a_m)\in[k]^m$. For any collection $\{\gamma_{\ma} \in \{0,1\}|\ma\in[k]^m \}=:\gamma\in\Gamma$, we define the function $h^{\gamma}$ by 
\begin{align}
    h^{\gamma}(x_1,\ldots,x_m) =\sum_{\ma \in [k]^m}\gamma_{\ma}\phi_{\ma}(x_1,\ldots,x_m).
\end{align}
Note that the functions $\phi_{\ma}$ and $\phi_{\ma'}$ have disjoint supports for $\ma\neq \ma'\in[k]^m$. Therefore, we have $h^\gamma \in \tF(1,1)$ by Lemma~\ref{lem:phi}. We further define a function by
\[ 
f^{\gamma} = 5g+  h^{\gamma},\quad \text{where}\quad g(x_1,\ldots,x_m) =\sum_{i = 1}^mx_i.
\]
Note that $g\in \tF(1,1) \cap \tM$, where $\tM$ can be either Lipschitz monotonic class~\eqref{eq:monotonic} or the isotonic class~\eqref{eq:iso}. Next, we show that $f^\gamma \in \tF(1,7)\cap \tM$ for all $\gamma\in \Gamma$. Specifically, for any two points $(x_1,x_2,\ldots,x_m)$ and $(x'_1,x_2,\ldots,x_m)$ with $x_1< x'_1$, we have
\begin{align}\label{eq:differ}
& f^\gamma(x'_1,x_2,\ldots,x_m)-f^\gamma(x_1,x_2,\ldots,x_m) \notag \\
 = &\ 5(x'_1-x_1)+{1 \over k} \sum_{\ma\in[k]^m} \left[K\left(kx'_1-a_1+{1\over 2}\right)-K\left(kx_1-a_1+{1\over 2}\right)\right]\gamma_{\ma} \prod_{i=2}^m K\left(kx_i-a_i + {1\over 2}\right) \notag \\
 \stackrel{(*)}{\geq} &\  5(x'_1-x_1)-(x'_1-x_1) 4\max_{\ma\in[k]^m}\gamma_{\ma}\notag \\
\geq &\ 5(x'_1-x_1)- 4 (x_1'-x_1) = (x'_1-x_1).
 \end{align}
 The line (*) holds due to the disjointness of the supports $\phi_{\ma}$ and $\phi_{\ma'}$ for $\ma\neq \ma'$, and the Lipschitz smoothness of function $K(\cdot)$: 
 \begin{align}
 K\left(kx'_1-a_1+{1\over 2}\right)-K\left(kx_1-a_1+{1\over 2}\right) \geq & -2k|x'_1-x_1| \times \\
 &\quad \mathds{1}\{(a_1-1)\leq kx'_1\leq a_1 \text{ or }(a_1-1)\leq kx_1\leq a_1\}.
 \end{align}
The inequality~\eqref{eq:differ} shows that $f^\gamma\in \tM$, where $\tM$ can be either the Lipschitz monotonic class~\eqref{eq:monotonic} or the isotonic class~\eqref{eq:iso}. Furthermore, the additivity of Lipschitz smoothness implies that $f^{\gamma} \in \tF(1,7)\cap \tM$. 

Now we prove~\eqref{eq:nonpara2} using the subspace $\tF^{\text{sub}}=\{f^{\gamma}\colon \gamma\in\Gamma\}\subset \tF(1,7)\cap \tM$. To apply Lemma~\ref{prop:minmax}, we upper bound $\sup_{f,f'\in \tF^{\text{sub}}}D(\mathbb{P}_f,\mathbb{P}_{f'})$ and lower bound the packing number of $\tF^{\text{sub}}$. Let $D(\cdot | \cdot)$ denote the KL divergence defined in~\eqref{eq:KL}. Under the Gaussian noise model, we have that, 
 \begin{align}
 \sup_{f,f'\in\tF^{\text{sub}}}  D(\mathbb{P}_f|\mathbb{P}_{f'})&\leq \frac{1}{2\sigma^2} \sup_{f,f'\in\tF^{\text{sub}}} \sum_{\omega\in[d]^m}  \left(f\left({\omega/ d}\right)-f'\left({\omega/ d}\right)\right)^2\\
 &\leq \frac{1}{2\sigma^2} \sup_{h,h'\in \{h^\gamma:\gamma\in \Gamma\}} \sum_{\omega\in[d]^m}\left(h\left({\omega/ d}\right)-h'\left({\omega/d}\right)\right)^2\\&\leq \frac{1}{2\sigma^2}d^m k^{-2}.
\end{align}
where the first inequality holds by Lemma~\ref{prop:kl}.

Next we lower bound the packing number of $\tF^{\text{sub}}$. For any $f^\gamma,f^{\gamma'}\in\tF^{\text{sub}}$, we have
\begin{align}
    \rho^2(f^{\gamma},f^{\gamma'})&\geq \frac{1}{d^m}\sum_{\omega\in[d]^m}\sum_{\ma\in[ k]^m}(\gamma_{\ma}-\gamma'_{\ma})^2\phi_{\ma}^2(\omega/d)\\
    &\geq C(m) k^{-2-m}\rho_H(\gamma,\gamma'),
\end{align}
where the last inequality uses Lemma~\ref{lem:phi} and $\rho_H$ is defined in Lemma~\ref{lem:covering}. By Lemma~\ref{lem:covering}, we can select a subset $S\subset \Gamma$ such that $|S|\geq \exp(k^m/16)$ and $\rho_H(\gamma,\gamma')\geq k^m/8$ for any $\gamma\neq \gamma'\in S$. Therefore setting $\epsilon^2 = C(m)k^{-2}$, we have
\begin{align}
    \log \tM(\epsilon,\tF^{\text{sub}},\rho)\geq k^m/16.
\end{align}
The choice of $k\asymp d^{m\over m+2}$ ensures $\sup_{f,f'\in \tF^{\text{sup}}}D(\mathbb{P}_f | \mathbb{P}_{f'})\lesssim \log \tM(\varepsilon,\tF^{\text{sup}},\rho)$. By Lemma~\ref{prop:minmax}, we have
\begin{align}
    \inf_{\hat\Theta}\sup_{f\in\tF(1,L)\cap \tM}\mathbb{P}\left(\frac{1}{d^m}\sum_{\omega\in[d]^m}\left(\hat\Theta(\omega)-f(\omega/d)\right)^2 \gtrsim d^{-{2m\over m+2}} \right) \geq 0.8,
\end{align}
which completes the proof of \eqref{eq:nonpara2}.
 \end{proof}

\section{Extension of Theorem~\ref{thm:BC} to isotonic functions}\label{sec:isotonic}
We provide an upper bound under isotonic assumption similar to Theorem~\ref{thm:BC}. Define isotonic function class $\tM$ 
\begin{align}
    \tM = \{f\colon[0,1]^m\rightarrow\mathbb{R}\ \big| \ f(x_1,\ldots,x_m)\leq f(x_1',\ldots,x_m')\text{ when } x_i\leq x_i' \text{ for } i\in[m]\}.
\end{align}
Unlike $\beta$-monotonic functions, the permutation $\pi$ per se may not be accurately estimated for isotonic functions. However, we find that the composition $\Theta\circ \pi$ can still be accurately estimated due to the joint monotonicity. 
\begin{lem}[Permutation error of Borda count algorithm under isotonic assumption]\label{lem:isotonic}
Consider the sub-Gaussian tensor model~\eqref{eq:obs} with $f\in\tM$. Let $\hat\pi$ be the permutation such that the permuted empirical score function $\tau\circ\hat\pi^{-1}$ is monotonically increasing. Then with high probability
\begin{align}
    \frac{1}{d^m}\|\Theta\circ\hat\pi-\Theta\circ\pi\|_F^2\lesssim \frac{\log d}{d^{m-1}}.
\end{align}
\end{lem}
The proof of Lemma~\ref{lem:isotonic} is provided in Section~\ref{sec:tech}. 
Now, we obtain the same statistical accuracy of the Borda count estimator for the isotonic functions as in Theorem~\ref{thm:BC}. 
\begin{thm}[Estimation error for Borda count algorithm under isotonic assumption]\label{thm:BC2} Consider the sub-Gaussian tensor model with $f\in \tF(\alpha,L)\cap \tM$.
Let $(\hat\Theta^{\textup{BC}},\hat\pi^{\textup{BC}})$ be the Borda count estimator in \eqref{eq:permute}-\eqref{eq:bclse} with {\small $\ell^* = \min(\lfloor \alpha\rfloor,(m-2)(m+1)/2)$} and $k^*\asymp d^{m/(m+2\min(\alpha,\ell^*+1))}$. Then, we have in high probability,
\begin{align}\label{eq:rates2}
     \textup{MSE}(\hat \Theta^{\textup{BC}}\circ\hat \pi^{\textup{BC}},\Theta\circ \pi
     ) \lesssim 
     \begin{cases} 
     d^{-\frac{2m\alpha}{m+2\alpha}}, & \text{ when } \alpha < m(m-1)/2,\\
     d^{-(m-1)}\log d, &\text{ when } \alpha \geq m(m-1)/2.
    \end{cases}
\end{align}
\end{thm}

\begin{proof}[Proof of Theorem~\ref{thm:BC2}]
The proof of Theorem~\ref{thm:BC2} is the exactly same as in Theorem~\ref{thm:BC}, except that we now use Lemma~\ref{lem:isotonic} in place of Lemma~\ref{lem:permute}. We omit the proof for brevity. 
\end{proof}

\section{Challenges for extending Theorem~\ref{thm:mminimax} to arbitrary~$\alpha$}\label{subsec:generala}
Our Theorem~\ref{thm:mminimax} is based on Lipschitz monotone ($\alpha = \beta=1$). Here we provide a high level explanation for technical challenges of constructing minimax lower bound under arbitrary $\alpha>0$ and extra monotonicity. When $\alpha \leq 1$, the nonparametric rate is dominant. For general $\alpha$, both the nonparametric and permutation rates must be controlled.

{\bf Challenge in nonparametric rates.} In the general proof for the nonparametric rate, we partition the domain into $k^m$ uniform blocks and consider a set of $\alpha$-smooth functions $\{\phi_{\mathbf{a}}: \mathbf{a} \in [k]^m\}$, each supported on one of these partitions within $[0,1]^m$. This construction is feasible when $0 < \alpha \leq 1$. To ensure marginal $\beta$-monotonicity, we add a simple linear function $g(x_1, \ldots, x_m) = \sum_{i \in [m]} x_i$. For $g$ to dominate the $\alpha$-smooth function, we require $\alpha = 1$. For arbitrary $\alpha \neq 1$, constructing such function classes remains an open problem. 

{\bf Challenge in permutation rate.} For permutation rates, we construct a smooth function that mimics a block-wise constant tensor over $d/2$ points, reducing the problem to estimating a block-constant tensor. For monotonic $\alpha$-smooth functions, constructing subregions that mimic marginally monotonic tensors complicates the derivation since these tensors are not block-constant. The monotonicity constraint also prevents the use of smooth convolution functions. Extending this approach to general $\alpha$ and $\beta$ remains an open problem.

\section{Extension of Section~\ref{sec:lse} to random designs}\label{sec:rdesign}
We now extend the grid designs~\eqref{eq:rep} to general designs~\eqref{eq:randommodel}. Specifically, suppose that the signal tensor is generated from the following model
    \begin{align}\label{eq:random}
        \Theta(i_1,\ldots,i_m) = f(x_{i_1},\ldots,x_{i_m}),\quad \text{for all}\quad (i_1,\ldots,i_d) \in [d]^m,
    \end{align}
    where the design points $\{x_i\}_{i\in[d]}$ are random variables sampled from distribution $\mathbb{P}_x$. We write $\{x_i\}_{i\in[d]}\sim \mathbb{P}_x$ and require $\mathbb{P}_x$ to be in the following class:
    \begin{equation}\label{eq:Px}
    \tP_x=\left\{\mathbb{P}_x: \mathbb{P}_{x}\left(\sum_{i\in[d]} \mathds{1}\left\{x_i\in\left[{a-1\over k},{a\over k}\right)\right\} \right)\asymp {d\over k}, \text{ for all $a\in[k]$}\} \geq 1 - \exp(-d^\delta)\right\},
    \end{equation}   
    for some arbitrary small constant $\delta\in(0,1)$. This distribution class has been considered in earlier work for smooth graphons~\citep{gao2015rate}. Intuitively speaking, each interval $[{a-1\over k},{a\over k})$ contains roughly $d/k$ observations. Note that both the fixed balanced design $\{1/d,2/d,\ldots,d/d\}$ and the i.i.d.\ uniform distribution belong to this class $\tP_x$.
 
 {\bf Statistical results.}  We extend Theorems~\ref{thm:LSE} and~\ref{thm:minimax} to the above setting. We redefine the $\tP(\alpha,L)$ in~\eqref{eq:parameterP}, $\caliB(k,\ell)$ in~\eqref{eq:polynomial}, and LSE in~\eqref{eq:lseopt} as follows to allow general designs:
 \begin{align}
& \tP^{\text{gen}}(\alpha,L)=\{\Theta: \Theta \text{ is from~\eqref{eq:random} with latent variable $\{x_i\}_{i\in[d]}$ and $f\in \tF(\alpha,L)$}\},\\
 &\caliB^{\text{gen}}(k,\ell)=\bigcup_{z\in \Pi(d,k)}\bigg\{\tB\colon \tB(\omega)=\sum_{(j_1,\ldots,j_m)\in[k]^m}\text{Poly}_{\ell,\Delta}(\omega)\mathds{1}\{z\circ \omega=(j_1,\ldots,j_m)\}\bigg\},\\
 &(\hat \Theta^{\text{LSE,gen}},\hat \pi^{\text{LSE,gen}})=\argmin_{\Theta\in\caliB^{\text{gen}}(k,\ell),\pi\in\Pi(d,k)\}}\FnormSize{}{\tY-\Theta\circ \pi},
  \end{align}
 where in the definition of $\caliB^{\text{gen}}(k,\ell)$ we have used $z \in \Pi(d,k)$ to denote any clustering from $[d]$ to $[k]$ which may be unbalanced, the shorthand notation $\omega:=(i_1,\ldots,i_m)$, and $z\circ \omega := (z(i),\ldots,z(i_m))$. The upper script ``gen'' denotes the general designs. 
  
 \begin{thm}[Extension of Theorems~\ref{thm:LSE} and~\ref{thm:minimax} for general designs]\label{lem:approx_random}
Consider the general design~\eqref{eq:random} and same conditions as in Theorems~\ref{thm:LSE} and~\ref{thm:minimax}. Assume the design points $\{x_i\}_{i\in[d]}\sim \mathbb{P}_x$ where $\mathbb{P}_x$ is in the distribution class $P_x$~\eqref{eq:Px}. 
\begin{itemize}
\item Upper bound: 
\[
\mathbb{E}\textup{MSE}(\hat \Theta^{\textup{LSE,gen}}\circ\hat \pi^{\textup{LSE,gen}},\Theta\circ\pi)\leq \textup{Rate}(d),
\] 
where the expectation is with respect to both $\{x_i\}_{i\in[d]}\sim \mathbb{P}_x$ and the observation noise. 
\item Lower bound: 
\[
\inf_{(\hat \Theta,\hat \pi)}\sup_{(\Theta,\pi)\in \tP^{\textup{gen}}(\alpha,L) \times \Pi(d,d)}  \sup_{\mathbb{P}_{x}\in \tP_x} \mathbb{P} \left( \textup{MSE}(\hat \Theta \circ \hat \pi, \Theta\circ \pi)\gtrsim \textup{Rate}(d)\right)\geq 0.8.
\]
\end{itemize}
    \end{thm}

    \begin{proof}[Proof of Theorem~\ref{lem:approx_random}]
The proof for upper bound is similar to that in Theorem~\ref{thm:LSE}, except that we now assess the mean-squared-error conditional on $\{x_i\}_{i\in[d]}$. Theorem~\ref{thm:LSE} uses Proposition~\ref{lem:approx} and the union bound over the set $(\Theta,\pi)\in \tP(\alpha,L)\times \Pi(d,k)$. Now we first extend the conclusion of Proposition~\ref{lem:approx} to allow general designs; specifically, the following approximation error holds
\[
\inf_{\tB\in \caliB^{\textup{gen}}(k,\ell)} {1\over d^m}\FnormSize{}{\Theta-\tB}^2\leq {L^2\over k^{2\min(\alpha,\ell+1)}}.
\]
To see the above line, recall that the original proof of Proposition~\ref{lem:approx} uses the equal-sized clustering such that $z(i)=\lceil ki /d\rceil$ for all $i\in[d]$. To accommodate random designs, we use the clustering function $z^\text{gen}\colon[d]\rightarrow[k]$ conditional $\{x_i\}_{i\in[d]}$ such that
    \begin{align}
        z^{\text{gen}}(i) = \lceil k*x_i\rceil.
    \end{align}
The proof of Proposition~\ref{lem:approx} carries over with $z$ replaced by $z^\text{gen}$. Second, we find that the union bound over the new set $(\Theta,\pi)\in \tP^{\text{gen}}(\alpha,L)\times \Pi(d,d)$ has the same order of degree-of-freedom compared to the original set $(\Theta,\pi)\in \tP(\alpha,L)\times \Pi(d,d)$. This is because the extra degree-of-freedom in unknown $\{x_i\}_{i\in[d]}$ is dominated by the existing complexity of $\pi\in\Pi(d,d)$, and thus has no effect in the final complexity.

The lower bound holds trivially because the fixed balanced designs $\{1/d,2/d,\ldots,d/d\}$ is a spacial case of the distribution class $\tP_x$ in~\eqref{eq:Px}.
\end{proof}   
 {\bf Computational results.} The extension of computational results to random designs is challenging due to the uncertainty in the design points. Although polynomial-time algorithms have been proposed under certain scenarios~\citep{balasubramanian2021nonparametric}, their optimality remains unknown. Additionally, our Borda count algorithm uses equal-sized blocks, which is suitable for balanced designs. It is an open question whether the Borda count algorithm can achieve computational optimality for random designs. We leave this question for future research.

\section{Technical lemmas}\label{sec:tech}

\begin{proof}[Proof of Lemma~\ref{lem:permutation}]
We provide the proof for $m=3$ only. The extension to higher orders $(m\geq 4)$ uses exactly the same techniques and thus is omitted. For notational simplicity, we do not distinguish the fractional number and its rounding integer. For example, we simply write $k/2$, instead of $\lceil k/2\rceil$, to represent its rounding integer. 

Let us pick $\gamma_1,\ldots,\gamma_{k/2}\in\{0,1\}^{k^2/4}$ such that $\rho_H(\gamma_p,\gamma_q)\geq k^2/16$ for all $p\neq q\in [k/2]$. This selection is possible by lemma~\ref{lem:covering}.
Fixing such $\gamma_1,\ldots,\gamma_{k/2}$, we define a core tensor $\tS\in\bbR^{k\times k\times k}$ as follows:
\begin{align}\label{eq:structure}
\tS([k/2],[k/2],r)&=\sqrt{\sigma^2\log k \over d^2} \times \textup{Mat}(\gamma_r),\quad \text{for all }r\in[k/2], \notag \\
\tS(\{k/2+1,\ldots,k\},\{k/2+1,\ldots,k\},r)&=\sqrt{\sigma^2\log k \over d^2} \times I_{(k/2)\times (k/2)},\quad \text{for all }r\in\{k/2+1,\ldots,k\}, \notag \\
\tS(i,j,r)&=0, \ \text{others}.
\end{align}
where $\text{Mat}(\gamma)$ is the operation to reshape the vector $\gamma\in \{0,1\}^{k^2/4}$ into a $(k/2)$-by-$(k/2)$ matrix, and $I_{(k/2)\times (k/2)}$ represents the $(k/2)$-by-$(k/2)$ identity matrix. Notice that for any $p\neq q\in[k/2]$, the corresponding slices in the tensor $\tS$ obey the lower bound:
\begin{align}\label{eq:sr}
    \FnormSize{}{\tS(:,:,p)-\tS(:,:,q)}^2\gtrsim \frac{\sigma^2k^2\log k}{d^2},
\end{align}

Define a subset of permutation set $\Pi(d,k)$ by
\begin{align}
    \tZ = \left\{z\in \Pi(d,k)\colon |z^{-1}(a)| = \frac{d}{k} \text{ for all } a\in[k], z^{-1}(a) = \left\{\frac{(a-1)d}{k}+1,\ldots,\frac{ad}{k}\right\}  \text{ for } a\in[k/2] \right\}.
\end{align}
Each $z\in\tZ$ induces a block structure in $\tB(k,0)$, with known indices for first $k/2$ clusters and unknown indies for the remaining $k/2$ clusters. We consider the collection of block tensors induced by $\tZ$ and $\tS$; i.e.,
\begin{align}
    \tB(\tZ) = \{\Theta^z\in\bbR^{d\times d\times d}\colon \Theta^z(i,j,k) = \tS(z(i),z(j),z(k)) \text{ for } z\in\tZ\}.
\end{align}
To apply Lemma~\ref{prop:minmax}, we seek for the upper bound $\sup_{\Theta,\Theta'\in\tB(\tZ)}D(\mathbb{P}_\Theta|\mathbb{P}_{\Theta'})$ and the lower bound for $\log \tM(\epsilon,\tB(\tZ),\rho)$, where $\tM(\cdot,\cdot,\cdot)$ denotes the packing number of $\tB(\tZ)$ under metric defined by
$\rho(\Theta,\Theta') = \frac{1}{d^3}\FnormSize{}{\Theta-\Theta'}^2.$
Let $D(\cdot| \cdot)$ denote the KL divergence defined in~\eqref{eq:KL}. Under the Gaussian noise model, we have
\begin{align}\label{eq:kldistance}
    D(\mathbb{P}_\Theta|\mathbb{P}_{\Theta'})\lesssim \frac{1}{\sigma^2}\FnormSize{}{\Theta-\Theta'}^2\lesssim \frac{1}{\sigma^2}d^3 \frac{\sigma^2\log k}{d^2}=d\log k,
\end{align}where the first inequality holds for any $\Theta,\Theta'\in\tB(\tZ)$ by Lemma~\ref{prop:kl}. 

Next we provide a lower bound of the packing number $\log \tM(\epsilon,\tB(\tZ),\rho)$ with $\epsilon^2 \asymp {\sigma^2\log k \over d^2}$. From the construction of $\tS$ in \eqref{eq:structure}, we have one-to-one correspondence between $\tZ$ and $\tB(\tZ)$. Thus $\tM(\epsilon,\tB(\tZ),\rho) = \tM(\epsilon,\tZ,\rho')$ for the metric $\rho'$ on $\tZ$ defined by $\rho'(z_1,z_2) = \rho(\Theta^{z_1},\Theta^{z_2})$. Let $P$ be the packing set in $\tZ$ with the same cardinality of $\tM(\epsilon,\tZ,\rho').$ Given any $z\in\tZ$, define its $\epsilon$-neighbor by $\tN(z,\epsilon) = \{z'\in\tZ\colon \rho'(z,z')\leq \epsilon\}.$ Then, we have $\cup_{z\in P} \tN(z,\epsilon) = \tZ$, because the cardinality of $P$ is same as packing number $\tM(\epsilon,\tZ,\rho')$. Therefore, we have
\begin{align}\label{eq:packing}
    |\tZ|\leq\sum_{z\in P}|\tN(z,\epsilon)|\leq |P|\max_{z\in P}|\tN(z,\epsilon)|.
\end{align}
It remains to find the upper bound of $\max_{z\in P}|\tN(z,\epsilon)|.$ For any $z_1,z_2\in\tZ$, $z_1(i) = z_2(i)$ for $i\in[d/2]$ and $|z_1^{-1}(p)| = |z_2^{-1}(q)|=d/k$ for all $p,q\in[k].$ Therefore,
\begin{align}
    \rho'^2(z_1,z_2) &\geq \frac{1}{d^3}\sum_{i_1,i_2\in[d/2], d/2<i_3\leq d}\left( \tS\left(z_1(i_1),z_1(i_2),z_1(i_3)\right)-\tS\left(z_2(i_1),z_2(i_2),z_2(i_3)\right)\right)^2\\
    &= \frac{1}{d^3}\sum_{d/2< i_3\leq d}\sum_{p,q \in[k/2]}\sum_{i_1\in z_1^{-1}(p),i_2\in z_1^{-1}(q)}\left( \tS\left(p,q,z_1(i_3)\right)-\tS\left(p,q,z_2(i_3)\right)\right)^2\\
    &=
    \frac{1}{d^3}\sum_{d/2<i_3\leq d}\sum_{p,q\in[k/2]}\left(d\over k\right)^2\left( \tS\left(p,q,z_1(i_3)\right)-\tS\left(p,q,z_2(i_3)\right)\right)^2\\&=
    \frac{1}{d^3}\sum_{d/2<i_3\leq d}\left(d\over k\right)^2\FnormSize{}{\tS(:,:,z_1(i_3))-\tS(:,:,z_2(i_3))}^2
    \\&\gtrsim  \frac{\sigma^2\log k}{d^3}|\{j\colon z_1(j)\neq z_2(j)\}|,
\end{align}where the last inequality is from \eqref{eq:sr}. Hence with the choice of $\epsilon^2 \asymp \frac{\sigma^2\log k}{d^2}$, we have $|\{j\colon z(j)\neq z'(j)\}|\leq d/6$ for any $z'\in\tN(z,\epsilon)$.
This implies
\begin{align}\label{eq:fpacking}
    \tM(\epsilon,\tB(\tZ),\rho) = |P| \geq\frac{|\tZ|}{\max_{z\in P}\tN(z,\epsilon)} \geq {\frac{(d/2)!}{[(d/k)!]^{k/2}} \over {d\choose d/6}k^{d/6}}\gtrsim \exp\left(d\log k\right).
\end{align}
Finally, applying Lemma~\ref{prop:minmax} based on  \eqref{eq:kldistance} and \eqref{eq:fpacking} gives
\begin{align}
    \inf_{\hat\Theta}\sup_{\Theta\in\tB(\tZ)}\mathbb{P}\left(\frac{1}{d^3}\FnormSize{}{\hat\Theta-\Theta}^2\gtrsim \frac{\sigma^2\log k}{d^2}\right)=\inf_{\hat\Theta}\sup_{z\in\tZ}\mathbb{P}\left(\frac{1}{d^3}\FnormSize{}{\hat\Theta-\tS\circ z}^2\gtrsim \frac{\sigma^2\log k}{d^2}\right)\geq 0.9.
\end{align}
\end{proof}

\begin{proof}[Proof of Lemma~\ref{lem:permute}]
Without loss of generality, assume that $\pi$ is the identity permutation and the tensor $\Theta$ is symmetric. 
Notice that $g(i)-\tau(i)$ is the sample average of $d^{m-1}$ independent mean-zero Gaussian random variables with the variance proxy $\sigma^2$. Based on the maxima of sub-Gaussian random variables, we have
\begin{align}\label{eq:concentration}
   \max_{i\in[d]} |g(i)-\tau(i)| < 2\sigma d^{-(m-1)/2}\sqrt{\log d},
\end{align}
with probability $1-\frac{2}{d^2}$. 

By the $\beta$-monotonicity of the function $g$, we have
\begin{align}\label{eq:mon1}
    g(1)\leq g(2)  \leq \cdots \leq g(d-1),
\end{align}
The estimated permutation $\hat\pi$ is defined for which
\begin{align}\label{eq:mon2}
    \tau \circ \hat\pi^{-1}(1) \leq\tau \circ \hat\pi^{-1}(2) \leq \cdots\leq \tau \circ \hat\pi^{-1}(d-1) \leq \tau \circ \hat\pi^{-1}(d).
\end{align}

For any given index $i$, we examine the error $|i-\hat\pi(i)|$. By \eqref{eq:mon1} and \eqref{eq:mon2}, we have
\begin{align}
    i = \underbrace{|\{j\colon g(j)\leq g(i)\}|}_{=:\textup{I}},\quad\text{and}\quad \hat\pi(i) = \underbrace{|\{j\colon \tau(j)\leq \tau(i)\}|}_{=:\textup{II}},
\end{align}
where $|\cdot|$ denotes the cardinality of the set. We claim that the sets I and II differ only in at most $d^{(m-1)\beta/2}$ elements. To prove this, we partition the indices in $[d]$ in two cases.
\begin{enumerate}
    \item Long-distance indices in $\{j\colon|j-i|/d\geq 
C \left(\sigma d^{-(m-1)/2}\sqrt{\log d}\right)^{\beta} \}$ for some sufficient large constant $C>0$. In this case, the ordering of $(i,j)$ remains the same in \eqref{eq:mon1} and \eqref{eq:mon2}, i.e.,
    \begin{align}\label{eq:equiv}
        g(i) < g(j) \Longleftrightarrow \tau(i) < \tau(j).
    \end{align}
    We only prove the right side direction in  \eqref{eq:equiv} here. The other direction can be similarly proved.
    Suppose that $g(i) < g(j)$. Then we have
    \begin{align}
        \tau(j)-\tau(i) &\geq -|g(j)-\tau(j)|-|g(i)-\tau(i)| + g(j)-g(i)
       \\& > -4\sigma d^{-(m-1)/2}\sqrt{\log d} + g(j)-g(i)
       \\&\geq 0,
    \end{align}
    where the second inequality is from \eqref{eq:concentration} with probability at least $(1-2/d^2)^d$ and the last inequality uses $\beta$-monotonicity of $g(\cdot)$, and the assumption
    $|j-i|/d\geq C\left(\sigma d^{-(m-1)/2}\sqrt{\log} d\right)^{\beta}$. Therefore we show that $g(i)<g(j)$ implies $\tau(i)<\tau(j).$
    In this case, we conclude that none of long-distance indices belongs to I$\Delta$II.
    \item Short-distance indices in $\{j\colon|j-i|/d< \left(\sigma d^{-(m-1)/2}\sqrt{\log d}\right)^{\beta} \}$. In this case, \eqref{eq:mon1} and \eqref{eq:mon2} may yield different ordering of $(i,j).$
\end{enumerate}
Combining the above two cases gives that 
\begin{align}
   \left \{j\colon\frac{1}{d}|j-i|\leq \left(4\sigma d^{-(m-1)/2}\sqrt{\log d}\right)^{\beta}\right \}\supset \textup{I}\Delta\textup{II}.
\end{align}
Finally, we have
\begin{align}
    \textup{Loss}(\pi,\hat\pi) := \frac{1}{d}\max_{i\in[d]}|\pi(i)-\hat\pi(i)|\leq\frac{1}{d}|\textup{I}\Delta\textup{II}|\leq 2\left(4\sigma d^{-(m-1)/2}\sqrt{\log d}\right)^{\beta},
\end{align}
with high probability.
\end{proof}

 \begin{proof}[Proof of Lemma~\ref{lem:phi}]
 We prove the two properties separately. 
 \begin{enumerate}
 \item $\phi_{\ma}\in\tF(1,2):$
 We first claim that 
 \begin{align}\label{eq:claim}
     \left|\prod_{i\in[m]}K(x_i)-\prod_{i\in[m]}K(y_i)\right|\leq 2\sum_{i\in[m]}|x_i-y_i|.
 \end{align}
 We prove this by induction. When $m = 1$, we have 
 \begin{align}\label{eq:as}
     |K(x_1)-K(y_1)|&\leq 2|x_1-y_1|.
 \end{align}
 Suppose~\eqref{eq:claim} holds for $\ell \in\{1,2,\ldots,(m-1)\}$. Then, we have
 \begin{align}
      &\left|K(x_{\ell+1})\prod_{i\in[l]}K(x_i) -K(y_{\ell+1}) \prod_{i\in[l]}K(y_i)\right|\\
      \leq &\ \left| K(x_{\ell+1})-K(y_{\ell+1})\right| \prod_{i\in[l]}K(x_i)+K(y_{\ell+1})\left|\prod_{i\in[l]}K(x_i)-\prod_{i\in[l]}K(y_i)\right|\\
      \leq &\ 2|x_{\ell+1}-y_{\ell+1}|+2\sum_{i\in[\ell]}^\ell |x_i-y_i|\\
      \leq &\  2 \sum_{i\in[\ell+1]}|x_i-y_i|,
 \end{align}
where the third line uses \eqref{eq:as} and the induction assumption for $\ell$. By mathematical induction, \eqref{eq:claim} holds true for any integer $m\in\mathbb{N}_{+}$. 

For any $\mx=(x_1,\ldots,x_m),\my = (y_1,\ldots,y_m)$ in the support of $\phi_{\ma}$, we have
\begin{align}
    |\phi_{\ma}(x_1,\ldots,x_m)-\phi_{\ma}(y_1,\ldots,y_m)|&\leq {2\over km} \sum_{i=1}^m|kx_i-ky_i|\\&\leq2 \|\mx-\my\|_\infty.
\end{align}
Therefore $\phi_{\ma}\in\tF(1,2).$

 \item $\sum_{\omega \in[d]^m}\phi^2_{\ma}(\omega/d)\geq C(m)d^mk^{-2-m}$ for some constant $C(m)>0$. We write
 \begin{align}\label{eq:sum}
     \sum_{\omega\in[d]^m}\phi_{\ma}^2\left(\omega/d\right) &= {1\over k^2m^2}\sum_{(i_1,\ldots,i_m)\in[d]^m} K^2\left(\frac{k i_1}{d}-a_1+\frac{1}{2}\right)\times \cdots \times K^2\left(\frac{k i_m}{d}-a_m+\frac{1}{2}\right)
 \end{align}
For any integers $a\in[k]$ and $i\in[d]$, we have
 \begin{align}\label{eq:tem}
\sum_{i\in[d]} K^2\left({ki\over d}-a+{1\over 2}\right) &= \sum_{i\in[d]}\left(1-2\left|{ki \over d} - a +{1\over 2}\right|\right)^2 \mathds{1}\{ d(a-1) \leq ki\leq da\}\\
& = \sum_{\frac{d(a-1)}{k}\leq i\leq\frac{da}{k}}\left(1-\frac{2k}{d}\left|i-\frac{(a-1/2)d}{k}\right|\right)^2\\
&\geq 2\int_0^{d/2k}\left(1-\frac{2k}{d}t\right)^2 dt= \frac{d}{3k}.
 \end{align}
 Plugging~\eqref{eq:tem} into~\eqref{eq:sum} gives the desired conclusion
 \[
 \sum_{\omega\in[d]^m}\phi^2_{\ma}(\omega/d)\gtrsim {d^m \over 3^mk^{m+2}m^2}.
 \]
 \end{enumerate}
 \end{proof}

\begin{proof}[Proof of Lemma~\ref{lem:isotonic}]
Without loss of generality, assume that $\pi$ is the identity permutation. Notice that $g(i)-\tau(i)$ is the sample average of $d^{m-1}$ independent mean-zero sub-Gaussian random variables with the variance proxy $\sigma^2$. Based on the maxima of sub-Gaussian random variables, we have 
\begin{align}
     \max_{i\in[d]}|g(i)-\tau(i)| \lesssim 2\sigma\left( \frac{\log d}{d^{m-1}}\right)^{1/2},
\end{align}
with probability $1-\frac{2}{d^2}$. Denote $\delta=2\sigma\sqrt{\log d / d^{m-1}}$. 
The estimated permutation $\hat\pi$ is defined for which 
\begin{align}\label{eq:mon3}
    \tau \circ \hat\pi^{-1}(1) \leq\tau \circ \hat\pi^{-1}(2) \leq \cdots\leq \tau \circ \hat\pi^{-1}(d-1) \leq \tau \circ \hat\pi^{-1}(d).
\end{align}
By definition of the $\hat\pi$, we have for any $i>j$ but $\hat\pi(i)= j$, 
\begin{align}\label{eq:pi1}
    \frac{1}{d^{m-1}}\sum_{(i_2,\ldots,i_m)\in[d]^{m-1}}\left|\Theta(i,i_2,\ldots i_m)-\Theta(\hat\pi(i),i_2,\ldots,i_m)\right|\lesssim \delta.
\end{align}
Similarly for any $i<j$ but $\hat\pi(i) = j$, we have
\begin{align}\label{eq:pi2}
    \frac{1}{d^{m-1}}\sum_{(i_2,\ldots,i_m)\in[d]^{m-1}}\left|\Theta(\hat\pi(i),i_2,\ldots i_m)-\Theta(i,i_2,\ldots,i_m)\right|\lesssim \delta.
\end{align}
Therefore we obtain
\begin{align}
    &\frac{1}{d^m}\sum_{(i_1,\ldots,i_m)\in[d]^m}|\Theta(\hat\pi(i_1),\hat\pi(i_2),\ldots,\hat\pi(i_m))-\Theta(i_1,i_2\ldots,i_m)|\\&\leq \frac{1}{d^m}\sum_{(i_1,\ldots,i_m)\in[d]^m}\bigg(|\Theta(\hat\pi(i_1),\hat\pi(i_2),\ldots,\hat\pi(i_m))-\Theta(i_1,\hat\pi(i_2),\ldots,\hat\pi(i_m))|\\&\hspace{3cm}+|\Theta(i_1,\hat\pi(i_2),\ldots,\hat\pi(i_m))-\Theta(i_1,i_2,\ldots,\hat\pi(i_m))|\\&\hspace{3cm}+\cdots+|\Theta(i_1,i_2,\ldots,\hat\pi(i_m))-\Theta(i_1,i_2,\ldots,i_m)|\bigg)\\&\leq\frac{m}{d^m}\sum_{i_1\in[d]}\sum_{(i_2,\ldots,i_m)}|\Theta(\hat\pi(i_1),i_2,\ldots,i_m)-\Theta(i_1,i_2,\ldots,i_m)|\\&\lesssim m\delta,
\end{align}
where the second inequality uses the symmetricity of the tensor and the last inequality uses \eqref{eq:pi1} and \eqref{eq:pi2} in the case of wrong permutations. Since the Frobenius norm is bounded by $\ell_1$ norm, we complete the proof.
\end{proof}

\begin{defn}[sub-Gaussian random vectors]
Let $\vnormSize{}{\cdot}$ denote the vector 2-norm.  A random vector $\me \in \mathbb{R}^d$ is said to be sub-Gaussian with variance proxy $\sigma^2$ if $\mathbb{E}\me=0$ and $\langle {\mc\over \vnormSize{}{\mc}}, \me\rangle$ is sub-Gaussian with variance proxy $\sigma^2$ for any vector $\mc\in\mathbb{R}^d$.   Note that the entries in $\me$ need not be independently nor identically distributed. 
 \end{defn}
 
\begin{lem}[Sub-Gaussian maxima under embedding]\label{lem:embedding}
Let $\mA\in\bbR^{d_1\times d_2}$ be a deterministic matrix with rank $r\leq\min(d_1,d_2)$. Let $\my\in\bbR^{d_1}$ be a vector consisting of independent sub-Gaussian entries with variance proxy $\sigma^2$. Then, there exists a sub-Gaussian random vector $\mx\in\bbR^r$  with variance proxy $\sigma^2$ such that  
\begin{align}
    \max_{\mp\in\bbR^{d_2}}\left\langle \frac{\mA\mp}{\|\mA\mp\|_2},\my\right\rangle = \max_{\mq\in\bbR^r}\left\langle \frac{\mq}{\|\mq\|_2},\mx\right\rangle.
\end{align}
\end{lem}

\begin{proof}[Proof of Lemma~\ref{lem:embedding}]
Let $\bmu_i\in\bbR^{d_1},\mv_j\in\bbR^{d_2}$ singular vectors and $\lambda_i\in\bbR$ be singular values of of $\mA$ such that
 $\mA = \sum_{i=1}^r \lambda_i\bmu_i\mv_i^T$.
Then  for any $\mp\in\bbR^{d_2}$, we have 
\begin{align}
    \mA\mp = \sum_{i=1}^r\lambda_i \bmu_i\mv_i^T\mp = \sum_{i=1}^r \lambda_i(\mv_i^T\mp) \bmu_i = \sum_{i = 1}^r \alpha_i\bmu_i,
\end{align}
where $\malpha(\mp) = (\alpha_1,\ldots,\alpha_r)^T := \left(\lambda_1(\mv_1^T\mp),\ldots,\lambda_r(\mv_r^T\mp)\right)^T\in\bbR^{r}$. Notice that $\malpha(\mp)$ covers $\bbR^r$ in the sense that $\{\malpha(\mp)\colon \mp\in\bbR^{d_2}\} = \bbR^r.$
Therefore, we have 
\begin{align}
    \max_{\mp\in\bbR^{d_2}}\left\langle \frac{\mA\mp}{\|\mA\mp\|_2},\my\right\rangle &=\max_{\mp\in\bbR^{d_2}}\sum_{i=1}^r\frac{\alpha_i}{\|\malpha(\mp)\|_2}\bmu_i^T\my\\&=\max_{\mp\in\bbR^{d_2}}\left\langle \frac{\malpha(\mp)}{\|\malpha(\mp)\|_2},\mx\right\rangle\\&=\max_{\mq\in\bbR^r}\left\langle \frac{\mq}{\|\mq\|_2},\mx\right\rangle,
\end{align}
where we define $\mU=[\bmu_1,\bmu_2,\ldots,\bmu_r]\in\mathbb{R}^{d_2\times r}$ and $\mx = \mU^T\my\in \bbR^{r}$. Since $\mU$ has orthonormal columns, the vector $\mx$ is a sub-Gaussian random vector with variance proxy $\sigma^2$. 
\end{proof}

\begin{lem}[Theorem 1.19 in \citet{rigollet2015high}]~\label{lem:subga} Let $\me\in\mathbb{R}^d$ be a sub-Gaussian vector with variance proxy $\sigma^2$. Then, 
\begin{align}
    \mathbb{P}\left(\max_{\mc\in\mathbb{R}^d}\left\langle \frac{\mc}{\|\mc\|_2}, \me\right\rangle\geq t\right)\leq \exp\left(-\frac{t^2}{8\sigma^2}+d\log 6\right).
\end{align}  
\end{lem}

\begin{lem}[Varshamov-Gilbert bound]~\label{lem:covering} There exists a subset $\{\gamma_1,\ldots,\gamma_N\}\subset\{0,1\}^d$ such that 
\begin{align}
    \rho_H(\gamma_i,\gamma_j):= \FnormSize{}{\gamma_i-\gamma_j}^2 \geq \frac{d}{4} \text{ for any } i\neq j\in[N],
\end{align}
for some $N\geq \exp(d/8).$
\end{lem}

For any two probability measures $\mathbb{P}$ and $\mathbb{Q}$, define the Kullback-Leibler divergence by
\begin{equation}\label{eq:KL}
D(\mathbb{P}|\mathbb{Q})=\int\left(\log {d \mathbb{P} \over d \mathbb{Q}}\right)d\mathbb{P} .
\end{equation}

\begin{lem}[Proposition 4.1 in \citet{gao2015rate}]\label{prop:minmax}
Let $(\Xi,\rho)$ be a metric space and $\{\mathbb{P}_\xi\colon\xi\in\Xi\}$ be a collection of probability measure. 
For any totally bounded $T\subset\Xi$, define the Kullback-Leibler diameter of $T$ by $d_{KL}(T) = \sup_{\xi,\xi'\in T}D(\mathbb{P}_{\xi}|\mathbb{P}_{\xi'})$. Then, for any $\epsilon>0$,
\begin{align}
    \inf_{\hat\xi}\sup_{\xi\in\Xi}\mathbb{P}_\xi\left\{\rho^2(\hat\xi,\xi)\geq \frac{\epsilon^2}{4}\right\}\geq 1-\frac{d_{KL}(T)+\log 2}{\log \tM(\epsilon,T,\rho)},
\end{align}
where $\tM(\epsilon,T,\rho)$ is the packing number of $T$ by $\epsilon$-radius balls with respect to the metric $\rho$.
\end{lem}

\begin{lem} [Proposition 17 in \citet{gao2016optimal}]\label{prop:kl} Let $\mathbb{P}_{\Theta}$ (respectively, $\mathbb{P}_{\Theta'}$) denote the probability distribution of Gaussian tensors with mean $\Theta$ (respectively, $\Theta'$) and entrywise i.i.d.\ noise from $N(0,\sigma^2)$. Then we have
\begin{align}
    D(\mathbb{P}_\Theta|\mathbb{P}_{\Theta'})\leq \frac{1}{2\sigma^2}\sum_{\omega\in[d]^m}(\Theta(\omega)-\Theta'(\omega))^2.
\end{align}
\end{lem}

\end{document}